\documentclass{tac}

\usepackage{amssymb}
\usepackage{amsfonts}
\usepackage{amsmath}
\usepackage{url}                

\usepackage[all,2cell]{xy}
\UseAllTwocells
\usepackage{mathtools}
\usepackage{colortbl}

\usepackage{verbatim} 

\input diagxy

\author {G.S.H. Cruttwell and Michael A.\ Shulman}

\thanks{The first author was supported by a PIMS Calgary postdoctoral
  fellowship, and the second author by a National Science Foundation
  postdoctoral fellowship during the writing of this paper.}

\address{Department of Computer Science, University of Calgary\\
  2500 University Dr. NW, Calgary, Alberta, T2N 1N4 Canada\\[5pt]
  Department of Mathematics, University of California, San Diego\\
  9500 Gilman Dr \#0112, La Jolla, CA, 92023-0112, USA\\
}

\title {A unified framework for generalized multicategories}

\copyrightyear{2010}


\eaddress{gscruttw@ucalgary.ca \CR mshulman@ucsd.edu}

\newtheoremrm{nonexample}{Non-Example}

\mathrmdef{Hom}
\mathrmdef{ob}
\mathrmdef{Id}
\mathbfdef{Set}
\mathbfdef{Cat}

\newcommand{\V}{\ensuremath{\mathbf{V}}}

\newcommand{\C}{\ensuremath{\mathbf{C}}}
\newcommand{\calV}{\ensuremath{\mathcal{V}}}
\newcommand{\calH}{\ensuremath{\mathcal{H}}}
\newcommand{\calK}{\ensuremath{\mathcal{K}}}
\newcommand{\calL}{\ensuremath{\mathcal{L}}}
\newcommand{\mb}[1]{{\mathbf{#1}}}

\newcommand{\qed}{\ \hfill \rule{1ex}{1ex}}
\newcommand{\iso}{\cong}
\newcommand{\maps}{\colon}

\newcommand{\bbtwo}{\ensuremath{\mathbf{2}}}
\newcommand{\eR}{\ensuremath{\overline{\mathbf{R}_+}}}

\newcommand{\mat}[1]{\mathbf{#1}\mbox{-}\mathcal{M}\mathit{at}}
\newcommand{\Mat}{\mbox{-}\mathbb{M}\mathsf{at}}
\newcommand{\MAT}[1]{\mathbf{#1}\mbox{-}\mathbb{M}\mathsf{at}}
\newcommand{\prof}[1]{\mathbf{#1}\mbox{-}\mathcal{P}\mathit{rof}}
\newcommand{\cat}[1]{\mathbf{#1}\mbox{-}\mathcal{C}\mathit{at}}
\newcommand{\icat}[1]{\mathcal{C}\mathit{at}(\mathbf{#1})}
\newcommand{\PROF}[1]{\mathbf{#1}\mbox{-}\mathbb{P}\mathsf{rof}}
\newcommand{\Span}[1]{\mathcal{S}\mathit{pan}(\mb{#1})}
\newcommand{\SPAN}[1]{\mathbb{S}\mathsf{pan}\!\left(\mb{#1}\right)}

\newcommand{\iPROF}[1]{\mathbb{P}\mathsf{rof}(\mb{#1})}
\newcommand{\DFib}{\mathbb{D}\mathsf{Fib}}
\newcommand{\MOD}{\mathbb{M}\mathsf{od}}
\newcommand{\Mod}{\ensuremath{\mathbb{M}\mathsf{od}}}
\newcommand{\KMod}{\mathbb{K}\mathsf{Mod}}
\newcommand{\nMod}{\mathsf{n}\mathbb{K}\mathsf{Mod}}

\newcommand{\dMod}{\mathsf{d}\mathbb{K}\mathsf{Mod}}
\newcommand{\Mon}{\mathcal{M}\mathit{on}}
\newcommand{\KMon}{\mathcal{K}\mathit{Mon}}
\newcommand{\nKMon}{\mathit{n}\mathcal{K}\mathit{Mon}}
\newcommand{\dKMon}{\mathit{d}\mathcal{K}\mathit{Mon}}

\newcommand{\pro}{\mbox{pro-}}
\newcommand{\Rel}{\mathbb{R}\mathsf{el}}
\newcommand{\twoCat}{2\mbox{-}\mathcal{C}\mathit{at}}
\newcommand{\SQ}{\mathbb{S}\mathsf{q}}

\newcommand{\vEquip}{\ensuremath{\mathit{v}\mathcal{E}\!\mathit{quip}}}
\newcommand{\vDbl}{\ensuremath{\mathit{v}\mathcal{D}\mathit{bl}}}
\newcommand{\vDblr}{\ensuremath{\mathit{v}\mathcal{D}\mathit{bl}_{r}}}
\newcommand{\vDbln}{\ensuremath{\mathit{v}\mathcal{D}\mathit{bl}_{n}}}

\newcommand{\op}{\mbox{\scriptsize op}}

\newcommand{\ten}{\otimes}

\newcommand{\Hkl}[2]{\mathbb{H}\text{-}\mathsf{Kl}(\mathbb{#1},#2)}
\newcommand{\HKl}{\mathbb{H}\text{-}\mathsf{Kl}}

\def\htopppp/#1/<#2>^#3_#4{\:%
\ifnum#2=0%
   \setwdth{#3}{#4}\deltax=\wdth \divide \deltax by \ul%
   \advance \deltax by \defaultmargin  \ratchet{\deltax}{100}%
\else \deltax #2%
\fi%
\xy\ar@{#1}|-@{|}^{#3}_{#4}(\deltax,0) \endxy%
\:}%
\def\htoppp/#1/<#2>^#3{\ifnextchar_{\htopppp/#1/<#2>^{#3}}{\htopppp/#1/<#2>^{#3}_{}}}%
\def\htopp/#1/<#2>{\ifnextchar^{\htoppp/#1/<#2>}{\htoppp/#1/<#2>^{}}}%
\def\htoop/#1/{\ifnextchar<{\htopp/#1/}{\htopp/#1/<0>}}%
\def\hto{\ifnextchar/{\htoop}{\htoop/>/}}%

\makeatletter
\let\c@equation\c@subsection
\makeatother
\numberwithin{equation}{section}

\renewcommand{\theenumi}{(\roman{enumi})}

\makeatletter                   
\let\your@state\state
\def\state#1{\gdef\currthmtype{#1}\your@state{#1}}
\let\your@staterm\staterm
\def\staterm#1{\gdef\currthmtype{#1}\your@staterm{#1}}
\def\currthmtype{}
\newif\ifhyperref
\@ifpackageloaded{hyperref}{\hyperreftrue}{\hyperreffalse}
\ifhyperref
\def\autoref#1{\ref*{label@name@#1}~\ref{#1}}
\else
\def\autoref#1{\ref{label@name@#1}~\ref{#1}}
\fi
\AtBeginDocument{
  \let\old@label\label
  \def\label#1{%
    {\let\your@currentlabel\@currentlabel%
      \edef\@currentlabel{\currthmtype}%
      \old@label{label@name@#1}}%
    \old@label{#1}}
}
\makeatother                    

\newenvironment{packeditem}{%
  \begin{itemize}%
    \setlength{\itemsep}{0pt}%
    \setlength{\topsep}{0pt}%
}{\end{itemize}}

\begin{document}

\maketitle
\begin{abstract}
  Notions of generalized multicategory have been defined in numerous
  contexts throughout the literature, and include such diverse
  examples as symmetric multicategories, globular operads, Lawvere
  theories, and topological spaces.  In each case, generalized
  multicategories are defined as the ``lax algebras'' or ``Kleisli
  monoids'' relative to a ``monad'' on a bicategory.
  However, the meanings of these words differ from author to author,
  as do the specific bicategories considered.  We propose a unified
  framework: by working with monads on double categories and related
  structures (rather than bicategories), one can define generalized
  multicategories in a way that unifies all previous examples, while
  at the same time simplifying and clarifying much of the theory.
\end{abstract}

\tableofcontents

\section{Introduction}
\label{sec:introduction}

A \emph{multicategory} is a generalization of a category, in which the
domain of a morphism, rather than being a single object, can be a
finite list of objects.  A proto\-typical example is the
multicategory $\mathbf{Vect}$ of vector spaces, in which a morphism
$(V_1,\dots,V_n)\to W$ is a multilinear map.  In fact, any monoidal
category gives a multicategory in a canonical way, where the morphisms
$(V_1,\dots,V_n)\to W$ are the ordinary morphisms $V_1\ten\dots\ten
V_n\to W$.  The multicategory $\mathbf{Vect}$ can be seen as arising
in this way, but it is also natural to
view its multicategory structure as more basic, with the tensor product then
characterized as a representing object for ``multimorphisms.''
This is also the case for many other multicategories; in fact,
monoidal categories can be identified with multicategories
satisfying a certain representability property
(see~\cite{HermidaRepresentable} and \S\ref{sec:vert-alg}).

In addition to providing an abstract formalization of the passage from
``multilinear map'' to ``tensor product,'' multicategories provide a
convenient way to present certain types of finitary algebraic theories
(specifically, \emph{strongly regular} finitary theories, whose axioms
involve no duplication, omission, or permutation of variables).
Namely, the objects of the multicategories are the \emph{sorts} of the
theory, and each morphism $(X_1,\dots,X_n)\to Y$ represents an
algebraic operation of the theory.  When
viewed in this light, multicategories (especially those with one
object, which correspond to one-sorted theories) are often called
\emph{non-symmetric operads} (see \cite{may:goils}).
The original definition of multicategories in
\cite{lambek:dedsys-ii} (see also \cite{lambek:multi-revis}) was also
along these lines (a framework for sequent calculus).  The two
viewpoints are related by the observation that when $A$ is a small
multicategory representing an algebraic theory, and $\mb{C}$ is a
large multicategory such as $\mb{Vect}$, a \emph{model} of the theory
$A$ in $\mb{C}$ is simply a \emph{functor} of multicategories
$A\to\mb{C}$.  This is a version of the \emph{functorial
  semantics} of~\cite{lawvere:functsem}.

Our concern in this paper is with \emph{generalized multicategories},
a well-known idea which generalizes the basic notion in two ways.
Firstly, one allows a
change of ``base context,'' thereby including both \emph{internal}
multicategories and \emph{enriched} multicategories.  Secondly,
and more interestingly, one allows the finite lists of objects serving
as the domains of morphisms to be replaced by ``something else.''
From the first point of view, this is desirable since there are many
other contexts in which one would like to analyze the relationship
between structures with coherence axioms (such as monoidal categories)
and structures with universal or ``representability'' properties.
From the second point of view, it is desirable since not all
algebraic theories are strongly regular.

For example, generalized multicategories include \emph{symmetric}
multicategories, in which the finite lists can be arbitrarily
permuted.  ``Representable'' symmetric multicategories correspond to
symmetric monoidal categories.
Enriched symmetric multicategories with one object can be identified
with the \emph{operads} of~\cite{may:goils,kelly:operads,km:oamm}.
These describe algebraic theories in whose axioms variables can be
permuted (but not duplicated or omitted).  In most applications of
operads (see~\cite{em:rma-infloop,brgmoer:operads} for some recent
ones), both symmetry and enrichment are essential.

An obvious variation of symmetric multicategories is \emph{braided}
multicategories.  If we
allow duplication and omission in addition to permutation of inputs,
we obtain (multi-sorted) \emph{Lawvere
  theories}~\cite{lawvere:functsem}; a slight modification also
produces the \emph{clubs} of~\cite{kelly:mv-funct-calc,kelly:clubs}.
There are also important generalizations to ``algebraic theories'' on
more complicated objects; for instance, the \emph{globular operads} of
\cite{batanin:monglob,leinster:higher-opds} describe a certain sort of
algebraic theory on globular sets that includes many notions of
weak $n$-category.

A very different route to generalized multicategories begins with the
observation of~\cite{RelationalAlgebras} that \emph{topological
  spaces} can be viewed as a type of generalized multicategory, when
finite lists of objects are replaced by ultrafilters, and morphisms
are replaced by a convergence relation.  Many other sorts of
topological structures, such as metric spaces,
closure spaces, uniform spaces, and approach spaces, can also be seen
as generalized multicategories;
see~\cite{LawvereMetric,CommonApproach,OneSetting}.

With so many different faces, it is not surprising that generalized
multicategories have been independently considered by many authors.
They were first studied in generality by \cite{burroni:t-cats},
but have also been considered by many other
authors, including \cite{leinster:higher-opds},
\cite{GeneralizedEnrichmentofCategories}, \cite{HermidaCoherent},
\cite{CommonApproach}, \cite{OneSetting}, \cite{RelationalAlgebras},
\cite{weber:operads}, \cite{bd:hda3}, \cite{cheng:opetopic}, and
\cite{ds:lmpo-conv}.  While
all these authors are clearly doing ``the same thing'' from an intuitive
standpoint, they work in different frameworks at different
levels of generality, making the formal definitions difficult to compare.
Moreover, all of these approaches share a
certain \emph{ad hoc} quality, which, given the naturalness and
importance of the notion, ought to be avoidable.

In each case, the authors observe that the ``something else''
serving as the domain of morphisms in a generalized multicategory
should be specified by some sort of \emph{monad}, invariably denoted
$T$.  For example, ordinary multicategories appear when $T$ is the
``free monoid'' monad, globular operads appear when $T$ is the
``free strict $\omega$-category'' monad, and topological spaces appear
when $T$ is the ultrafilter monad.  All the difficulties then center
around what sort of thing $T$ is a monad \emph{on}.

Leinster~\cite{GeneralizedEnrichmentofCategories,leinster:higher-opds}
takes it to be a \emph{cartesian monad} on an ordinary category
$\mathbf{C}$, i.e.\ $\mathbf{C}$ has pullbacks, $T$ preserves them,
and the naturality squares for its unit and multiplication are
pullback squares.  Burroni~\cite{burroni:t-cats}, whose approach is
basically the same, is able to deal with any monad on a category with pullbacks.
Hermida~\cite{HermidaCoherent} works with a cartesian 2-monad on a
suitable 2-category.  Barr and Clementino
et.~al.~\cite{RelationalAlgebras,CommonApproach,OneSetting} work with a monad on
\Set\ equipped with a ``lax extension'' to the bicategory of matrices
in some monoidal category.  Weber~\cite{weber:operads} works with a
``monoidal pseudo algebra'' for a 2-monad on a suitable 2-category.
Baez-Dolan~\cite{bd:hda3} and Cheng~\cite{cheng:opetopic} (see
also~\cite{fghw:gen-species}) use a monad on \Cat\ extended to the
bicategory of profunctors (although they consider only the ``free
symmetric strict monoidal category'' monad).

Inspecting these various definitions and looking for commonalities, we
observe that in all cases,
the monads involved naturally live on a \emph{bicategory}, be it a
bicategory of spans (Burroni, Leinster), two-sided fibrations
(Hermida), relations (Barr), matrices (Clementino et.~al., Weber), or profunctors
(Baez-Dolan, Cheng).  What causes problems is that the monads of interest are frequently
\emph{lax} (preserving composition only up to a noninvertible
transformation), but there is no obvious general notion of
lax monad on a bicategory, since there is no good 2-category
(or even tricategory) of bicategories that contains lax or oplax
functors.

Furthermore, merely knowing the bicategorical monad (however one
chooses to formalize this) is insufficient for the theory, and in
particular for the definition of functors and transformations between
generalized multicategories.  Leinster, Burroni, Weber, and Hermida
can avoid this problem because their bicategorical monads are induced
by monads on some underlying category or 2-category.  Others resolve
it by working with an \emph{extension} of a given monad on \Set\ or
\Cat\ to the bicategory of matrices or profunctors, rather than merely
the bicategorical monad itself.  However, the various definitions of
such extensions are tricky to compare and have an \emph{ad hoc}
flavor.

Our goal in this paper (and its sequels) is to give a common framework
which includes \emph{all} previous approaches to generalized
multicategories, and therefore provides a natural context in which to
compare them.  To do this, instead of considering monads on
bicategories, we instead consider monads on types of double categories.
This essentially solves both problems mentioned above: on the one hand
there is a perfectly good 2-category of double categories and lax
functors (allowing us to define monads on a double category), and on the
other hand the vertical arrows of the double categories (such as
morphisms in the cartesian category $\mathbf{C}$, functions in \Set, or
functors between categories) provide the missing data with which to
define functors and transformations of generalized multicategories.

The types of double categories we use are neither strict or pseudo
double categories, but instead an even weaker notion, for the following
reason.  An important intermediate step in the definition of generalized
multicategories is the \emph{horizontal Kleisli} construction of a monad
$T$, whose (horizontal) arrows $X\hto Y$ are arrows $X\hto TY$.  Without
strong assumptions on $T$, such arrows cannot be composed associatively,
and hence the horizontal Kleisli construction does not give a pseudo
double category or bicategory.  It does, however, give a weaker
structure, which we call a \emph{virtual double category}.

Intuitively, virtual double categories generalize pseudo double
categories in the same way that multicategories generalize monoidal
categories.  There is no longer a horizontal composition operation, but
we have cells of shapes such as the following:
\[
\bfig
\square/`>`>`{@{>}|-*@{|}}/<2000,375>[X_0`X_n`Y_0`Y_1;`f`g`q]
\morphism(0,375)/@{>}|-*@{|}/<500,0>[X_0`X_1;p_1]
\morphism(500,375)/@{>}|-*@{|}/<500,0>[X_1`X_2;p_2]
\morphism(1000,375)/@{>}|-*@{|}/<500,0>[X_2`\cdots;p_3]
\morphism(1500,375)/@{>}|-*@{|}/<500,0>[\cdots`X_n;p_n]
\place(950,225)[\twoar(0,-1)]
\place(1050,175)[\scriptstyle{\alpha}]
\efig
\]
We will give an explicit definition in \S\ref{sec:fcmulticats}.  Virtual
double categories have been studied by \cite{burroni:t-cats} under the
name of \emph{multicat\'egories} and by \cite{leinster:higher-opds}
under the name of \textbf{fc}-\emph{multicategories}, both of whom
additionally described a special case of the horizontal Kleisli
construction.  They are, in fact, the generalized multicategories
relative to the ``free category'' or ``free double category'' monad
(depending on whether one works with spans or profunctors).
In~\cite{dpp:paths-dbl} virtual double categories were called \emph{lax
  double categories}, but we believe that name belongs properly to lax
algebras for the 2-monad whose strict algebras are double categories.
(We will see in \autoref{eg:vdblcat} that \emph{oplax} double categories
in this ``2-monadically correct'' sense can be identified with a
restricted class of virtual double categories.)

Next, in \S\S\ref{sec:2cat-vdc}--\ref{sec:genmulti} we will show that
for any monad $T$ on a virtual double category $\mathbb{X}$, one can
define a notion which we call a \textbf{$T$-monoid}.  In fact, we will
construct an entire new virtual double category $\KMod(\mathbb{X},T)$
whose objects are $T$-monoids, by composing the ``horizontal Kleisli''
construction mentioned above with the ``monoids and bimodules''
construction $\MOD$ described in~\cite[\S5.3]{leinster:higher-opds},
which can be applied to any virtual double category.  Then in
\S\ref{sec:2cat-tmon} we will construct from $\KMod(\mathbb{X},T)$ a
2-category $\KMon(\mathbb{X},T)$ of $T$-monoids, $T$-monoid functors,
and transformations, generalizing the analogous definition
in~\cite[\S5.3]{leinster:higher-opds}.  This requires a notion of when a
virtual double category has \emph{units}, which we define in
\S\ref{sec:composites-units} along with the parallel notion of when it
has \emph{composites}.  (These definitions generalize those
of~\cite{HermidaRepresentable} and can also be found
in~\cite{dpp:paths-dbl}; they are also a particular case of the
``representability'' of~\cite{HermidaCoherent} and our
\S\ref{sec:vert-alg}.)

For particular $\mathbb{X}$ and $T$, the notion of $T$-monoid
specializes to several previous definitions of generalized
multicategories.  For example, if $\mathbb{X}$ consists of objects and
spans in a cartesian category \C\ and $T$ is induced from a monad on \C,
we recover the definitions of Leinster, Kelly, and Burroni.  And if
$\mathbb{X}$ consists of sets and matrices enriched over some monoidal
category \V\ and $T$ is a ``canonical extension'' of a taut set-monad to
$\mathbb{X}$, then we recover the definitions of Clementino et.\ al.

However, the other definitions of generalized multicategory cannot quite
be identified with $T$-monoids for any $T$, but rather with only a
restricted class of them.  For instance, if $\mathbb{X}$ consists of
categories and profunctors, and $T$ extends the ``free symmetric
monoidal category'' monad on \Cat\ (this is the situation of Baez-Dolan
and Cheng), then $T$-monoids are not quite the same as ordinary
symmetric multicategories.  Rather, a $T$-monoid for this $T$ consists
of a category $A$, a symmetric multicategory $M$, and a
bijective-on-objects functor from $A$ to the underlying ordinary
category of $M$.  There are two ways to restrict the class of such
$T$-monoids to obtain a notion equivalent to ordinary symmetric
multicategories: we can require $A$ to be a discrete category (so that
it is simply the set of objects of $M$), or we can require the functor
to also be fully faithful (so that $A$ is simply the underlying ordinary
category of $M$).  We call the first type of $T$-monoid
\emph{object-discrete} and the second type \emph{normalized}.

In order to achieve a full unification, therefore, we must give general
definitions of these classes of $T$-monoid and account for their
relationship.  It turns out that this requires additional structure on
our virtual double categories: we need to assume that horizontal arrows
can be ``restricted'' along vertical ones, in a sense made precise in
\S\ref{sec:pfrbi}.  \emph{Pseudo} double categories with this property
were called \emph{framed bicategories} in~\cite{FramedBicats}, where
they were also shown to be equivalent to the \emph{proarrow equipments}
of~\cite{wood:proarrows-i} (see also \cite{DomThesis}).  Accordingly, if
a virtual double category $\mathbb{X}$ has this property, as well as all
units, we call it a \emph{virtual equipment}.

Our first result in \S\ref{sec:normalization}, then, is that if $T$ is a
well-behaved monad on a virtual equipment, object-discrete and
normalized $T$-monoids are equivalent.  However, normalized $T$-monoids
are defined more generally than object-discrete ones, and moreover when
$T$ which are insufficiently well-behaved, it is the normalized
$T$-monoids which are of more interest.  Thus, we subsequently discard
the notion of object-discreteness.  (Hermida's generalized
multicategories also arise as normalized $T$-monoids, where $\mathbb{X}$
consists of discrete fibrations in a suitable 2-category \calK\ and $T$
is an extension of a suitable 2-monad on \calK.  Weber's definition is a
special case, since as given it really only makes sense for generalized
operads, for which normalization is automatic; see
\S\ref{sec:mon-psalg}.)  In Table~\ref{tab:egs} we summarize the
meanings of $T$-monoids and normalized $T$-monoids for a number of
monads $T$.

\begin{table}
  \centering
  \small
\noindent\begin{tabular}{p{3cm}p{2.5cm}p{3cm}p{3cm}p{3cm}}

\textbf{Monad} $T$ & \textbf{on} & $T$\textbf{-monoid} & \textbf{Normalized} \mbox{$T$-\textbf{Monoid}} & \textbf{Pseudo} \mbox{$T$ \textbf{-algebra}} \\
\rowcolor[gray]{0.8} $\Id$ & $\MAT{V}$ & \V-enriched \mbox{category} & Set & Set \\
$\Id$ & $\SPAN{C}$ & Internal category in \C & Object of \C\ & Object of \C\ \\
\rowcolor[gray]{0.8}$\Id$ & $\Rel$ & Ordered Set & Set & Set \\
$\Id$ & $\MAT{\eR}$ & Metric Space & Set & Set \\
\rowcolor[gray]{0.8}Powerset & $\Rel$ & Closure Space & $T_1$ Closure Space & Complete \mbox{Semilattice} \\
$\MOD(\text{powerset})$ & $\PROF{2}$ & Modular Closure Space  & Closure Space & Meet-Complete Preorder \\
\rowcolor[gray]{0.8}Ultrafilter & $\Rel$ & Topological Space & $T_1$ space & Compact \mbox{Hausdorff} space \\
$\MOD(\text{ultrafilter})$ & $\PROF{2}$ & \raggedright Modular Top.\ Space & Topological Space & Ordered Compact Hausdorff space \\
\rowcolor[gray]{0.8}Ultrafilter & $\MAT{\eR}$ & Approach space & ? & Compact \mbox{Hausdorff} space \\
Free monoid & $\SPAN{Set}$ & Multicategory & ? & Monoid \\
\rowcolor[gray]{0.8}$\MOD(\text{free monoid})$ & $\PROF{Set}$ & ``Enhanced'' multi\-category & Multicategory & Monoidal \mbox{category} \\
Free sym.\ strict mon.\ cat. & $\PROF{Set}$ & ``Enhanced'' sym.\ multicategory & Symmetric \mbox{multicategory} & Symmetric mon.\ cat. \\
\rowcolor[gray]{0.8}Free category & $\SPAN{Grph}$ & \raggedright Virtual double category & ? & Category \\
$\MOD(\text{free category})$ & $\iPROF{Grph}$ & ? & \raggedright Virtual double category & Pseudo double category \\
\rowcolor[gray]{0.8}\raggedright Free cat.\ w/ finite products & $\PROF{Set}$ & ? & Multi-sorted Lawvere theory & Cat.\ w/ finite products \\
\raggedright Free cat.\ w/ small products & $\PROF{Set}_{\mathrm{ls}}$ & ? & Monad on $\mb{Set}$ & Cat.\ w/ small products \\
\rowcolor[gray]{0.8}\raggedright Free presheaf $\mb{S}^{\op}\to \Set$ & $\SPAN{\mb{Set}^{\ob(\mb{S})}}$ & Functor $\mb{A}\to\mb{S}$ & Functor $\mb{A} \to \mb{S}$ w/ discrete fibers & Functor $\mb{S}^{\op} \to \mb{Set}$ \\
$\MOD(\text{free presheaf})$ & $\iPROF{\mb{Set}^{\ob(S)}}$ & ? & Functor $\mb{A} \to \mb{S}$ & Pseudofunctor $\mb{S}^{\op} \to \mb{Cat}$ \\
\rowcolor[gray]{0.8}\raggedright Free strict \mbox{$\omega$-category} & $\SPAN{Globset}$ & Multi-sorted globular operad & ? & Strict $\omega$-category\\
$\MOD(\text{free $\omega$-cat.})$ & $\iPROF{Globset}$ & ? & Multi-sorted globular operad & Monoidal \mbox{globular} cat.\\
\rowcolor[gray]{0.8}Free $M$-set ($M$ a monoid) & $\SPAN{Set}$ & $M$-graded \mbox{category} & ? & $M$-set\\
\end{tabular}
\caption{Examples of generalized multicategories.
  The boxes marked ``?''\ do not have any established name; in most
  cases they also do not seem very interesting.}
\label{tab:egs}
\end{table}

Now, what determines whether the ``right'' notion of generalized
multicategory is a plain $T$-monoid or a normalized one?  The obvious
thing distinguishing the situations of Leinster, Burroni, and Clementino
et.\ al.\ from those of Baez-Dolan, Cheng, and Hermida is that in the
former case, the objects of $\mathbb{X}$ are ``set-like,'' whereas in
the latter, they are ``category-like.''  However, some types of
generalized multicategory can be constructed starting from two different
monads on two different virtual equipments, one of which belongs to the
first group and the other to the second.

For example, observe that an ordinary (non-symmetric) multicategory has
an underlying ordinary category, containing the same objects but only
the morphisms $(V_1) \to W$ with unary source.  Thus, such a
multicategory can be defined in two ways: either as extra structure on
its \emph{set} of objects, or as extra structure on its underlying
\emph{category}.  In the second case, normalization is the requirement
that in the extra added structure, the multimorphisms with unary source
do no more than reproduce the originally given category.  Thus, ordinary
multicategories arise both as $T$-monoids the ``free monoid'' monad on
sets and spans, and as normalized $T$-monoids for the ``free monoidal
category'' monad on categories and profunctors.

Our second result in \S\ref{sec:normalization} is a generalization of
this relationship.  We observe that the virtual equipment of categories
and profunctors results from applying the ``monoids and modules''
construction $\Mod$ to the virtual equipment of sets and spans.  Thus,
we generalize this situation by showing that for any monad $T$ on a
virtual equipment, plain $T$-monoids can be identified with normalized
$\Mod(T)$-monoids.  That this is so in the examples can be seen by
inspection of Table~\ref{tab:egs}.  Moreover, it is sensible because
application of $\Mod$ takes ``set-like'' things to ``category-like''
things.

It follows that the notion of ``normalized $T$-monoid'' is actually
\emph{more general} than the notion of $T$-monoid, since arbitrary
$T$-monoids for some $T$ can be identified with the normalized
$S$-monoids for some $S$ (namely $S=\Mod(T)$), whereas normalized
$S$-monoids cannot always be identified with the arbitrary $T$-monoids
for any $T$.  (For instance, this is not the case when $S$ is the ``free
\emph{symmetric} monoidal category'' monad on categories and
profunctors.)  This motivates us to claim that the ``right'' notion of
generalized multicategory is \emph{a normalized $T$-monoid}, for some
monad $T$ on a virtual equipment.

Having reached this conclusion, we also take the opportunity to propose
a new naming system for generalized multicategories which we feel is
more convenient and descriptive.  Namely, if (pseudo) $T$-algebras are
called \emph{widgets}, then we propose to call normalized $T$-monoids
\emph{virtual widgets}.  The term ``virtual double category'' is of
course a special case of this: virtual double categories are the
normalized $T$-monoids for the monad $T$ on $\iPROF{Grph}$ whose
algebras are double categories.

Of course, ``virtual'' used in this way is a ``red herring''
adjective\footnote{The \emph{mathematical red herring principle} states
  that an object called a ``red herring'' need not, in general, be
  either red or a herring.} akin to ``pseudo'' and ``lax'', since a
virtual widget is not a widget.  The converse, however, is true: every
widget has an underlying virtual widget, so the terminology makes some
sense. For example, the observation above that every monoidal category
has an underlying multicategory is an instance of this fact.  Moreover,
it often happens that virtual widgets share many of the same properties
as widgets, and many theorems about widgets can easily be extended to
virtual widgets.  Thus, it is advantageous to use a terminology which
stresses the close connection between the two.  Another significant
advantage of ``virtual widget'' over ``$T$-multicategory'' is that
frequently one encounters monads $T$ for which $T$-algebras have a
common name, such as ``double category'' or ``symmetric monoidal
category,'' but $T$ itself has no name less cumbersome than ``the free
double category monad'' or ``the free symmetric monoidal category
monad.''  Thus, it makes more sense to name generalized multicategories
after the algebras for the monad than after the monad itself.

By the end of \S\ref{sec:normalization}, therefore, we have unified all
existing notions of generalized multicategory under the umbrella of
\emph{virtual $T$-algebras}, where $T$ is a monad on some virtual
equipment.  Since getting to this point already takes us over 40 pages,
we leave to future work most of the development of the theory and its
applications, along with more specific comparisons between existing
theories (see~\cite{cs:functoriality,cs:laxidem}).

However, we do spend some time in \S\ref{sec:vert-alg} on the topic of
\emph{representability}.  This is a central idea in the theory of
generalized multicategories, which states that any pseudo $T$-algebra
(or, in fact, any oplax $T$-algebra) has an underlying virtual
$T$-algebra.  Additionally, one can characterize the virtual
$T$-algebras which arise in this way by a ``representability'' property.
This can then be interpreted as an alternate \emph{definition} of pseudo
$T$-algebra which replaces ``coherent algebraic structure'' by a
``universal property,'' as advertised in~\cite{HermidaCoherent}.  In
addition to the identification of monoidal categories with
``representable'' multicategories, this also includes the fact that
compact Hausdorff spaces are $T_1$ spaces with additional properties,
and that fibrations over a category $\mb{S}$ are equivalent to
pseudofunctors $\mb{S}^{\op} \to \mb{Cat}$.  In~\cite{cs:laxidem} we
will extend more of the theory of representability
in~\cite{HermidaCoherent} to our general context.

Finally, the appendices are devoted to showing that all existing notions of
generalized multicategory are included in our framework.  In
Appendix~\ref{sec:composites-mod-hkl} we prepare the way by giving
sufficient conditions for our constructions on virtual double
categories to preserve composites, which is important since most
existing approaches use bicategories.  Then in
Appendix~\ref{sec:comparisons} we summarize how each existing theory
we are aware of fits into our context.

We have chosen to postpone these comparisons to the end, so that the
main body of the paper can present a unified picture of the subject,
in a way which is suitable also as an introduction for a reader
unfamiliar with any of the existing approaches.  It should be noted,
though, that we claim no originality for any of the examples or
applications, or the ideas of representability in
\S\ref{sec:vert-alg}.  Our goal is to show that all of these examples
fall into the same framework, and that this general framework allows for
a cleaner development of the theory.

\subsection{Acknowledgements}
\label{sec:acknowledgements}

The first author would like to thank Bob Par\'e for his suggestion to
consider ``double triples'', as well as helpful discussions with Maria
Manuel Clementino, Dirk Hofmann, and Walter Tholen.  The second author
would like to thank David Carchedi and the patrons of the $n$-Category
Caf\'e blog for several helpful conversations.  Both authors would like
to thank the editor and the referee for helpful suggestions.


\section{Virtual double categories}
\label{sec:fcmulticats}

The definition of virtual double category may be somewhat imposing, so
we begin with some motivation that will hopefully make it seem
inevitable.  We seek a framework which includes all sorts of
generalized multicategories.  Since \emph{categories} themselves are a
particular sort of generalized multicategory (relative to an identity
monad), our framework should in particular include all sorts of
generalized categories.  In particular, it should include both
categories \emph{enriched} in a monoidal category \V\ and categories
\emph{internal} to a category \C\ with pullbacks, so let us begin by
considering how to unify these two situations.

We start by recalling that both \V-enriched categories and \C-internal
categories are particular cases of \emph{monoids in a monoidal
  category}.  On the one hand, if \V\ is a cocomplete closed monoidal
category and $\mathcal{O}$ is a fixed set, then \V-enriched categories
with object set $\mathcal{O}$ can be identified with monoids in the
monoidal category of \emph{$\mathcal{O}$-graphs} in \V---i.e.\
$(\mathcal{O}\times\mathcal{O})$-indexed families of objects of \V, with
monoidal structure given by ``matrix multiplication.''  On the other
hand, if \C\ is a category with pullbacks and $O$ is an object of \C,
then \C-internal categories with object-of-objects $O$ can be identified
with monoids in the monoidal category of \emph{$O$-spans} in \C---i.e.\
diagrams of the form $O \leftarrow A \rightarrow O$, with monoidal
structure given by span composition.

Now, both of these examples share the same defect: they require us to
fix the objects (the set $\mathcal{O}$ or object $O$).  In particular,
the \emph{morphisms} of monoids in these monoidal categories are
functors which are the identity on objects.  It is well-known that one
can eliminate this fixing of objects by combining all the monoidal
categories of graphs and spans, respectively, into a \emph{bicategory}.
(In essense, this observation dates all the way back to~\cite{benabou:bicats}.)
In the first case
the relevant bicategory $\mat{V}$ consists of \emph{\V-matrices}: its
objects are sets, its arrows from $X$ to $Y$ are $(X\times Y)$-matrices
of objects in \V, and its composition is by ``matrix multiplication.''
In the second case the relevant bicategory $\Span{C}$ consists of
\emph{\C-spans}: its objects are objects of \C, its arrows from $X$
to $Y$ are spans $X\leftarrow A\rightarrow Y$ in \C, and its
composition is by pullback.  It is easy to define monoids in a
bicategory to generalize monoids in a monoidal
category\footnote{Monoids in a bicategory are usually called
  \emph{monads}.  However, we avoid that term for these sorts of
  monoids for two reasons.  Firstly, the morphisms of monoids we are
  interested in are not the same as the usual morphisms of monads
  (although they are related; see \cite[\S2.3--2.4]{ls:ftm2}).
  Secondly, there is potential for
  confusion with the monads \emph{on} bicategories and related
  structures which play an essential role in the theory we present.}.

However, we still have the problem of functors.  There is no way to
define morphisms between monoids in a bicategory so as to recapture
the correct notions of enriched and internal functors in $\mat{V}$ and
$\Span{C}$.  But we can solve this problem if instead of bicategories
we use \emph{(pseudo) double categories}, which come with objects, two
different kinds of arrow called ``horizontal'' and ``vertical,'' and
2-cells in the form of a square:
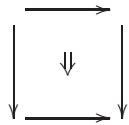
\[\bfig\scalefactor{.7}\square[\;`\;`\;`\;;```] \place(250,300)[\twoar(0,-1)]\efig\]
Both $\mat{V}$ and $\Span{V}$ naturally enlarge to pseudo double
categories, interpreting their existing arrows and composition as
horizontal and adding new vertical arrows.   For $\mat{V}$ the new
vertical arrows are functions between sets, while for $\Span{C}$ the
new vertical arrows are morphisms in \C.  We can now define monoids in
a double category (relative to the horizontal structure) and morphisms
between such monoids (making use of the vertical arrows) so as to
recapture the correct notion of functor in both cases (see
\autoref{def:mod}).

The final generalization from pseudo double categories to virtual double
categories is more difficult to motivate at the moment, but as remarked
in the introduction, we will find it essential in \S\ref{sec:genmulti}.
Conceptually (and, in fact, formally), a virtual double category is
related to a pseudo double category in the same way that a multicategory
is related to a monoidal category.  Thus, just as one can define monoids
in any multicategory, one can likewise do so in any virtual double
category.

\begin{definition}
A \textbf{virtual double category} $\mathbb{X}$ consists of the
following data.
	\begin{itemize}
	\item A category $\mb{X}$ (the objects and vertical arrows), with the arrows written vertically:
		\[
			\bfig
			\morphism<0,-250>[X`Y;f]
			\efig
		\]
	\item For any two objects $X, Y \in \mb{X}$, a class of
          horizontal arrows, written horizontally with a slash
          through the arrow:
		\[  \bfig \morphism/@{>}|-*@{|}/[X`Y;p] \efig \]
	\item Cells, with vertical source and target, and horizontal multi-source and target, written as follows:
          \begin{equation}\label{eq:cell}
			\bfig
			\square/`>`>`{@{>}|-*@{|}}/<2000,375>[X_0`X_n`Y_0`Y_1;`f`g`q]
			\morphism(0,375)/@{>}|-*@{|}/<500,0>[X_0`X_1;p_1]
			\morphism(500,375)/@{>}|-*@{|}/<500,0>[X_1`X_2;p_2]
			\morphism(1000,375)/@{>}|-*@{|}/<500,0>[X_2`\cdots;p_3]
			\morphism(1500,375)/@{>}|-*@{|}/<500,0>[\cdots`X_n;p_n]
			\place(950,225)[\twoar(0,-1)]
			\place(1050,175)[\scriptstyle{\alpha}]
			\efig
          \end{equation}
	Note that this includes cells with source of length 0, in which case
  we must have $X_0=X_n$; such cells are visually represented as follows:
		\[
			\bfig
			\morphism<-250,-375>[X`Y_0;f]
			\morphism<250,-375>[X`Y_1;g]
			\morphism(-250,-375)|b|/@{>}|-*@{|}/<500,0>[Y_0`Y_1;q]
			\place(-30,-175)[\twoar(0,-1)]
			\place(50,-225)[\scriptstyle{\alpha}]
			\efig
		\]

	\item For the following configuration of cells,
		\[
			\bfig
			\square/`>`>`{@{>}|-*@{|}}/<3000,250>[Y_0`Y_m`Z_0`Z_1;`g_0`g_1`r]
			\square(0,250)/-->`>`>`{@{>}|-*@{|}}/<750,250>[X_0`X_{n_1}`Y_0`Y_1;p_{11} \ldots p_{1n_1}`f_0`f_1`q_1]
			\square(750,250)/-->`>`>`{@{>}|-*@{|}}/<750,250>[X_{n_1}`X_{n_2}`Y_1`Y_2;p_{21} \ldots p_{2n_1}``f_2`q_2]
			\morphism(1500,500)/-->/<750,0>[X_{n_2}`\cdots;\cdots]
			\morphism(2250,500)/-->/<750,0>[\cdots`X_{n_m};\cdots]
			\morphism(1500,250)|b|/-->/<1500,0>[Y_2`Y_m;q_3 \cdots q_m]
			\morphism(3000,500)|r|<0,-250>[X_{n_m}`Y_m;f_m]
			
			\place(400,375)[\scriptstyle{\alpha_1}]
			\place(300,425)[\twoar(0,-1)]
			\place(1150,375)[\scriptstyle{\alpha_2}]
			\place(1050,425)[\twoar(0,-1)]
			
			\place(2100,425)[\twoar(0,-1)]
			\place(2275,375)[\scriptstyle{\alpha_3} \cdots \scriptstyle{\alpha_m}]
			\place(1425,175)[\twoar(0,-1)]
			\place(1525,125)[\scriptstyle{\beta}]
			
			\efig
		\]
	a composite cell
		\[
			\bfig
			\square/`>`>`{@{>}|-*@{|}}/<3000,250>[Y_0`Y_m`Z_0`Z_1;`g_0`g_1`r]
			\square(0,250)/-->`>``/<750,250>[X_0`X_{n_1}`Y_0`;p_{11} \ldots p_{1n_1}`f_0``]
			\square(750,250)/-->```/<750,250>[X_{n_1}`X_{n_2}``;p_{21} \ldots p_{2n_1}```]
			\morphism(1500,500)/-->/<750,0>[X_{n_2}`\cdots;\cdots]
			\morphism(2250,500)/-->/<750,0>[\cdots`X_{n_m};\cdots]
			\morphism(3000,500)|r|<0,-250>[X_{n_m}`Y_m;f_m]
			
			\place(1750,250)[\scriptstyle{\beta (\alpha_m \boxdot \cdots \boxdot \alpha_1)}]
			\place(1450,300)[\twoar(0,-1)]

			\efig
		\]
			
	\item For each horizontal arrow $p$, an identity cell
	
		\[
			\bfig
			\node A(0,250)[X]
			\node B(500,250)[Y]
			\node C(0,0)[X]
			\node D(500,0)[Y]
			\arrow/@{>}|-*@{|}/[A`B;p]
			\arrow|b|/@{>}|-*@{|}/[C`D;p]
			\arrow/=/[A`C;]
			\arrow/=/[B`D;]
			\place(200,175)[\twoar(0,-1)]
			\place(300,125)[\scriptstyle{1_p}]
			\efig
		\]
	\item Associativity and identity axioms for cell composition.
          The associativity axiom states that
          \[\Big(\gamma(\beta_m \boxdot \dots \boxdot \beta_1)\Big)
          (\alpha_{m k_m}\boxdot \dots\boxdot \alpha_{1 1}) =
          \gamma \Big(
          \beta_m(\alpha_{m k_m}\boxdot\dots\boxdot \alpha_{m 1})
          \boxdot \dots\boxdot
          \beta_1(\alpha_{1 k_1}\boxdot\dots\boxdot \alpha_{1 1})
          \Big)
          \]
          while the identity axioms state that
          \[\alpha(1_{p_1} \boxdot\dots\boxdot 1_{p_n}) = \alpha
          \qquad\text{and}\qquad
          1_q(\alpha_1) = \alpha_1
          \]
          whenever these equations make sense.
	\end{itemize}
\end{definition}

\begin{remark}
  As mentioned in the introduction, virtual double categories have
  also been called \emph{fc-multicategories} by Leinster~\cite{leinster:higher-opds} and
  \emph{multicat\'egories} by Burroni~\cite{burroni:t-cats}.  Our terminology is chosen to
  emphasize their close relationship with double categories, and to
  fit into the general naming scheme of \S\ref{sec:vert-alg}.
\end{remark}

\begin{remark}
  In much of the double-category literature, it is common for the
  ``slashed'' arrows (spans, profunctors, etc.) to be the
  \emph{vertical} arrows.  We have chosen the opposite convention
  purely for economy of space: the cells in a virtual double category
  fit more conveniently on a page when their multi-source is drawn
  horizontally.
\end{remark}

\begin{examples}
  As suggested by the discussion at the beginning of this section,
  monoidal categories, bicategories, 2-categories, multicategories, and
  pseudo double categories can each be regarded as examples of virtual
  double categories, by trivializing the vertical or horizontal
  structure in various ways; see
  \cite[p.~4]{GeneralizedEnrichmentofCategories}
  or~\cite[\S5.1]{leinster:higher-opds} for details.
\end{examples}

We now present the two virtual double categories that will serve as
initial inputs for most our examples: spans and matrices.
(Both are also described in~\cite[Ch.~5]{leinster:higher-opds}.)
For consistency, we name all of our virtual double categories by their
horizontal arrows, rather than their vertical arrows or objects.

\begin{example}\label{eg:mat}
  Let $(\mb{V}, \otimes, I)$ be a monoidal category.  The virtual double
  category $\MAT{V}$ is defined as follows: its objects are sets, its
  vertical arrows are functions, its horizontal arrows $X \hto^p Y$ are
  families $\{p(y,x)\}_{x\in X, y\in Y}$ of objects of \V\ (i.e.\
  $(X\times Y)$-matrices), and a cell of the form~\eqref{eq:cell}
  consists of a family of \V-arrows
  \[ p_1(x_1, x_0) \otimes p_2(x_2, x_1) \otimes \cdots \otimes p_n(x_{n}, x_{n-1}) \to^{\alpha} q(fx_n, gx_0), \]
  one for each tuple $(x_0,\dots,x_n) \in X_0\times\dots\times X_n$.
  When $n=0$, of course, the $n$-ary tensor product on the left is to be
  interpreted as the unit object of \V.

  In particular, if $\mb{V}$ is the 2-element chain $0\le
  1$, with $\otimes$ given by $\wedge$, then the horizontal arrows of
  $\MAT{V}$ are relations.  In this case we denote $\MAT{V}$ by $\Rel$.
\end{example}

It is well-known that $\mb{V}$-matrices form a bicategory (and, in
fact, a pseudo double category) as long as $\mb{V}$ has coproducts
preserved by $\otimes$.  However, if we merely want a virtual double
category, we see that this requirement is unnecessary.  (In fact, \V\
could be merely a multicategory itself.)

\begin{example}
For a category $\mb{C}$ with pullbacks, the virtual double category
$\SPAN{C}$ is defined as follows: its objects and vertical arrows are
those of \C, its horizontal arrows $X\hto^p Y$ are spans $X\toleft P \to
Y$, and a cell of the form~\eqref{eq:cell} is a morphism of spans
\[ P_1 \times_{X_1} P_2 \times_{X_2} \cdots \times_{X_{n-1}} P_n \to^{\alpha} Q \]
lying over $f$ and $g$.
When $n=0$, the $n$-ary span composite in the domain is to be
interpreted as the identity span $X_0 \toleft X_0 \to X_0$.
\end{example}

Note that in this case, we \emph{do} need to require that $\mb{C}$ have
pullbacks.  If
$\mb{C}$ does not have pullbacks, a more natural setting would be to
consider $\SPAN{C}$ as a ``co-virtual double category'', in which the
horizontal \emph{target} of a cell is a string of horizontal arrows.
However, co-virtual double categories do not provide the structure
necessary to define generalized multicategories.

We now recall the construction of \emph{monoids and modules} in a
virtual double category from \cite[\S5.3]{leinster:higher-opds}.

\begin{definition}\label{def:mod}
Let $\mathbb{X}$ be a virtual double category.  The virtual double
category $\MOD(\mathbb{X})$ has the following components:
	\begin{itemize}
	\item The objects (\textbf{monoids}) consist of four parts
          $(X_0, X, \bar{x}, \hat{x})$: an object $X_0$ of $\mathbb{X}$,
          a horizontal endo-arrow  $X_0 \hto^{X} X_0$ in $\mathbb{X}$,
          and multiplication and unit cells
	\[
		\bfig 
		\node a(0,250)[X_0]
		\node b(500,250)[X_0]
		\node c(1000,250)[X_0]
		\node d(0,0)[X_0]
		\node e(1000,0)[X_0]
		
		\arrow/@{>}|-*@{|}/[a`b;X]
		\arrow/@{>}|-*@{|}/[b`c;X]
		\arrow/=/[a`d;]
		\arrow/=/[c`e;]
		\arrow|b|/@{>}|-*@{|}/[d`e;X]
		
		\place(475, 175)[\twoar(0,-1)]
		\place(550,125)[\scriptstyle{\bar{x}}]
		
		\place(1500,100)[\mbox{and}]

		\morphism(2000,300)/=/<-250,-375>[X_0`X_0;]
		\morphism(2000,300)/=/<250,-375>[X_0`X_0;]
		\morphism(1750,-75)|b|/@{>}|-*@{|}/<500,0>[X_0`X_0;X]
		\place(1975,100)[\twoar(0,-1)]
		\place(2050,50)[\scriptstyle{\hat{x}}]

		\efig
	\]
	satisfying associativity and identity axioms.
	\item The vertical arrows (\textbf{monoid homomorphisms})
          consist of two parts $(f_0,f)$: a vertical arrow $X_0
          \to^{f_0} Y_0$ in $\mathbb{X}$ and a cell in $\mathbb{X}$:
	\[
		\bfig
		\node a(0,250)[X_0]
		\node b(500,250)[X_0]
		\node c(0,0)[Y_0]
		\node d(500,0)[Y_0]
		\arrow/@{>}|-*@{|}/[a`b;X]
		\arrow|b|/@{>}|-*@{|}/[c`d;Y]
		\arrow[a`c;f_0]
		\arrow|r|[b`d;f_0]
		\place(225,175)[\twoar(0,-1)]
		\place(300,125)[\scriptstyle{f}]
		\efig
	\]
	which is compatible with the multiplication and units of $X$ and $Y$.
	\item The horizontal arrows (\textbf{modules}) consist of three
          parts $(p,\bar{p}_r, \bar{p}_l)$: a horizontal arrow $X_0
          \hto^p Y_0$ in $\mathbb{X}$ and two cells in $\mathbb{X}$:
	\[
		\bfig 
		\node a(0,250)[X_0]
		\node b(500,250)[X_0]
		\node c(1000,250)[Y_0]
		\node d(0,0)[X_0]
		\node e(1000,0)[Y_0]

		\arrow/@{>}|-*@{|}/[a`b;X]
		\arrow/@{>}|-*@{|}/[b`c;p]
		\arrow/=/[a`d;]
		\arrow/=/[c`e;]
		\arrow|b|/@{>}|-*@{|}/[d`e;p]
		\place(450,175)[\twoar(0,-1)]
		\place(550,125)[\scriptstyle\bar{p}_r]
		
		\place(1500,100)[\mbox{and}]
		
		\node a(2000,250)[X_0]
		\node b(2500,250)[Y_0]
		\node c(3000,250)[Y_0]
		\node d(2000,0)[X_0]
		\node e(3000,0)[Y_0]

		\arrow/@{>}|-*@{|}/[a`b;p]
		\arrow/@{>}|-*@{|}/[b`c;Y]
		\arrow/=/[a`d;]
		\arrow/=/[c`e;]
		\arrow|b|/@{>}|-*@{|}/[d`e;p]
		\place(2450,175)[\twoar(0,-1)]
		\place(2550,125)[\scriptstyle\bar{p}_l]
		
		\efig
	\]
	which are compatible with the multiplication and units of $X$ and $Y$.
	\item The cells are cells in $\mathbb{X}$:
		\[
			\bfig
			\square/`>`>`/<2000,400>[(X^0)_0`(X^n)_0`(Y^0)_0`(Y^1)_0;`f_0`g_0`]
			\morphism|b|/@{>}|-*@{|}/<2000,0>[(Y^0)_0`(Y^1)_0;q]
			\morphism(0,400)/@{>}|-*@{|}/<500,0>[(X^0)_0`(X^1)_0;p_1]
			\morphism(500,400)/@{>}|-*@{|}/<500,0>[(X^1)_0`(X^2)_0;p_2]
			\morphism(1000,400)/@{>}|-*@{|}/<500,0>[(X^2)_0`\cdots;p_3]
			\morphism(1500,400)/@{>}|-*@{|}/<500,0>[\cdots`(X^n)_0;p_n]
			\place(950,250)[\twoar(0,-1)]
			\place(1050,200)[\scriptstyle{\alpha}]
			\efig
		\]
	which are compatible with the left and right actions of the
        horizontal cells.
	\end{itemize}
\end{definition}

Note that, as observed in~\cite[\S5.3]{leinster:higher-opds}, we can
define $\MOD(\mathbb{X})$ without requiring \emph{any} hypotheses on the
virtual double category $\mathbb{X}$, in contrast to the situation for
bicategories or pseudo double categories.

\begin{example}\label{eg:prof}
We denote the virtual double category $\MOD(\MAT{V})$ by $\PROF{V}$; its
objects are $\V$-enriched categories, its vertical arrows are
$\V$-functors, its horizontal arrows are $\V$-profunctors, and its cells
are a generalization of the ``forms'' of~\cite{ds:monbi-hopfagbd}
(including, as a special case, natural transformations between
profunctors).
When \V\ is closed (hence enriched over itself) and symmetric,
\V-profunctors $C\hto^H D$ can be identified with \V-functors
$D^{op}\times C \to \V$.
\end{example}

Again, note that because we are working with virtual double
categories, we do not require that $\mb{V}$ have any colimits (in fact,
$\mb{V}$ could be merely a multicategory).

\begin{example}\label{eg:iprof}
Let $\mb{C}$ be a category with pullbacks.  We denote the virtual double
category $\MOD(\SPAN{C})$ by
$\iPROF{C}$; it consists of \emph{internal} categories,
functors, profunctors, and transformations in $\mb{C}$.
\end{example}

Note that $\MAT{Set}\cong \SPAN{Set}$ and thus $\PROF{Set} \cong \iPROF{Set}$.


\section{Monads on a virtual double category}
\label{sec:2cat-vdc}

We claimed in \S\ref{sec:introduction} that the ``inputs'' of a
generalized multicategory are parametrized by a monad.  Why should
this be so?  Suppose that we have an operation $T$ which, given a set
(or object) of objects $X$, produces a set (or object) $TX$ intended
to parametrize such inputs.  For ``ordinary'' multicategories, $TX$
will be the set of finite lists of elements of $X$.

Now, from the perspective of the previous section, the data of a
category includes an object $A_0$ and a horizontal arrow
$A_0\hto^A A_0$ in some virtual double category.  For example, if we work
in $\MAT{Set}$, then $A$ is a matrix consisting of the hom-sets
$A(x,y)$ for every $x,y\in A_0$.  Now, instead, we want to have
hom-sets $A(\mathfrak{x},y)$ whose domain $\mathfrak{x}$ is an element
of $TA_0$.  Thus, it makes sense to consider a horizontal arrow
$A_0\hto^A TA_0$ as part of the data of a $T$-multicategory.
(We use $A_0\hto TA_0$ rather than $TA_0\hto A_0$,
for consistency with Examples \ref{eg:prof} and \ref{eg:iprof}: the codomain of the
horizontal arrow datum of a monoid specifies the \emph{domains} of the
arrows \emph{in} that monoid.)

However, we now need to specify the units and composition of our
generalized multicategory.  The unit should be a cell into $A$ with
0-length domain, but its source and target vertical arrows can no
longer both be identities because $A_0 \neq TA_0$.  In an ordinary
multicategory, the identities are morphisms $(x)\to^{1_x} x$ whose
domain is a singleton list; in terms of $\MAT{Set}$ this can be
described by a cell
\[\bfig\scalefactor{.7}
\Atriangle/=`->`->/[A_0`A_0`TA_0;`\eta`A]
\place(500,300)[\twoar(0,-1)]
\efig
\]
where $A_0\to^{\eta} TA_0$ is the inclusion of singleton lists.

Regarding composition, in an ordinary multicategory we can compose a
morphism $(y_1,\dots,y_n)\to^g z$ not with a single morphism, but with
a \emph{list} of morphisms $(f_1,\dots, f_n)$ where $(x_{i1},\dots,
x_{ik_i}) \to^{f_i} y_i$.  In terms of $\MAT{Set}$ this represents the
fact that we cannot ask for a multiplication cell with domain $\hto^A
\hto^A$, since the domain of $A$ does not match its codomain, but
instead we can consider a cell with domain $\hto^A \hto^{TA}$, where
we extend $T$ to act on $\Set$-matrices in the obvious way.  Now,
however, the codomain of $TA$ is $T^2A_0$; in order to have a cell
with codomain $A$ we need to ``remove parentheses'' from the resulting
domain $((x_{11},\dots,x_{1k_1}),\dots,(x_{n1},\dots,x_{nk_n}))$ to
obtain a single list.  Thus the composition should be a cell
\[\bfig\scalefactor{.8}\square/`=`->`@{>}|-*@{|}/<1200,500>[A_0`T^2A_0`A_0`TA_0;``\mu`A]
\morphism(0,500)/@{>}|-*@{|}/<600,0>[A_0`TA_0;A]
\morphism(600,500)/@{>}|-*@{|}/<600,0>[TA_0`T^2A_0;TA]
\place(600,300)[\twoar(0,-1)]
\efig
\]
where $\mu$ is the ``remove parentheses'' function.  Of course, these
functions $\eta$ and $\mu$ are the structure maps of the ``free
monoid'' monad on $\Set$.  Thus we see that in order to define
ordinary multicategories, what we require is an ``extension'' of this
monad to $\MAT{Set}$.

In order to have a good notion of a \emph{monad on a virtual double
  category}, we need at least a 2-category of virtual double categories.
Since virtual double categories are themselves a special case of
generalized multicategories, it suffices to observe that generalized
multicategories of any sort form a 2-category.  However, since we have
not yet defined generalized multicategories in our context, at this
point in our exposition it is appropriate to give an explicit
description of the 2-category \vDbl.

Of course, the objects of \vDbl\ are virtual double categories, and its
1-morphisms (called \emph{functors} of virtual double categories) are
the obvious structure-preserving maps: functions from the objects,
vertical and horizontal arrows, and cells of the domain to those of the
codomain, preserving all types of source, target, identities, and
composition.  However, the definition of 2-morphisms in \vDbl\ is
slightly less obvious.

\begin{definition}\label{def:transforms}
Given functors $\mathbb{X} \two^{F}_{G} \mathbb{Y}$ of virtual double categories, a
\textbf{transformation} $F \to^{\theta} G$ consists of the following data.
	\begin{itemize}
	\item For each object $X$ in $\mathbb{X}$, a vertical arrow
    $FX \to^{\theta_X} GX$, which form the components of a natural
    transformation between the vertical parts of $F$ and $G$.
	\item For each horizontal arrow $X \hto^p Y$ in $\mathbb{X}$, a cell in $\mathbb{Y}$:
	\begin{equation}\label{eq:transf-2cell}
		\bfig
		\node a(0,250)[FX]
		\node b(500,250)[FY]
		\node c(0,0)[GX]
		\node d(500,0)[GY]
		\arrow/@{>}|-*@{|}/[a`b;Fp]
		\arrow|b|/@{>}|-*@{|}/[c`d;Gp]
		\arrow[a`c;\theta_X]
		\arrow|r|[b`d;\theta_Y]
		\place(215,175)[\twoar(0,-1)]
		\place(300,125)[\scriptstyle{\theta_p}]
		\efig
	\end{equation}
	\item An axiom asserting that $\theta$ is ``cell-natural,'' meaning that
          \[
          \theta_q (F\alpha) = (G\alpha) (\theta_{p_1}\boxdot \dots\boxdot \theta_{p_n})
          \]
          whenever this makes sense.
	\end{itemize}
  Virtual double categories, functors, and transformations form a 2-category denoted \vDbl.
\end{definition}

\begin{example}
  Any lax monoidal functor $\mb{V} \to^N \mb{W}$ induces functors
  $\MAT{V} \to^{N_*} \MAT{W}$ and $\PROF{V} \to^{N_*} \PROF{W}$ in an
  evident way.
  Moreover, any monoidal natural transformation $N \to^{\psi} M$ between
  lax monoidal functors induces transformations $N_* \to^{\psi_*} M_*$
  in both cases.
  In this way $\MAT{(-)}$ and $\PROF{(-)}$ become 2-functors from the
  2-category of monoidal categories to \vDbl.
\end{example}

\begin{example}
  Similarly, any pullback-preserving functor $\mb{C} \to^N \mb{D}$
  between categories with pullbacks induces functors $\SPAN{C}\to^{N_*}
  \SPAN{D}$ and $\iPROF{C}\to^{N_*} \iPROF{D}$, and any natural
  transformation $N \to^{\psi} M$ between such functors induces
  transformations $N_* \to^{\psi_*} M_*$, thereby making $\SPAN{-}$ and
  $\iPROF{-}$ into 2-functors as well.
\end{example}

\begin{example}
  When restricted to bicategories or pseudo double categories,
  functors of virtual double categories are equivalent to the usual
  notions of \emph{lax} functor.
  
  This is a special case of a general fact; see \autoref{thm:oplaxalg-funct}.
  When transformations of virtual double categories are similarly
  restricted, they become \emph{icons} in the sense of \cite{lack:icons}
  (for bicategories) and \emph{vertical transformations} (for pseudo
  double categories).
\end{example}

By a \emph{monad on a virtual double category} $\mathbb{X}$, we will
mean a monad in the 2-category $\vDbl$.  Thus, it consists of a
functor $T\maps \mathbb{X}\to \mathbb{X}$ and transformations
$\eta\maps \Id\to T$ and $\mu\maps TT\to T$ satisfying the usual axioms.
We now give the examples of such monads that we will be interested in.

\begin{example}\label{example:cartesian}
  Since 2-functors preserve monads, any pullback-preserving monad on a
  category \C\ with pullbacks induces monads on $\SPAN{C}$ and
  $\iPROF{C}$.
  Examples of such monads include the following.
	\begin{packeditem}
	\item The ``free monoid'' monad on $\mb{Set}$ (or, more
          generally, on any countably lextensive category).
	\item The ``free $M$-set'' monad $(M \times -)$ on $\mb{Set}$, for any monoid
          $M$ (or more generally, for any monoid object in a category
          with finite limits).
        \item The monad $(-)+1$ on any lextensive category.
	\item The ``free category'' monad on the category of directed graphs.
	\item The ``free strict $\omega$-category'' monad on the category of globular sets.
	\end{packeditem}
  Many more examples can be found in
  \cite[pp.~103--107]{leinster:higher-opds}; see also
  \S\ref{section:cartesian}.  By the argument above, each of these monads extends to a
  monad on a virtual double category of spans.
\end{example}

The assignments $\V\mapsto \MAT{V}$ and $\V\mapsto\PROF{V}$ are also
2-functorial, but the monads we obtain in this way from monads on
monoidal categories are not usually interesting for defining
multicategories.  However, there are some general ways to construct
monads on virtual double categories of matrices, at least when \V\ is
a preorder.  The following is due to \cite{SealCanonical}, which in
turn expands on \cite{OneSetting}.

\begin{example}\label{example:tv}
  By a \textbf{quantale} we mean a closed symmetric monoidal complete
  lattice.  A quantale is \textbf{completely distributive} if for any
  $b\in \V$ we have $b = \bigvee\{a \mid a\prec b\}$, where $a\prec b$
  means that whenever $b \le \bigvee S$ then there is an $s\in S$ with
  $a\le s$.  (If in this definition $S$ is required to be
  \emph{directed}, we obtain the weaker notion of a \emph{continuous
    lattice}.)  For us, the two most important completely
  distributive quantales are the following.
  \begin{packeditem}
  \item The two-element chain $\bbtwo = (0\le 1)$.
  \item The extended nonnegative reals $\eR=[0,\infty]$ with the
    reverse of the usual ordering and $\otimes = +$.
  \end{packeditem}

  A functor said to be \textbf{taut} if it preserves pullbacks of
  monomorphisms (and therefore also preserves monomorphisms).  A monad
  is \textbf{taut} if its functor part is taut, and moreover the
  naturality squares of $\eta$ and $\mu$ for any \emph{monomorphism}
  are pullbacks.
  Some important taut monads on \Set\ are the identity monad, the
  powerset monad (whose algebras are complete lattices), the filter
  monad, and the ultrafilter monad (whose algebras are compact Hausdorff
  spaces).

  Now let \V\ be a completely distributive quantale and $T$ a taut
  monad on \Set.  For a \V-matrix $X \hto^p Y$ and elements
  $\mathcal{F}\in TX$ and $\mathcal{G}\in TY$, define
  \[
  Tp({\cal G}, {\cal F}) =
  \bigvee \Bigl\{ v \in \mb{V} \Bigm|
  \forall B \subseteq Y: \bigl({\cal G} \in TB \Rightarrow {\cal F} \in T(p_v[B])\bigr) \Bigr\},
  \]
  where
  \[ p_v[B] = \Bigl\{ x \in X \Bigm| \exists y \in B : v \leq p(y,x) \Bigl\}.
  \]
  It is proven in~\cite{SealCanonical} that this
  action on horizontal arrows extends $T$ to a monad on $\MAT{V}$.
  (Actually, Seal shows that it is a ``lax extension of $T$ to
  $\mat{V}$ with op-lax unit and counit''; we will show in
  \S\ref{sectiontv} that this is the same as a monad on $\MAT{V}$.)

  In \cite{SealCanonical} this monad on $\MAT{V}$ is called the
  ``canonical extension'' of $T$ (note, however, that it is written
  backwards from his definition, as our Kleisli arrows will be $X \hto
  TY$, whereas his are $TX \hto Y$).  Since $\MAT{V}$ is isomorphic to
  its ``horizontal opposite,'' there is also an ``op-canonical
  extension'', which is in general distinct
  (although in some cases, such as for the ultrafilter monad, the two are identical).
  There are also many
  other extensions: for more detail, see \cite{ss:extn-laxalg}.
\end{example}

Another general way of constructing monads on virtual double categories
is to apply the construction $\MOD$ from the previous section, which
turns out to be a 2-functor.  Its 1-functoriality is fairly obvious and
was described in \cite[\S5.3]{leinster:higher-opds}; its action on
2-morphisms is given as follows.

\begin{definition}
Let $\mathbb{X} \two^{F}_G \mathbb{Y}$ be functors between virtual
double categories, and $F \to^{\theta} G$ a transformation.  One can
define a transformation
	\[ \MOD(F) \to^{\MOD( \theta )} \MOD(G) \]
whose vertical-arrow component at an object $(X_0, X, \bar{x}, \hat{x})$
is the monoid homomorphism
	\[ 
		\bfig
		\node a(0,375)[FX_0]
		\node b(500,375)[FX_0]
		\node c(0,0)[GX_0]
		\node d(500,0)[GX_0]
		\arrow/@{>}|-*@{|}/[a`b;F(X)]
		\arrow|b|/@{>}|-*@{|}/[c`d;G(X)]
		\arrow[a`c;\theta_{X_0}]
		\arrow|r|[b`d;\theta_{X_0}]
		\place(175,230)[\twoar(0,-1)]
		\place(325,180)[\scriptstyle{\theta_{X}}]
		\efig
	\]
and whose cell component at a horizontal arrow $(p, \bar{p}_r, \bar{p}_l)$ is given by
	\[
		\bfig
		
		\node a(0,375)[FX_0]
		\node b(500,375)[FY_0]
		\node c(0,0)[GX_0]
		\node d(500,0)[GY_0]
		\arrow/@{>}|-*@{|}/[a`b;Fp]
		\arrow|b|/@{>}|-*@{|}/[c`d;Gp]
		\arrow[a`c;\theta_{X_0}]
		\arrow|r|[b`d;\theta_{Y_0}]
		\place(200,230)[\twoar(0,-1)]
		\place(310,180)[\scriptstyle{\theta_p}]
		\efig
	\]
\end{definition}

\begin{proposition}
  With action on objects, 1-cells, and 2-cells described above, $\MOD$
  is an endo-2-functor of $\vDbl$.\qed
\end{proposition}

Note that the 2-functors $\PROF{(-)}$ and $\iPROF{-}$ can now be seen as
the composites of \Mod\ with $\MAT{(-)}$ and $\SPAN{-}$, respectively.

\begin{corollary}
  A monad $T$ on a virtual double category $\mathbb{X}$ induces a monad
  $\MOD(T)$ on $\MOD(\mathbb{X})$.\qed
\end{corollary}

\begin{example}
  Any monad $T$ on $\MAT{V}$ induces a monad on $\PROF{V}$.  For
  instance, this applies to the monads constructed in
  \autoref{example:tv}.
\end{example}

\begin{example}\label{eg:free-monoid-on-mat}
  Let \V\ be a symmetric monoidal category with an initial object
  $\emptyset$ preserved by $\otimes$.  Then the ``free monoid'' monad
  $T$ on \Set\ extends to a monad on $\MAT{V}$ as follows: a \V-matrix
  $X\hto^p Y$ is sent to the matrix $TX \hto^{Tp} TY$ defined by
  \[Tp\Big((y_1,\dots,y_m),(x_1,\dots,x_n)\Big) =
  \begin{cases}
    p(y_1,x_1)\ten \dots \ten p(y_n,x_n) & \text{if } n = m\\
    \emptyset &\text{if } n\neq m
  \end{cases}
  \]
  Applying $\MOD$, we obtain an extension of the ``free strict monoidal
  \V-category'' monad from $\cat{V}$ to $\PROF{V}$.
\end{example}

\begin{example}
  Likewise, any monad $T$ on $\SPAN{C}$ extends to a monad on
  $\iPROF{C}$.  But most interesting monads on $\SPAN{C}$ are induced from
  \C, so this gains us little beyond the observation that
  $\iPROF{-}$ is a 2-functor.
\end{example}

Not every monad on $\PROF{V}$ or $\iPROF{C}$ is induced by one on
$\MAT{V}$ or $\SPAN{C}$, however.  The following examples are also
important.

\begin{example}\label{eg:free-ssmc}
  Let \V\ be a symmetric monoidal category with finite
  colimits preserved by $\ten$ on both sides.  Then there is a ``free
  symmetric strict monoidal \V-category'' monad $T$ on $\cat{V}$,
  defined by letting the objects of $TX$ be finite lists of objects of
  $X$, with
  \[TX\Big((x_1,\dots,x_n),(y_1,\dots,y_m)\Big) =
  \begin{cases}
    \displaystyle\sum_{\sigma \in S_n} \bigotimes_{1\le i\le n} X(x_{\sigma(i)},y_i) & \text{if } n = m\\
    \emptyset &\text{if } n\neq m.
  \end{cases}
  \]
  A nearly identical-looking definition for profunctors extends this
  $T$ to a monad on $\PROF{V}$.  A similar definition applies for
  braided monoidal \V-categories.
\end{example}

\begin{example}\label{eg:free-safp}
  For \V\ as in \autoref{eg:free-ssmc}, there is also a ``free
  \V-category with strictly associative finite products'' monad on
  $\cat{V}$.  The objects of this $TX$ are again finite lists of
  objects of $X$, but now we have
  \[TX\Big((x_1,\dots,x_m),(y_1,\dots,y_n)\Big)
  = \prod_{1\le i\le n} \sum_{1\le j\le m} X(x_j,y_i).\]
  If \V\ is cartesian monoidal, then this can equivalently be written as
  \[TX\Big((x_1,\dots,x_m),(y_1,\dots,y_n)\Big)
  = \sum_{\alpha\colon n\rightarrow m} \prod_{1\le i\le n} X(x_{\alpha(i)},y_i).\]
  Again, a nearly identical definition for profunctors extends this to
  a monad on $\PROF{V}$.
\end{example}

\begin{example}\label{eg:free-coten}
  Monads that freely adjoin other types of limits and colimits also
  extend from $\cat{V}$ to $\PROF{V}$ in a similar way.  For instance,
  if \V\ is a locally finitely presentable closed monoidal category as
  in~\cite{kelly:enr-lfp}, there is a ``free \V-category with
  cotensors by finitely presentable objects'' monad on $\cat{V}$.  An
  object of $TX$ consists of a pair\footnote{To be precise, this
    definition only gives a pseudomonad on $\cat{V}$.  It is, however,
    easy to modify it to make a strict monad.} $(v;x)$ where $x\in
  X$ and $v\in \V$ is finitely presentable.  On homs we have
  \[TX\bigl((v;x),(w;y)\bigr) = \bigl[w,X(x,y)\ten v\bigr].\]
  As before, a nearly identical definition extends this to $\PROF{V}$.
\end{example}


\section{Generalized multicategories}
\label{sec:genmulti}

We now lack only one final ingredient for the definition of
generalized multicategories.  Since multicategories are like
categories, we expect them to also be monoids in some virtual double
category.  However, as we have seen in \S\ref{sec:2cat-vdc}, their
underlying data should include a horizontal arrow $A_0\hto TA_0$ rather
than $A_0\hto A_0$.  Thus we need to construct, given $T$ and
$\mathbb{X}$, a virtual double category in which the horizontal arrows
are horizontal arrows of the form $A_0\hto TA_0$ in $\mathbb{X}$.  This
is the purpose of the following definition.

\begin{definition}
  Let $T$ be a monad on a virtual double category $\mathbb{X}$.
  Define the \textbf{horizontal Kleisli virtual double category of $T$},
  $\Hkl{X}{T}$, as follows.
	\begin{itemize}
	\item Its vertical category is the same as that of $\mathbb{X}$.
	\item A horizontal arrow $X \hto^p Y$ is a horizontal arrow $X \hto^p TY$ in $\mathbb{X}$.
	\item A cell with nullary source uses the unit of the monad, so that a cell  
	\[
				\bfig
				\morphism<-250,-375>[X`Y;f]
				\morphism<250,-375>[X`Z;g]
				\morphism(-250,-375)|b|/@{>}|-*@{|}/<500,0>[Y`Z;p]
				\place(-35,-175)[\twoar(0,-1)]
				\place(50, -225)[\scriptstyle{\alpha}]
				\efig
	\]
	in $\Hkl{X}{T}$ is a cell
	\[
				\bfig
				\morphism<-250,-600>[X`Y;f]
				\morphism<125,-300>[X`TX;\eta]
				\morphism(125,-300)|r|<125,-300>[TX`TZ;Tg]
				\morphism(-250,-600)|b|/@{>}|-*@{|}/<500,0>[Y`TZ;p]
				\place(-35,-350)[\twoar(0,-1)]
				\place(50,-400)[\scriptstyle{\alpha}]
				\efig
	\]
	in $\mathbb{X}$ (note that $Tg \circ \eta = \eta\circ g$ by naturality).
	\item A cell with non-nullary source uses the multiplication of the monad, so that a cell
		\[
			\bfig
			\square/`>`>`/<2000,250>[X_0`X_n`Y_0`Y_1;`f`g`]
			\morphism|b|/@{>}|-*@{|}/<2000,0>[Y_0`Y_1;q]
			\morphism(0,250)/@{>}|-*@{|}/<500,0>[X_0`X_1;p_1]
			\morphism(500,250)/@{>}|-*@{|}/<500,0>[X_1`X_2;p_2]
			\morphism(1000,250)/@{>}|-*@{|}/<500,0>[X_2`\cdots;p_3]
			\morphism(1500,250)/@{>}|-*@{|}/<500,0>[\cdots`X_n;p_n]
			\place(965,150)[\twoar(0,-1)]
			\place(1050,100)[\scriptstyle{\alpha}]
			\efig
		\]
		in $\Hkl{X}{T}$ is a cell
		\[
			\bfig
			\morphism(0,500)|l|<0,-750>[X_0`Y_0;f]
			\morphism(0,-250)|b|/@{>}|-*@{|}/<2000,0>[Y_0`TY_1;q]
			\morphism(0,500)/@{>}|-*@{|}/<500,0>[X_0`TX_1;p_1]
			\morphism(500,500)/@{>}|-*@{|}/<500,0>[TX_1`T^2X_2;Tp_2]
			\morphism(1000,500)/@{>}|-*@{|}/<500,0>[T^2X_2`\cdots;T^2p_3]
			\morphism(1500,500)/@{>}|-*@{|}/<500,0>[\cdots`T^nX_n;T^{n-1}p_n]
			\morphism(2000,500)|r|<0,-300>[T^n X_n`\vdots;\mu]
			\morphism(2000,200)|r|<0,-200>[\vdots`TX_n;\mu]
			\morphism(2000,0)|r|<0,-250>[TX_n`TY_1;Tg]
			\place(965,175)[\twoar(0,-1)]
			\place(1050,125)[\scriptstyle{\alpha}]
			\efig
		\]
		in $\mathbb{X}$ (note that $Tg\circ \mu^{n-1} = \mu^{n-1} \circ T^n
    g$, by naturality).
	\item The composite of 
	\[
			\bfig
			\morphism(0,500)<200,0>[`;]
			\place(300,500)[\cdots]
			\morphism(400,500)<200,0>[`;]
			\morphism(0,500)<0,-250>[`;]
			\morphism(600,500)<0,-250>[`;]
			\morphism(0,250)<600,0>[`;]
			\place(300,375)[\alpha_1]
			
			\morphism(600,500)<200,0>[`;]
			\place(900,500)[\cdots]
			\morphism(1000,500)<200,0>[`;]
			\morphism(600,500)<0,-250>[`;]
			\morphism(1200,500)<0,-250>[`;]
			\morphism(600,250)<600,0>[`;]
			\place(900,375)[\alpha_2]
			
			\morphism(1800,500)<200,0>[`;]
			\place(2100,500)[\cdots]
			\morphism(2200,500)<200,0>[`;]
			\morphism(1800,500)<0,-250>[`;]
			\morphism(2400,500)<0,-250>[`;]
			\morphism(1800,250)<600,0>[`;]
			\place(2100,375)[\alpha_n]
			
			\place(1500,375)[\cdots]
			
			\morphism(0,250)<0,-250>[`;]
			\morphism(2400,250)<0,-250>[`;]
			\morphism<2400,0>[`;]
			\place(1200,125)[\beta]

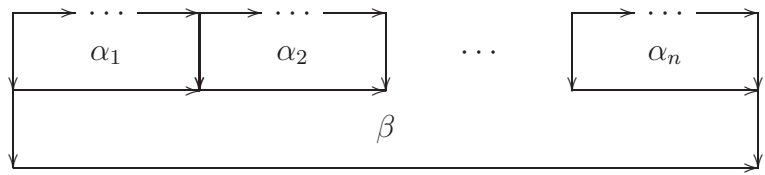
			\efig
	\]
	is given by the composite of
	\[
			\bfig
			\morphism(0,1800)<200,0>[`;]
			\place(300,1100)[\alpha_1]
			\place(300,1800)[\cdots]
			\morphism(400,1800)<200,0>[`;]
			
			\square(600,1600)<200,200>[```;```]
			\place(700,1700)[\mu]
			
			\place(700,1500)[\cdots]
			\square(600,1200)<200,200>[```;```]
			\place(700,1300)[\mu]
			
			\square(600,400)/`>``>/<600,800>[```;```]
			\morphism(0,400)<600,0>[`;]
			
			\square(1000,1600)<200,200>[```;```]
			\place(1100,1700)[\mu]
			
			\square(1000,1200)<200,200>[```;```]
			\place(1100,1300)[\mu]
			
			\place(900,1700)[\cdots]
			\place(1200,1500)[\cdots]
			\place(900,1300)[\cdots]

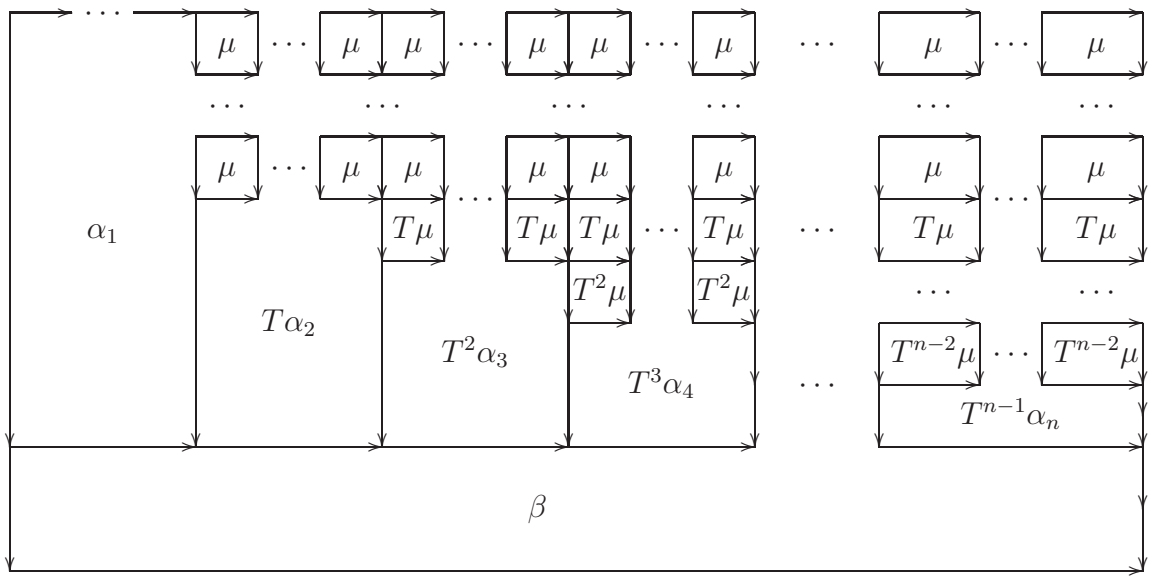
			\square(1200,1600)<200,200>[```;```]
			\place(1300,1700)[\mu]
						
			\square(1200,1200)<200,200>[```;```]
			\place(1300,1300)[\mu]
			
			\square(1200,1000)<200,200>[```;```]
			\place(1300,1100)[T\mu]
			
			\morphism(1200,1000)<0,-600>[`;]
			\place(900,800)[T \alpha_2]
			
			\place(1500,1200)[\cdots]
			\square(1600,1000)<200,200>[```;```]
			\place(1700,1100)[T \mu]
			\square(1600,1200)<200,200>[```;```]
			\place(1700,1300)[\mu]
			\square(1600,1600)<200,200>[```;```]
			\place(1700,1700)[\mu]
			
			\square(1800,1600)<200,200>[```;```]
			\place(1900,1700)[\mu]
			\square(1800,1200)<200,200>[```;```]
			\place(1900,1300)[\mu]
			\square(1800,1000)<200,200>[```;```]
			\place(1900,1100)[T\mu]
			\square(1800,800)<200,200>[```;```]
			\place(1900,900)[T^2\mu]
			
			\morphism(1800,800)<0,-400>[`;]
			\morphism(1000,400)<800,0>[`;]
			\place(1500,700)[T^2\alpha_3]
			
			\place(2100,1700)[\cdots]
			\place(1500,1700)[\cdots]
			\place(2300,1500)[\cdots]
			\place(1800,1500)[\cdots]
			\place(2100,1100)[\cdots]
			
			\square(2200,1600)<200,200>[```;```]
			\place(2300,1700)[\mu]
			\square(2200,1200)<200,200>[```;```]
			\place(2300,1300)[\mu]
			\square(2200,1000)<200,200>[```;```]
			\place(2300,1100)[T\mu]
			\square(2200,800)<200,200>[```;```]
			\place(2300,900)[T^2\mu]
			
			\morphism(2400,800)<0,-200>[`;]
			\morphism(2400,600)<0,-200>[`;]
			\morphism(1800,400)<600,0>[`;]
			\place(2100,600)[T^3\alpha_4]
			
			\square(2800,1600)<325,200>[```;```]
			\place(2975,1700)[\mu]
			\square(3325,1600)<325,200>[```;```]
			\place(3500,1700)[\mu]
			
			\square(2800,1200)<325,200>[```;```]
			\place(2975,1300)[\mu]
			\square(3325,1200)<325,200>[```;```]
			\place(3500,1300)[\mu]
			
			\square(2800,1000)<325,200>[```;```]
			\place(2975,1100)[T\mu]
			\square(3325,1000)<325,200>[```;```]
			\place(3500,1100)[T\mu]
			
			\place(2975,900)[\cdots]
			\place(3500,900)[\cdots]
			\place(3225,1200)[\cdots]
			\place(3225,1700)[\cdots]
			\place(2975,1500)[\cdots]
			\place(3500,1500)[\cdots]
			
			\place(2600,1700)[\cdots]
			\place(2600,1100)[\cdots]
			\place(2600,600)[\cdots]
			
			\square(2800,600)<325,200>[```;```]
			\place(2975,700)[T^{n-2}\mu]
			\square(3325,600)<325,200>[```;```]
			\place(3500,700)[T^{n-2}\mu]
			\place(3225,700)[\cdots]
			
			\morphism(2800,600)<0,-200>[`;]
			\morphism(2800,400)<850,0>[`;]
			\morphism(3650,600)<0,-100>[`;]
			\morphism(3650,500)<0,-100>[`;]
			\place(3225,500)[T^{n-1}\alpha_n]

			\morphism(3650,400)<0,-200>[`;]
			\morphism(3650,200)<0,-200>[`;]
			\morphism(0,1800)<0,-1400>[`;]
			\morphism(0,400)<0,-400>[`;]
			\morphism<3650,0>[`;]
			\place(1700,200)[\beta]

			\efig
	\]
	in $\mathbb{X}$.
	\item Identity cells use those of $\mathbb{X}$:
	\[
		\bfig
		\node A(0,250)[X]
		\node B(500,250)[TY]
		\node C(0,0)[X]
		\node D(500,0)[TY]
		\arrow/@{>}|-*@{|}/[A`B;p]
		\arrow|b|/@{>}|-*@{|}/[C`D;p]
		\arrow/=/[A`C;]
		\arrow/=/[B`D;]
		\place(200,175)[\twoar(0,-1)]
		\place(300,125)[\scriptstyle{1_P}]
		\efig
	\]
	\end{itemize}
\end{definition}

In general, the associativity for $\Hkl{X}{T}$ is shown by using the
(cell) naturality of $\mu$ and $\eta$, as well as the monad axioms.  The
general associativity is too large a diagram to show here; instead, we
will demonstrate a sample associativity calculation, which is
representative of the general situation.  Consider the following cells
in $\Hkl{X}{T}$:
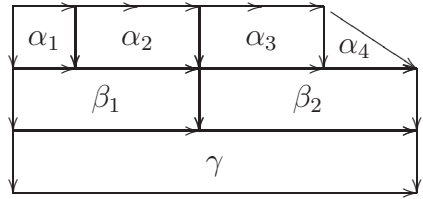
\[
	\bfig
	\square(0,400)<200,200>[```;```]
	\place(100,500)[\alpha_1]
	
	\morphism(200,600)<200,0>[`;]
	\morphism(400,600)<200,0>[`;]
	\square(200,400)/`>`>`>/<400,200>[```;```]
	\place(400,500)[\alpha_2]
	
	\morphism(600,600)<200,0>[`;]
	\morphism(800,600)<200,0>[`;]
	\square(600,400)/`>`>`>/<400,200>[```;```]
	\place(800,500)[\alpha_3]
	
	\morphism(1000,600)<300,-200>[`;]
	\morphism(1000,400)<300,0>[`;]
	\place(1100,475)[\alpha_4]
	
	\square(0,200)<600,200>[```;```]
	\place(300,300)[\beta_1]
	
	\square(600,200)<700,200>[```;```]
	\place(950,300)[\beta_2]
	
	\square<1300,200>[```;```]
	\place(650,100)[\gamma]
	
	\efig
\]
There are two possible ways to compose these cells: either composing the bottom first:
 \[(\gamma(\beta_1 \boxdot \beta_2))(\alpha_1 \boxdot \alpha_2 \boxdot \alpha_3 \boxdot \alpha_4) \]
or the top two first, followed by composition with the bottom:
 \[ \gamma ((\beta_1(\alpha_1 \boxdot \alpha_2)) \boxdot (\beta_2(\alpha_3 \boxdot \alpha_4))) \]
The first composite is given by the following composite in $\mathbb{X}$:

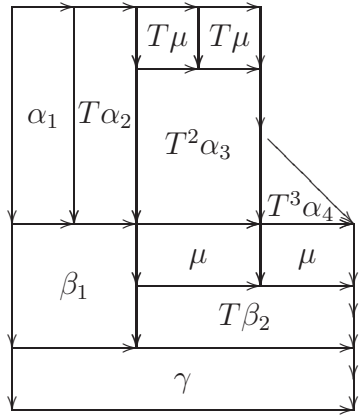
\[
	\bfig
	\square(0,600)<200,700>[```;```]
	\place(100,950)[\alpha_1]
	
	\square(200,600)/>`>``>/<200,700>[```;```]
	\morphism(400,1300)<0,-200>[`;]
	\morphism(400,1100)<0,-500>[`;]
	\place(300,950)[T\alpha_2]
	
	\square(400,1100)<200,200>[```;```]
	\place(500,1200)[T \mu]
	
	\square(600,1100)<200,200>[```;```]
	\place(700,1200)[T \mu]
	
	\square(400,600)/`>``>/<400,500>[```;```]
	\morphism(800,1100)<0,-200>[`;]
	\morphism(800,900)<0,-300>[`;]
	\place(600,850)[T^2\alpha_3]
	
	\morphism(800,900)<300,-300>[`;]
	\morphism(800,600)<300,0>[`;]
	\place(930,650)[T^3\alpha_4]
	
	\square(0,200)/`>``>/<400,400>[```;```]
	\morphism(400,600)<0,-200>[`;]
	\morphism(400,400)<0,-200>[`;]
	\place(200,400)[\beta_1]
	
	\square(400,400)<400,200>[```;```]
	\place(600,500)[\mu]
	
	\square(800,400)<300,200>[```;```]
	\place(950,500)[\mu]
	
	\square(400,200)/`>``>/<700,200>[```;```]
	\morphism(1100,400)<0,-100>[`;]
	\morphism(1100,300)<0,-100>[`;]
	\place(750,300)[T \beta_2]
	
	\square/`>``>/<1100,200>[```;```]
	\morphism(1100,200)<0,-100>[`;]
	\morphism(1100,100)<0,-100>[`;]
	\place(550,100)[\gamma]

	\efig
\]
By using cell naturality of $\mu$ twice, the above becomes
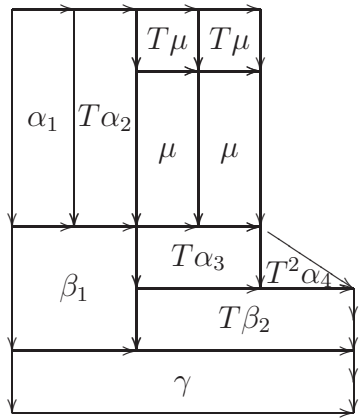
\[
	\bfig
	\square(0,600)<200,700>[```;```]
	\place(100,950)[\alpha_1]

	\square(200,600)/>`>``>/<200,700>[```;```]
	\morphism(400,1300)<0,-200>[`;]
	\morphism(400,1100)<0,-500>[`;]
	\place(300,950)[T\alpha_2]

	\square(400,1100)<200,200>[```;```]
	\place(500,1200)[T \mu]

	\square(600,1100)<200,200>[```;```]
	\place(700,1200)[T \mu]
	
	\square(400,600)<200,500>[```;```]
	\place(500,850)[\mu]
	
	\square(600,600)<200,500>[```;```]
	\place(700,850)[\mu]
	
	\square(0,200)/`>``>/<400,400>[```;```]
		\morphism(400,600)<0,-200>[`;]
		\morphism(400,400)<0,-200>[`;]
		\place(200,400)[\beta_1]
		
		\square(400,400)<400,200>[```;```]
		\place(600,500)[T \alpha_3]
		
		\morphism(800,600)<300,-200>[`;]
		\morphism(800,400)<300,0>[`;]
		\place(925,445)[T^2 \alpha_4]
		
		\square(400,200)/`>``>/<700,200>[```;```]
		\morphism(1100,400)<0,-100>[`;]
		\morphism(1100,300)<0,-100>[`;]
		\place(750,300)[T \beta_2]
		
		\square/`>``>/<1100,200>[```;```]
		\morphism(1100,200)<0,-100>[`;]
		\morphism(1100,100)<0,-100>[`;]
	\place(550,100)[\gamma]
	
	\efig
\]
we then use the monad axiom $T \mu \circ \mu = \mu \circ \mu$ to get
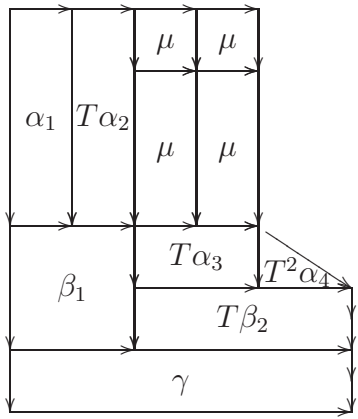
\[
	\bfig
	\square(0,600)<200,700>[```;```]
	\place(100,950)[\alpha_1]

	\square(200,600)/>`>``>/<200,700>[```;```]
	\morphism(400,1300)<0,-200>[`;]
	\morphism(400,1100)<0,-500>[`;]
	\place(300,950)[T\alpha_2]

	\square(400,1100)<200,200>[```;```]
	\place(500,1200)[\mu]

	\square(600,1100)<200,200>[```;```]
	\place(700,1200)[\mu]
	
	\square(400,600)<200,500>[```;```]
	\place(500,850)[\mu]
	
	\square(600,600)<200,500>[```;```]
	\place(700,850)[\mu]
	
	\square(0,200)/`>``>/<400,400>[```;```]
	\morphism(400,600)<0,-200>[`;]
	\morphism(400,400)<0,-200>[`;]
	\place(200,400)[\beta_1]

	\square(400,400)<400,200>[```;```]
	\place(600,500)[T \alpha_3]

	\morphism(800,600)<300,-200>[`;]
	\morphism(800,400)<300,0>[`;]
	\place(925,445)[T^2 \alpha_4]

	\square(400,200)/`>``>/<700,200>[```;```]
	\morphism(1100,400)<0,-100>[`;]
	\morphism(1100,300)<0,-100>[`;]
	\place(750,300)[T \beta_2]

	\square/`>``>/<1100,200>[```;```]
	\morphism(1100,200)<0,-100>[`;]
	\morphism(1100,100)<0,-100>[`;]
	\place(550,100)[\gamma]
	
	\efig
\]	
which is the second composite $\gamma ((\beta_1(\alpha_1 \boxdot \alpha_2)) \boxdot (\beta_2(\alpha_3 \boxdot \alpha_4)))$.  
	
\begin{remark}\label{thm:not-kleisli-object}
  Unfortunately, we do not know of any universal property satisfied by this construction.
  In particular, $\Hkl{X}{T}$ is not a Kleisli object for $T$ in \vDbl\
  in the sense of~\cite{street:ftm}; the latter would instead contain
  \emph{vertical} Kleisli arrows.
  In fact, for general $\mathbb{X}$ there need not even be a canonical
  functor $\mathbb{X}\to \Hkl{X}{T}$.
\end{remark}

We can now give our first preliminary definition of generalized
multicategories relative to a monad $T$.

\begin{definition}
  Let $T$ be a monad on a virtual double category
  $\mathbb{X}$.  A \textbf{$T$-monoid} is defined to be a monoid in
  $\Hkl{X}{T}$, and likewise for a \textbf{$T$-monoid homomorphism}.  We
  denote the virtual double category $\MOD(\Hkl{X}{T})$, whose objects
  are $T$-monoids, by $\KMod(\mathbb{X},T)$.
\end{definition}

As a reference, the data for a $T$-monoid consists of an object $X_0 \in
\mathbb{X}$, a horizontal arrow $X_0 \hto^{X} TX_0$ in $\mathbb{X}$, and
cells
\[	
	\bfig
	\square/`=`>`/<1000,250>[X_0`T^2X_0`X_0`TX_0;``\mu`]
	\morphism|b|/@{>}|-*@{|}/<1000,0>[X_0`TX_0;X]
	\morphism(0,250)/@{>}|-*@{|}/<500,0>[X_0`TX_0;X]
	\morphism(500,250)/@{>}|-*@{|}/<500,0>[TX_0`T^2X_0;TX]
	\place(475, 150)[\twoar(0,-1)]
	\place(550,100)[\scriptstyle{\bar{x}}]
	
	\place(1500,100)[\mbox{and}]
			
	\morphism(2000,300)/=/<-250,-375>[X_0`X_0;]
	\morphism(2000,300)<250,-375>[X_0`TX_0;\eta]
	\morphism(1750,-75)|b|/@{>}|-*@{|}/<500,0>[X_0`TX_0;X]
	\place(1975,100)[\twoar(0,-1)]
	\place(2050,50)[\scriptstyle{\hat{x}}]
	
	\efig
\]
Note that these cells have precisely the forms we predicted at the
beginning of \S\ref{sec:2cat-vdc}.

\begin{remark}\label{rmk:functoriality}
  We have seen in \S\ref{sec:fcmulticats} that $\MOD$ is a 2-functor.
  In fact, under suitable hypotheses (involving the notions of
  restriction and composites to be introduced in
  \S\ref{sec:composites-units} and \S\ref{sec:pfrbi}), $\HKl$ is also
  a (pseudo) functor, and thus so is $\KMod$.  In fact, $\HKl$ is a
  pseudo functor in two different ways, corresponding to the two
  different kinds of morphisms of monads: lax and colax.  This was
  observed by~\cite{leinster:higher-opds} in his context; we will
  discuss the functoriality of $\HKl$ and $\KMod$ in our framework in
  the forthcoming~\cite{cs:functoriality}.
\end{remark}

We now consider some examples.

\begin{example}
  Of course, if $T$ is the identity monad on any $\mathbb{X}$, then a
  $T$-monoid is just a monoid, and $\KMod(\mathbb{X},T) =
  \Mod(\mathbb{X})$.
\end{example}

Recall from \autoref{example:cartesian} that any pullback-preserving monad
$\mb{C} \to^T \mb{C}$ extends to a monad on $\SPAN{C}$.

\begin{example}
  If $T$ is the ``free monoid'' monad extended to $\SPAN{Set}\iso
  \MAT{Set}$, then a $T$-monoid consists of a set $A_0$, a \Set-matrix
  $\{A((x_1,\dots,x_n),y)\}_{x_i,y\in A_0}$, and composition and
  identity functions.  It is easy to see that this reproduces the
  notion of an ordinary multicategory.  Likewise, $T$-monoid
  homomorphisms are functors between multicategories.
\end{example}

\begin{example}
  If $T$ is the ``free category'' monad on directed graphs, then a
  $T$-monoid is a virtual double category.  (This is, of course, the
  origin of the name ``fc-multicategory,'' where $\mathbf{fc}$ is a
  name for this monad.)  The vertices and edges of
  the directed graph $A_0$ are the objects and horizontal arrows,
  respectively, while in the span $A_0 \leftarrow A \rightarrow TA_0$
  the vertices of $A$ are the vertical arrows and its edges are the
  cells.  Likewise, $T$-monoid homomorphisms are functors between
  virtual double categories.
\end{example}

\begin{example}
  Let $M$ be a monoid and $T=(M\times -)$ the ``free $M$-set'' monoid on
  $\SPAN{Set}$.  A $T$-monoid consists of a set $A_0$ and a family of
  sets $\{A(m;x,y)\}_{x,y\in A_0, m\in M}$.  The composition and
  identity functions make it into an \emph{$M$-graded category}, i.e.\
  a category in which every arrow is labeled by an element of $M$ in a
  way respecting composition and identities.  The case $M=\mathbb{Z}$
  may be most familiar.
\end{example}

\begin{example}
  Let $\mb{C}$ be a lextensive category and $T$ the monad $(-)+1$ on
  $\SPAN{C}$.  A $T$-monoid consists of an object $A_0$ and a span
  $A_0 \leftarrow A \rightarrow A_0+1$; since $\mb{S}$ is extensive,
  $A$ decomposes into two spans $A_0 \leftarrow A_1 \rightarrow A_0$
  and $A_0 \leftarrow B \rightarrow 1$.  The composition and identity
  functions then make the first span into an internal category in
  $\mb{C}$ and the second into an internal diagram over this category.
\end{example}

\begin{example}\label{eg:cat-over}
  Let $\mb{S}$ be a small category, and let $T$ be the monad on
  $\Set^{\ob(\mb{S})}$ whose algebras are functors
  $\mb{S}\to\Set$.  Thus, for a family $\{A_x\}_{x\in\ob(\mb{S})}$
  in $\Set^{\ob(\mb{S})}$ we have
  \[(TA)_x = \coprod_{y\in \ob(\mb{S})} A_y\times \mb{S}(x,y).
  \]
  This $T$ preserves pullbacks, so it induces a monad on
  $\SPAN{\Set^{\ob(\mb{S})}}$.  A $T$-monoid $A\hto^M TA$ can be
  identified with a category over $\mb{S}$, i.e\ a functor
  $\mb{A}\to\mb{S}$.  Namely, the elements of the set $A_x$ are the
  objects of the fiber of $\mb{A}$ over $x\in\ob(\mb{S})$, while $M$
  can be broken down into a collection of spans
  \[A_x \toleft M_{x,y} \to A_y\times \mb{S}(x,y),
  \]
  which together supply the arrows of $\mb{A}$ and their images in
  $\mb{S}$.  The morphisms of $T$-monoids are likewise the functors
  over $\mb{S}$.
\end{example}

\begin{example}\label{eg:glob-opd}
  If $T$ is the ``free strict $\omega$-category'' monad on
  $\SPAN{Globset}$, then a $T$-monoid of the form $1\hto T1$ can be
  identified with a \emph{globular operad} in the sense
  of~\cite{batanin:monglob}, as described
  in~\cite{leinster:higher-opds}.  General $T$-monoids are
  \emph{globular multicategories} (or \emph{many-sorted globular
    operads}) as considered
  in~\cite[p.273--274]{leinster:higher-opds}.
\end{example}

Recall from \autoref{example:tv} that any taut monad $T$ on \Set\
(such as the identity monad, the powerset monad, the filter monad, or
the ultrafilter monad) extends to a monad on $\MAT{V}$ for any
completely distributive quantale $\mb{V}$ (such as $\mb{2}$ or \eR).
We will show in \S\ref{sectiontv} that in such a case, our $T$-monoids
are the same as the $(T,\mb{V})$-algebras studied
by~\cite{CommonApproach,OneSetting,SealCanonical}, and others; thus we
have the following examples.

\begin{example}
If $T$ is the identity monad, then $\KMod(\MAT{V}, T) = \PROF{V}$.
Thus, for $\mb{V} = \mb{2}$, $T$-monoids are preorders; and for $\mb{V}
= \eR$, $T$-monoids are metric spaces (in the sense of
\cite{LawvereMetric}).
\end{example}

\begin{example}
  If $T$ is the ultrafilter monad, and $\mb{V} = \mb{2}$, then a
  $T$-monoid consists of a set equipped with a binary relation between
  ultrafilters and points satisfying unit and composition axioms.  If
  we call this relation ``convergence,'' then the axioms precisely
  characterize the convergence relation in a topological space; thus
  $T$-monoids are topological spaces, and $T$-monoid homomorphisms are
  continuous functions.  This observation is originally due to
  \cite{RelationalAlgebras} (note that although his construction of an
  ultrafilter monad on $\Rel$ looks quite different, it is in fact the
  same).
\end{example}

\begin{example}\label{eg:closure-spaces}
  If $T$ is the powerset monad, and $\mb{V} = \mb{2}$, then $T$-monoids
  are \emph{closure spaces}.  A closure space consists of a set $A$
  equipped with an operation $c(-)$ on subsets which is:
  \begin{packeditem}
  \item extensive: $X\subseteq c(X)$, 
  \item monotone: $Y\subseteq X \Rightarrow c(Y)\subseteq c(X)$, and
  \item idempotent: $c(c(X))\subseteq c(X)$.
  \end{packeditem}
\end{example}

\begin{example}
  If $T$ is the ultrafilter monad and $\mb{V} = \eR$, then $T$-monoids
  are equivalent to \emph{approach spaces}.  An approach space is a set $X$ equipped
  with a function $d: X \times PX \to [0,\infty]$ such that
	\begin{packeditem}
	\item $d(x,\{x\}) = 0$,
	\item $d(x,\emptyset) = \infty$,
	\item $d(x,A \cup B) = \mbox{min} \{ d(x,A), d(x,B) \}$, and
	\item $\forall \epsilon \geq 0, d(x,A) \leq d(x, \{x: d(x,A) \leq \epsilon \}) + \epsilon$.
	\end{packeditem}
Approach spaces have found applications in approximation theory,
products of metric spaces, and measures of non-compactness: for more
detail, see \cite{lowen:approach}.
\end{example}

Finally, we consider $T$-monoids relative to the additional examples
of monads on $\PROF{V}$ from the end of \S\ref{sec:2cat-vdc}.

\begin{example}\label{eg:sym-multicat}
  Let $T$ be the ``free symmetric strict monoidal \V-category'' monad
  on $\PROF{V}$ from \autoref{eg:free-ssmc}.  If $A_0$ is a discrete
  \V-category, then a $T$-monoid
  $A_0 \hto^A TA_0$ is a symmetric \V-enriched multicategory (known to
  some authors as simply a ``multicategory'').  Likewise, from the
  ``free braided strict monoidal \V-category'' monad we obtain
  braided multicategories.

  If $A_0$ is not discrete, then a $T$-monoid (for $\V=\Set$) is a
  symmetric multicategory in the sense of \cite{bd:hda3} and
  \cite{cheng:opetopic}: in addition to the multi-arrows, there is
  also another type of arrow between the objects of the multicategory
  which forms a category, and which acts on the multi-arrows.
\end{example}

\begin{example}\label{eg:lawvere-theories}
  Let $T$ be the ``free category with strictly associative finite
  products'' monad on $\PROF{Set}$ from \autoref{eg:free-safp}.  If
  $A_0$ is a one-object discrete category, then a $T$-monoid $A_0
  \hto^A TA_0$ is a Lawvere theory, as in~\cite{lawvere:functsem}.
  If $A_0$ has more than one object, but is still discrete, then a
  $T$-monoid $A_0 \hto^A TA_0$ is a ``multi-sorted'' Lawvere theory.

  This is a little different from the more usual definition of Lawvere
  theory, but the equivalence between the two is easy to see.  A
  Lawvere theory is commonly defined to be a category $A$ with object
  set $\mathbb{N}$ such that each object $n$ is the $n$-fold product
  $1^n$.  This implies that $A(m,n) \iso A(m,1)^n$, so it is
  equivalent to give just the collection of sets $A(m,1)$ with
  suitable additional structure.  Since $T1$ has object set
  $\mathbb{N}$, a $T$-monoid $1\hto T1$ also consists of sets $A(m,1)$
  for $m\in\mathbb{N}$, and it is then straightforward to verify that
  the additional structures in the two cases are in bijective
  correspondence.  Note, however, that the \emph{morphisms} between
  such $T$-monoids do not correspond to all of the morphisms between
  theories considered in~\cite{lawvere:functsem}, but only those of
  ``degree one;'' the others are only visible from the ``category with
  object set $\mathbb{N}$'' viewpoint.

  The relationship between these two definitions of Lawvere theory is
  analogous to the way in which an \emph{operad} can also be defined
  as a certain sort of monoidal category with object set $\mathbb{N}$.
  In fact, both arise from a very general \emph{adjunction} between
  $T$-algebras and $T$-monoids; see \autoref{rmk:free-models} and the
  forthcoming~\cite{cs:laxidem}.
\end{example}

\begin{example}\label{eg:enr-lawvth}
  If $T$ is the ``free \V-category with strictly associative finite
  products'' monad on $\PROF{V}$ from \autoref{eg:free-safp} and $A_0$
  is a one-object discrete \V-category, then a $T$-monoid $A_0 \hto^A
  TA_0$ is a ``\V-enriched finite product theory.''  If $A_0$ is
  unchanged but $T$ is instead the ``free \V-category with finitely
  presentable cotensors'' monad from \autoref{eg:free-coten}, then a
  $T$-monoid $A_0 \hto^A TA_0$ is a ``Lawvere \V-theory'' as defined
  in~\cite{power:enr-theories} (with the same caveat as in the
  previous example).  To obtain a ``multi-sorted Lawvere
  \V-theory'' we need $T$ to adjoin both finite products and finite
  cotensors.
\end{example}

\begin{example}\label{eg:clubs}
  If $T$ is any of
  \begin{packeditem}
    \item the ``free symmetric strict monoidal category'' monad,
    \item the ``free category with strictly associative finite
      products'' monad, or
    \item the ``free category with strictly associative finite
      coproducts'' monad,
  \end{packeditem}
  but now considered as a monad on $\SPAN{Cat}$, then a
  $T$-monoid with a discrete category of objects is a
  \emph{club} in the sense of~\cite{kelly:mv-funct-calc}
  and~\cite{kelly:abst-coh} (relative to $P$, $S$, or $S^{op}$, in
  Kelly's terminology).  See also \S\ref{sec:clubs}.
\end{example}

\begin{remark}\label{rmk:structure-types}
  When $T$ is the ``free symmetric strict monoidal category'' monad on
  $\PROF{Set}$, the horizontal arrows between discrete categories in
  $\Hkl{\PROF{Set}}{T}$ are the \emph{generalized species of
    structure} of~\cite{fghw:gen-species} (called \emph{structure types}
  in~\cite{bd:fin-feyn}).  The \emph{esp\'eces de structures}
  of~\cite{joyal:species,joyal:funcan-species} are the particular case
  of horizontal arrows $1\hto 1$ in $\Hkl{\PROF{Set}}{T}$.  Likewise,
  when $T$ is the analogous monad on $\SPAN{Gpd}$, the horizontal
  arrows in $\Hkl{\SPAN{Gpd}}{T}$ are the (generalized) \emph{stuff
    types} of~\cite{bd:fin-feyn,morton:catalg-qm}.
\end{remark}

We can see from these examples that for virtual double categories
whose objects are ``category-like,'' it is often the $T$-monoids whose objects are
discrete which are of particular interest.  We will make this notion
precise in \S\ref{sec:normalization}, and propose that often a better
solution is to consider ``normalized'' $T$-monoids.  

First, however, we must develop some additional machinery for virtual
double categories.  We will describe when horizontal arrows have units
and composites, as well as when horizontal arrows can be ``restricted''
along vertical arrows.  With this theory in hand, we can then return to
study ``object discrete'' and ``normalized'' $T$-monoids, as well as
when such $T$-monoids are ``representable''.

\begin{remark}
  If $\V$ is a complete and cocomplete closed symmetric monoidal
  category, then the virtual double category $\PROF{V}$ is itself
  ``almost'' of the form $\Hkl{X}{T}$.  We take $\mathbb{X}$ to be the
  double category $\SQ(\cat{V})$, whose objects are \V-categories,
  whose vertical and horizontal arrows are both \V-functors, and whose
  cells are \V-natural transformations, and we define $T A =
  \V^{A^{\op}}$ to be the enriched presheaf category of $A$.  The
  observation is then that a \V-profunctor $A\hto B$ can be identified
  with an ordinary \V-functor $A\to \V^{B^{\op}}$, so that $\PROF{V}$
  is almost the same as $\HKl(\SQ(\cat{V}),T)$.  This is not quite
  right, since $T$ is not really a monad due to size issues.  But
  these problems can be dealt with, for instance by using ``small
  presheaves'' as in~\cite{dl:lim-smallfr}.

  Assuming the functoriality of $\HKl$ mentioned in
  \autoref{rmk:functoriality}, this observation implies that if $S$ is
  another monad on $\cat{V}$ related to $T$ by a distributive
  law, or equivalently a monad in the 2-category of
  monads and monad morphisms (see~\cite{beck:dl,street:formalmonads}), then
  $S$ induces a monad on $\PROF{V}$, which we can in turn use to define
  generalized multicategories as $S$-monoids in $\PROF{V}$.  For
  example, since a symmetric
  monoidal structure on $A$ extends to $T A$ by Day convolution, the
  ``free symmetric monoidal \V-category'' monad distributes over $T$,
  inducing its extension to $\PROF{V}$ considered in
  \autoref{eg:free-ssmc}.  This is the argument used
  in~\cite{fghw:gen-species} to construct the bicategory
  $\calH(\HKl(\PROF{V},T))$.  Similar arguments apply to the monad
  from \autoref{eg:free-safp}.
\end{remark}


\section{Composites and units}
\label{sec:composites-units}

In \S\ref{sec:fcmulticats} we introduced (virtual) double categories
as a framework in which one can define monoids and monoid
homomorphisms so as to include both enriched and internal categories
with the appropriate notions of functor.
However, we would certainly like to be able to recover \emph{natural
  transformations} as well, but this requires more structure than is
present in a virtual double category.

It is not hard to see that the vertical category of any (pseudo) double
category can be enriched to a \emph{vertical 2-category}, whose
2-cells $f\Rightarrow g$ are the squares of the form
\[\vcenter{\xymatrix{A \ar@{=}[r]|-@{|}^-{} \ar[d]_g \ar@{}[dr]|{\twoar(0,-1)} & A \ar[d]^f\\
    B\ar@{=}[r]|-@{|}_-{} & B,}}
\]
and that in examples such as $\iPROF{C}$ and $\PROF{V}$ these 2-cells
recover the appropriate notion of natural transformation.  In a
virtual double category this definition is impossible, since there may
not be any horizontal identity arrows.  However, it turns out that we
can characterize those horizontal identities which \emph{do} exist by
means of a universal property.  In fact, it is not much harder to
characterize arbitrary horizontal composites (viewing identities as
0-ary composites).  In this section we study such composites; in the
next section we will use them to define vertical 2-categories.

\begin{definition}\label{defn:pfrbi-opcart}
  In a virtual double category, a cell
  \[\xymatrix{ \ar[r]|-@{|}^{p_1}\ar@{=}[d] \ar@{}[drrr]|{{\textstyle\Downarrow} \mathrm{opcart}} &
    \ar[r]|-@{|}^<>(.5){p_2} & \dots \ar[r]|-@{|}^<>(.5){p_n} &
    \ar@{=}[d]\\
    \ar[rrr]|-@{|}_q &&& }\]
  is said to be \textbf{opcartesian} if any cell
  \[\xymatrix{
    \ar[r]|-@{|}^{r_1}\ar[d]_f &
    \ar[r]|-@{|}^<>(.5){r_2} & \dots \ar[r]|-@{|}^<>(.5){r_m} &
    \ar[r]|-@{|}^{p_1}  \ar@{}[drrr]|{{\textstyle\Downarrow}} &
    \ar[r]|-@{|}^<>(.5){p_2} & \dots \ar[r]|-@{|}^<>(.5){p_n} &
    \ar[r]|-@{|}^{s_1} &
    \ar[r]|-@{|}^<>(.5){s_2} & \dots \ar[r]|-@{|}^<>(.5){s_k} &
    \ar[d]^g\\
    \ar[rrrrrrrrr]|-@{|}_t &&&&&&&&& }\]
  factors through it uniquely as follows:
  \[\xymatrix{
    \ar[r]|-@{|}^<>(.5){r_1}\ar@{=}[d] &
    \ar[r]|-@{|}^<>(.5){r_2} \ar@{=}[d] &
    \dots \ar[r]|-@{|}^<>(.5){r_m}
    \ar@{=}@<2mm>[d] \ar@{=}@<-2mm>[d] \ar@{}[d]|{\cdots} &
    \ar[r]|-@{|}^<>(.5){p_1} \ar@{=}[d]  \ar@{}[drrr]|{{\textstyle\Downarrow} \mathrm{opcart}} &
    \ar[r]|-@{|}^<>(.5){p_2} & \dots \ar[r]|-@{|}^<>(.5){p_n} &
    \ar[r]|-@{|}^<>(.5){s_1}\ar@{=}[d] &
    \ar[r]|-@{|}^<>(.5){s_2} \ar@{=}[d] &
    \dots \ar[r]|-@{|}^<>(.5){s_k}
    \ar@{=}@<2mm>[d] \ar@{=}@<-2mm>[d] \ar@{}[d]|{\cdots} &
    \ar@{=}[d]\\
    \ar[r]|-@{|}_<>(.5){r_1}\ar[d]_f &
    \ar[r]|-@{|}_<>(.5){r_2} & \dots \ar[r]|-@{|}_<>(.5){r_m} &
    \ar[rrr]|{q}  \ar@{}[drrr]|{{\textstyle\Downarrow}} &&&
    \ar[r]|-@{|}_<>(.5){s_1} &
    \ar[r]|-@{|}_<>(.5){s_2} &
    \dots \ar[r]|-@{|}_<>(.5){s_k} &
    \ar[d]^g\\
    \ar[rrrrrrrrr]|-@{|}_t &&&&&&&&& .}\]
\end{definition}

If a string of composable horizontal arrows is the source of some opcartesian
cell, we say that it \textbf{has a composite}.  We refer to the
target $q$ of that cell as a \textbf{composite} of the given string
and write it as $p_1\odot\dots\odot p_n$.  Similarly, if $n=0$ and
there is an opcartesian cell of the form
\[\xymatrix@R=1pc@C=1pc{ &X \ar@{}[dd]|(.6){{\textstyle\Downarrow}}
  \ar@{=}[ddl]\ar@{=}[ddr]\\
  \\ X \ar[rr]|-@{|}_{U_X} && X.}
\]
we say that $X$ \textbf{has a unit} $U_X$.

These universal properties make it easy to show that composites and
units, when they exist, behave like composites and units in a pseudo
double category.  For example, composites and units are unique up to
isomorphism: given two opcartesian arrows with the same
source, factoring each through the other gives an isomorphism between
their targets.  Likewise, the composite of opcartesian cells is
opcartesian, so composition is associative up to coherent isomorphism
whenever all relevant composites exist.  More precisely, if $p\odot q$
exists, then $(p\odot q)\odot r$ exists if and only if $p\odot q\odot
r$ exists, and in that case they are isomorphic.  It follows that if
$p\odot q$ and $q\odot r$ exist, then
\[(p\odot q)\odot r \iso p\odot (q\odot r),
\]
each existing if the other does.  Similarly, any string in which
all but one arrow is a unit:
\[\xymatrix{X \ar[r]|-@{|}^{U_X} &
  \dots \ar[r]|-@{|}^{U_X} &
  X \ar[r]|-@{|}^p &
  Y \ar[r]|-@{|}^{U_Y} &
  \dots \ar[r]|-@{|}^{U_Y}  & Y}\]
has a composite, which is (isomorphic to) $p$.

The following theorem, which was also observed in~\cite{dpp:paths-dbl},
is a straightforward generalization of the relationship between monoidal
categories and ordinary multicategories described in
\cite{HermidaRepresentable}.  It is also a special case of the general
relationship between pseudo algebras and generalized multicategories, as
in \cite[\S6.6]{leinster:higher-opds}, \cite{HermidaCoherent}, and
\S\ref{sec:vert-alg} of the present paper.

\begin{theorem}\label{dblcatcharacterization}
  A virtual double category is a pseudo double category if and only if
  every string of composable horizontal arrows (including zero-length
  ones) has a composite.
\end{theorem}
\begin{proof}[(sketch)]
  ``Only if'' is clear, by definition of how a pseudo double category
  becomes a virtual one.  For ``if'', we use the isomorphisms
  constructed above; we invoke again the universal property of
  opcartesian cells to show coherence.
\end{proof}

\begin{example}
  If \V\ has an initial object $\emptyset$ which is preserved by
  $\otimes$ on both sides, then $\MAT{V}$ has units: the unit of a set
  $X$ is the matrix
  \[U_X(x,x') =
  \begin{cases}
    I & x=x'\\
    \emptyset & x\neq x'.
  \end{cases}\]
  If \V\ has all small coproducts which are preserved by $\otimes$ on
  both sides, then $\MAT{V}$ has composites given by ``matrix
  multiplication.''  For instance, the composite of matrices
  $X\hto^p Y$ and $Y\hto^q Z$ is
  \[(p\odot q)(z,x) = \coprod_{y\in Y} p(y,x)\otimes q(z,y).\]
\end{example}

\begin{example}
  Since $\SPAN{C}$ is a pseudo double category, all composites and
  units always exist.  Composites are given by pullback, and the unit
  of $X$ is the unit span $X\leftarrow X \rightarrow X$.
\end{example}

Regarding units in $\PROF{V}$ and $\iPROF{C}$, we have the following.

\begin{proposition}\label{prop:mod-units}
  For any virtual double category $\mathbb{X}$, all units exist in
  $\MOD(\mathbb{X})$.  For any monoid $A$, its unit cell
  \[\bfig\scalefactor{.7} \Atriangle/=`=`->/[A_0`A_0`A_0;``A]
  \morphism(0,0)|m|<1000,0>[A_0`A_0;|]
  \place(500,200)[\Downarrow\scriptstyle \hat{a}]
  \efig
  \]
  is opcartesian in $\MOD(\mathbb{X})$.  Therefore, $U_A$ is $A$
  itself, regarded as an $A$-$A$-bimodule.
\end{proposition}
\begin{proof}
  Firstly, the unit axioms of $A$ show that $\hat{a}$ is, in fact, a
  cell in $\MOD(\mathbb{X})$.  Now we must show that composing with
  $\hat{a}$ gives a bijection between cells
  \[\bfig\scalefactor{.7}
  \square/`>`>`>/<1500,500>[B`C`D`E;`g`f`P]
  \morphism(0,500)/{-->}/<750,0>[B`A;p_1\dots p_m]
  \morphism(750,500)/{-->}/<750,0>[A`C;q_1\dots q_n]
  \place(750,250)[\Downarrow\scriptstyle \alpha]
  \efig \qquad\text{and}\qquad
  \bfig\scalefactor{.7}
  \square/`>`>`>/<1850,500>[B`C`D`E;`g`f`p]
  \morphism(0,500)/{-->}/<750,0>[B`A;p_1\dots p_m]
  \morphism(750,500)<350,0>[A`A;A]
  \morphism(1100,500)/{-->}/<750,0>[A`C;q_1\dots q_n]
  \place(925,250)[\Downarrow\scriptstyle \beta]
  \efig
  \]
  in $\MOD(\mathbb{X})$.  Clearly composing any cell $\beta$ of the
  second form with $\hat{a}$ gives a cell $\alpha$ of the first form.
  Conversely, given $\alpha$ of the first form, there are two cells
  of the second form defined by letting $A$ act first on $p_m$ from
  the right and $q_1$ from the left, respectively.  These are equal by
  one of the axioms for $\alpha$ to be a cell in $\MOD(\mathbb{X})$;
  we let $\beta$ be their common value.  (In the case when $m=0$ or
  $n=0$, we use the action of $f$ or $g$ instead.)  The other axioms
  for $\alpha$ then carry over to $\beta$ to show that it is a cell
  in $\MOD(\mathbb{X})$.

  The unit axioms for the action of $A$ on bimodules show that
  $\alpha\mapsto\beta\mapsto \alpha$ is the identity, while the
  equivariance axioms for $\beta$ regarding the two actions of $A$
  show that $\beta\mapsto\alpha\mapsto\beta$ is the identity.  Thus we
  have a bijection, as desired.
\end{proof}

Therefore, since $\PROF{V}=\MOD(\MAT{V})$ and $\iPROF{C}=\MOD(\SPAN{C})$, they both
always have all units.  By contrast, extra assumptions on $\mathbb{X}$
are required for composites to exist in $\MOD(\mathbb{X})$; here are
the two examples of greatest interest to us.

\begin{example}
  If \V\ has small colimits preserved by $\otimes$ on both sides, then
  $\PROF{V}$ has all composites; the composite of enriched profunctors
  $A\hto^p B$ and $B\hto^q C$ is given by the coend
  \[(p\odot q)(z,x) = \int^{\mathrlap{y\in B}} p(y,x)\otimes q(z,y).\]
\end{example}

\begin{example}
  If \C\ has coequalizers preserved by pullback, then $\iPROF{C}$ has
  all composites; the composite of internal profunctors is an
  ``internal coend.''
\end{example}

Together with \autoref{prop:mod-units}, these examples will suffice
for the moment.  In appendix \ref{sec:composites-mod-hkl} we will give
general sufficient conditions for composites to exist in
$\MOD(\mathbb{X})$, and for composites and units to exist in
$\Hkl{X}{T}$.

\begin{remark}\label{thm:weak-opcart}
  If \autoref{defn:pfrbi-opcart} is satisfied only for $m=k=0$,
  we say that the cell is \textbf{weakly opcartesian}.  We do not
  regard a weakly opcartesian cell as exhibiting its target as a
  composite of its source, since the weak condition is insufficient to
  prove associativity and unitality.  However, a weakly opcartesian
  cell does suffice to \emph{detect} its target as the composite
  of its source, if we already know that that composite exists.
  Furthermore, if any composable string in $\mathbb{X}$ is the source
  of a weakly opcartesian cell and moreover weakly opcartesian
  cells are closed under composition, then one can show, as in
  \cite{HermidaRepresentable}, that in fact every weakly opcartesian
  cell is cartesian; see also \S\ref{sec:vert-alg}.

  Virtual double categories having only weakly opcartesian cells seem
  to be fairly rare; one example is $\MAT{V}$ where \V\ has colimits
  which are not preserved by its tensor product.  Note that in this
  case, $\PROF{V}$ need not even have weakly opcartesian cells, since
  we require $\otimes$ to preserve coequalizers simply to make the
  composite of two profunctors into a profunctor.
\end{remark}

If $\mathbb{X}$ and $\mathbb{Y}$ have units, we say that a functor (or
monad) $\mathbb{X}\to^F\mathbb{Y}$ is \textbf{normal} if it preserves
opcartesian cells with nullary source, which is to say it preserves
units in a coherent way.  Likewise, if $\mathbb{X}$ and $\mathbb{Y}$
have all units and composites (i.e.\ are pseudo double categories), we
say that $\mathbb{X}\to^F\mathbb{Y}$ is \textbf{strong} if it
preserves \emph{all} opcartesian cells.

\begin{example}\label{eg:span-strong}
  Any functor $\SPAN{C}\to\SPAN{D}$ induced by a pullback-preserving
  functor $\mb{C}\to\mb{D}$ is strong, and in particular normal.
\end{example}

\begin{example}\label{eg:mod-normal}
  It is also easy to see that $\MOD(F)$ is normal for any functor $F$,
  by the construction of units in \autoref{prop:mod-units}.
\end{example}

\begin{examples}\label{eg:fmc-strong}
  If \V\ is a cocomplete symmetric monoidal category in which $\otimes$
  preserves colimits on both sides, then $\MAT{V}$ has all composites,
  and the extension of the ``free monoid'' monad to $\MAT{V}$ from
  \autoref{eg:free-monoid-on-mat} is easily seen to be strong.  Since
  the ``free strict monoidal \V-category'' monad on $\PROF{V}$ is
  obtained by applying $\MOD$ to this, it is normal by our above
  observation.  In fact, it is also strong, essentially because the
  tensor product of reflexive coequalizers is again a reflexive
  coequalizer (see, for example,~\cite[A1.2.12]{ptj:elephant1}).
\end{examples}

\begin{examples}\label{eg:fssmc-strong}
  The ``free symmetric strict monoidal \V-category'' monad on $\PROF{V}$
  from \autoref{eg:free-ssmc} is also normal, essentially by definition,
  as are the ``free \V-category with strictly associative finite
  products'' monad from \autoref{eg:free-safp} and its relatives from
  \autoref{eg:free-coten}.  A more involved computation with
  coequalizers shows that the first is actually strong, and the second
  is strong whenever \V\ is cartesian monoidal.  However, it seems that
  the others are not in general strong.
\end{examples}

\begin{examples}
  The monads on $\MAT{V}$ constructed in \autoref{example:tv} are not
  generally strong or even normal.  Two notable exceptions are the
  ultrafilter monads on $\Rel$ and $\MAT{\eR}$, of which the first is
  strong and the second is normal.
\end{examples}

We write $\vDbln$ for the locally full sub-2-category of $\vDbl$
determined by the virtual double categories with units and normal
functors between them; thus $\MOD$ is a 2-functor $\vDbl\to\vDbln$.
In fact, we have the following.

\begin{proposition}\label{prop:mod-adjt}
  $\Mod$ is right pseudo-adjoint to the forgetful 2-functor $\vDbln
  \to\vDbl$.
\end{proposition}
\begin{proof}
  ``Pseudo-adjoint'' means that we have a pseudonatural $\eta$ and
  $\epsilon$ that satisfy the triangle identities up to coherent
  isomorphism.  We take $\epsilon_{\mathbb{X}}$ to be the forgetful functor
  $\Mod(\mathbb{X})\to\mathbb{X}$ which sends a monoid to its
  underlying object and a module to its underlying horizontal arrow;
  this is strictly 2-natural.  If $\mathbb{X}$ has units, we take
  $\eta_{\mathbb{X}}$ to be the ``unit assigning'' functor $\mathbb{X} \to
  \Mod(\mathbb{X})$ which sends $X$ to $U_X$ (which has a unique
  monoid structure) and $X\hto^p
  Y$ to itself considered as a $(U_X,U_Y)$-bimodule; this is only
  pseudonatural since normal functors preserve units only up to
  isomorphism.  But if we choose the units in $\Mod(\mathbb{X})$
  according to \autoref{prop:mod-units}, then the triangle identities
  are satisfied on the nose.
\end{proof}

In particular, if $1$ denotes the terminal virtual double category, then
the category of \emph{normal} functors $1\to \mathbb{X}$ is equivalent
to the vertical category of $\mathbb{X}$.  It then follows from
\autoref{prop:mod-adjt} that the category of arbitrary functors $1\to
\mathbb{X}$ is equivalent to the vertical category of
$\Mod(\mathbb{X})$.  Thus, \autoref{prop:mod-adjt} generalizes the
well-known observation (which dates back to~\cite{benabou:bicats}) that
monoids in a bicategory $\mathcal{B}$ are equivalent to lax functors
$1\to\mathcal{B}$.

\begin{remark}\label{rmk:mod-monad}
  It follows that $\Mod$ is a \emph{pseudomonad} on the 2-category
  $\vDbln$, and so in particular it has a multiplication
  \begin{equation}
    \Mod(\Mod(\mathbb{X}))\to \Mod(\mathbb{X}).\label{eq:mod-mult}
  \end{equation}
  Inspection reveals that an object of $\Mod(\Mod(\mathbb{X}))$
  consists of an object $X$ of $\mathbb{X}$, two monoids $X\hto^A X$
  and $X\hto^M X$, and a monoid homomorphism $A\to M$ whose vertical
  arrow components are identities.  The
  multiplication~\eqref{eq:mod-mult} simply forgets the monoid $A$.
  This idea will be further discussed in \cite{cs:laxidem}.
\end{remark}

\begin{remark}\label{rmk:horiz-bicat}
  If $\mathbb{X}$ is a virtual double category in which all units and
  composites exist (equivalently, it is a pseudo double category),
  then it has a \emph{horizontal bicategory} $\calH(\mathbb{X})$
  consisting of its objects, horizontal arrows, and cells of the
  form
  \[\xymatrix@-.5pc{ \ar[r]|-@{|}^-{} \ar@{=}[d] \ar@{}[dr]|{\Downarrow} &  \ar@{=}[d]\\
    \ar[r]|-@{|}_-{} & .}\]
  Clearly when $\mb{C}$ and $\mb{V}$ satisfy the required conditions
  for all composites to exist in our examples, we recover in this way
  the usual bicategories of matrices, spans, and enriched and internal
  profunctors.  Any functor between pseudo double categories likewise
  induces a lax functor of horizontal bicategories, but
  this is not true of transformations without additional assumptions;
  see Remarks \ref{rmk:oplax-htrans} and \ref{rmk:hstrong-htrans}.
\end{remark}

\section{$2$-categories of $T$-monoids}
\label{sec:2cat-tmon}

As proposed in the previous section, we now use the notion of units
introduced there to define 2-categories of generalized
multicategories, generalizing the approach taken
in~\cite[\S5.3]{leinster:higher-opds}.

\begin{proposition}\label{prop:2cat}
  Let $\mathbb{X}$ be a virtual double category in which all units
  exist.  Then it has a \textbf{vertical 2-category} $\calV\mathbb{X}$
  whose objects are those of $\mathbb{X}$, whose morphisms are the
  vertical arrows of $\mathbb{X}$, and whose 2-cells
  $\xymatrix{A\rtwocell^f_g{\alpha} &B}$ are the cells
  \[\xymatrix{A \ar[r]|-@{|}^{U_A}\ar[d]_g \ar@{}[dr]|{{\textstyle\Downarrow}\alpha}
    & A \ar[d]^f\\
    B\ar[r]|-@{|}_{U_B} & B}\]
  in $\mathbb{X}$.
\end{proposition}
\begin{proof}
  This is straightforward; note that when composing 2-cells we must
  use the isomorphisms $U_A\iso U_A\odot U_A$.
\end{proof}

In particular, for any $\mathbb{X}$, $\MOD(\mathbb{X})$ has a vertical
2-category, which we denote $\Mon(\mathbb{X})$ and call the
\emph{2-category of monoids in $\mathbb{X}$}.  (This 2-category is
closely related to various 2-categories of monads in a bicategory; see
\cite[\S2.3--2.4]{ls:ftm2}.)  It turns out that in this case, the
description of the 2-cells of $\Mon(\mathbb{X})$ can be rephrased in a
way that looks much more like a natural transformation.

\begin{proposition}\label{thm:mon-2cell}
  Giving a 2-cell $\xymatrix{A\rtwocell^f_g{\alpha} &B}$ in
  $\Mon(\mathbb{X})$ is equivalent to giving a cell
  \[\bfig\scalefactor{.7} \Atriangle[A_0`B_0`B_0;g_0`f_0`B]
  \morphism(0,0)|m|<1000,0>[B_0`B_0;|]
  \place(500,200)[\Downarrow\scriptstyle \alpha_0]
  \efig\]
  in $\mathbb{X}$ such that
  \begin{equation}
    \bfig\scalefactor{.8}
    \square(500,500)|arrb|[A_0`A_0`B_0`B_0;A`f_0`f_0`B]
    \place(750,750)[\Downarrow\scriptstyle f]
    \morphism(500,1000)<-500,-500>[A_0`B_0;g_0]
    \place(300,650)[{\scriptstyle \alpha_0}\Downarrow]
    \morphism(0,500)|b|<500,0>[B_0`B_0;B]
    \place(500,250)[\Downarrow\scriptstyle \bar{b}]
    \square/`=`=`->/<1000,500>[B_0`B_0`B_0`B_0;```B]
    \efig \;=\;
    \bfig\scalefactor{.8}
    \square(0,500)|allb|[A_0`A_0`B_0`B_0;A`g_0`g_0`B]
    \place(250,750)[{\scriptstyle g} \Downarrow]
    \morphism(500,1000)<500,-500>[A_0`B_0;f_0]
    \place(700,650)[\Downarrow\scriptstyle \alpha_0]
    \morphism(500,500)|b|<500,0>[B_0`B_0;B]
    \place(500,250)[\Downarrow\scriptstyle \bar{b}]
    \square/`=`=`->/<1000,500>[B_0`B_0`B_0`B_0.;```B]
    \efig\label{eq:simple-mon-2cell}
  \end{equation}
\end{proposition}
\begin{proof}
  This is an instance of \autoref{prop:mod-units}, where
  $m=n=0$, $f=g=1_A$, and $P=U_A$.
\end{proof}

\begin{example}
  Recall that in $\PROF{V} = \MOD(\MAT{V})$, the objects are
  \V-enriched categories and the vertical morphisms are \V-enriched
  functors; thus these are the objects and morphisms of
  $\Mon(\MAT{V})$.  Recalling from \autoref{eg:mat} the definition
  of cells in $\MAT{V}$, \autoref{thm:mon-2cell} implies
  that a 2-cell $\xymatrix{A\rtwocell^f_g{\alpha}&B}$ in $\Mon(\MAT{V})$ is
  given by a family of morphisms $I\to^{\alpha_x} B(fx,gx)$ for $x$
  an object of $A$, such that for every $x,y$ the following square
  commutes:
  \[\bfig\square<1200,400>[A(x,y)`B(fy,gy)\otimes
  B(fx,fy)`B(gx,gy)\otimes B(fx,gx)`B(fx,gy).;\alpha_y\otimes
  f`g\otimes \alpha_x``] \efig
  \]
  This is precisely the usual definition of a \V-enriched natural
  transformation; thus we have $\Mon(\MAT{V})\simeq \cat{V}$.
\end{example}

\begin{example}
  Likewise, $\Mon(\SPAN{C})\simeq \icat{C}$ is the 2-category of
  internal categories, functors, and natural transformations in \C.
\end{example}

\begin{examples}
  On the other hand, the vertical 2-category of $\SPAN{S}$ is just
  $\mb{S}$, regarded as a 2-category with only identity 2-cells.  The
  vertical 2-category of $\MAT{V}$ depends on what $\mb{V}$
  is, but usually it is not very interesting.  Thus, in general, vertical
  2-categories are only interesting for virtual double categories
  whose objects are ``category-like'' rather than ``set-like.''
\end{examples}

Now, if $T$ is a monad on a virtual double category $\mathbb{X}$, we
write $\KMon(\mathbb{X},T)$ for the 2-category
$\calV(\KMod(\mathbb{X},T))$ and call it \emph{the 2-category of
  $T$-monoids in $\mathbb{X}$}.  Its objects are $T$-monoids, its
morphisms are $T$-monoid homomorphisms, and its 2-cells may be called
\emph{$T$-monoid transformations}.  By \autoref{thm:mon-2cell}
and the definition of $\Hkl{X}{T}$, a $T$-monoid transformation
$\alpha\maps f\to g\maps A\to B$ is specified by a cell
\[\bfig\scalefactor{.7} \Atriangle[A_0`B_0`TB_0;g_0`Tf_0\circ \eta_{A_0} =
\eta_{TB_0} \circ f_0`B]
\morphism(0,0)|m|<1000,0>[B_0`TB_0;|]
\place(500,200)[\Downarrow\scriptstyle \alpha_0]
\efig
\]
such that
\[\bfig\scalefactor{.8}
\square(500,1000)|arrb|[A_0`TA_0`TA_0`T^2A_0;A`\eta`\eta`B]
\place(750,1250)[\Downarrow\scriptstyle \eta_A]
\square(500,500)|bmrb|[TA_0`T^2A_0`TB_0`T^2B_0;`Tf_0`T^2f_0`TB]
\place(750,750)[\Downarrow\scriptstyle Tf]
\morphism(500,1500)<-500,-1000>[A_0`B_0;g_0]
\place(300,800)[{\scriptstyle \alpha_0}\!\Downarrow]
\morphism(0,500)|b|<500,0>[B_0`TB_0;B]
\place(500,250)[\Downarrow\scriptstyle \bar{b}]
\square/`=`->`->/<1000,500>[B_0`T^2B_0`B_0`TB_0;``\mu_B`B]
\efig \;=\;
\bfig\scalefactor{.8}
\square(0,500)|almb|[A_0`TA_0`B_0`TB_0;A`g_0`Tg_0`B]
\place(250,750)[{\scriptstyle g} \Downarrow]
\morphism(500,1000)/{@{>}@/^5pt/}/<500,-500>[TA_0`T^2B_0;T(\eta_{B_0}) \circ Tf_0]
\place(700,650)[\Downarrow\scriptstyle T\alpha_0]
\morphism(500,500)|b|<500,0>[TB_0`T^2B_0;TB]
\place(500,250)[\Downarrow\scriptstyle \bar{b}]
\square/`=`->`->/<1000,500>[B_0`T^2B_0`B_0`TB_0.;``\mu_B`B]
\efig
\]
Many authors have defined this 2-category $\KMon(\mathbb{X},T)$ in
seemingly ad-hoc ways, whereas it falls quite naturally out of the
framework of virtual double categories.  This was also observed
in~\cite[\S5.3]{leinster:higher-opds}; see \S\ref{section:cartesian}.

\begin{example}
  Let $T$ be the ``free monoid'' monad on $\MAT{Set}$, so that
  $T$-monoids are ordinary multicategories.  If $A\two^f_g B$ are
  functors, then according to the above, a transformation
  $f\to^{\alpha} g$ consists of, for each $x\in A$, a morphism
  $\alpha_x \in B(\eta(fx),gx)$ (that is to say, a morphism
  $(fx)\to^{\alpha_x} gx$ with source of length one) such that for
  any morphism $\xi\maps (x_1,\dots, x_n) \to y$ in $A$, we have
  \[\alpha_y \circ f(\xi) = g(\xi) \circ (\alpha_{x_1},\dots, \alpha_{x_n}).
  \]
  This is the usual notion of transformation for functors between
  multicategories.
\end{example}

\begin{example}
  When $T$ is the ``free category'' monad on directed graphs, so that
  $T$-monoids are virtual double categories, $T$-monoid
  transformations are the same as the transformations we defined in
  \ref{def:transforms}.
\end{example}

\begin{example}
  Let $U$ be the ultrafilter monad on $\Rel$, so that $U$-monoids are
  topological spaces.  If $A\two^f_g B$ are continuous maps (i.e.\
  $U$-monoid homomorphisms), then there exists a transformation $f\to
  g$ (which is necessarily unique) just when for all $x\in A$, the
  principal ultrafilter $\eta_{fx}$ converges to $gx$ in $B$.  This is
  equivalent to saying that $f\ge g$ in the pointwise ordering induced
  by the usual \emph{specialization order} on $B$.  The situation for
  other topological examples is similar.
\end{example}

Any \emph{normal} functor clearly induces a strict 2-functor between
vertical 2-categories.  In fact, if $\twoCat$ denotes the 2-category of
2-categories, strict 2-functors, and strict 2-natural transformations,
then we have:

\begin{proposition}
  There is a strict 2-functor $\calV\maps \vDbln \to \twoCat$.
\end{proposition}

In particular, any normal monad $T$ on $\mathbb{X}$ induces a strict 2-monad on
$\calV(\mathbb{X})$.  As we saw in \S\ref{sec:2cat-vdc}, most monads
on virtual double categories are ``extensions'' of a known monad on
their vertical categories (or vertical 2-categories), so this
construction usually just recovers the familiar monad we started with.
In \S\ref{sec:vert-alg}, we will show that $\calV(T)$-algebras are
closely related to $T$-monoids.

\section{Virtual equipments}
\label{sec:pfrbi}

If we have succeeded in convincing the reader that virtual double
categories are inevitable, she may be justified in wondering why they
have not been more studied.  Certainly, \emph{virtual} double
categories involve significant complexity above and beyond
\emph{pseudo} double categories, and the latter suffice to describe
the important examples $\SPAN{C}$, $\MAT{V}$, $\PROF{V}$, and
$\iPROF{C}$ as long as \V\ and \C\ are suitably cocomplete.  However,
even pseudo double categories have generally received less
publicity than bicategories.

One possible reason for this is that in most of the (pseudo or virtual) double
categories arising in practice, the vertical arrows are more tightly
coupled to the horizontal arrows that we have heretofore accounted for;
in fact they can almost be \emph{identified with} certain horizontal
arrows.  For example, a \V-functor $A\to^f B$ is determined, up to
isomorphism, by the \V-profunctor $A\hto^{B(1,f)} B$ defined by
$B(1,f)(b,a) = B(b,fa)$.  Furthermore, this coupling is very important
for many applications, such as the formal definition of weighted
limits and colimits (see \cite{street:enr-cohom,wood:proarrows-i}), so
a mere double category (pseudo or virtual) would be insufficient for
these purposes.  Because of this, many authors have been content to
work with bicategories, or bicategories with a collection of horizontal
arrows singled out (such as the ``proarrow equipments'' of~\cite{wood:proarrows-i}).

However, while not all pseudo double categories exhibit this type of
coupling, it is possible to characterize those that do (and they turn
out to be equivalent to the ``proarrow equipments'' mentioned above).  The
basic idea of this dates back to~\cite{bs:dblgpd-xedmod}, but it has
recently been revived in various equivalent forms; see for instance
\cite{DomThesis,gp:double-limits,gp:double-adjoints,dpp:spans,FramedBicats}
and the Notes at the end of~\cite[Ch.~5]{leinster:higher-opds}.
Since this type of coupling also plays an important role in the theory
of generalized multicategories, in this section we give the basic
definitions appropriate to the virtual case.

The basic idea is the following.  The profunctor $B(1,f)$ considered
above can be constructed in two stages: first we consider the
hom-profunctor $B(-,-)\maps B\hto B$, and then we \emph{precompose} it
with $f$ on one side.  We already know from
\S\ref{sec:composites-units} that hom-profunctors are the \emph{units}
in $\PROF{V}$, so it remains only to characterize \emph{precomposition
  with functors} in terms of $\PROF{V}$.  This is accomplished by the
following definition.

\begin{definition}
  A cell
  \begin{equation}
    \xymatrix{ \ar[r]|-@{|}^-{p}\ar[d]_f \ar@{}[dr]|{{\textstyle\Downarrow} \mathrm{cart}} &
      \ar[d]^g\\
      \ar[r]|-@{|}_q & }\label{eq:cart-2cell}
  \end{equation}
  in a virtual double category is \textbf{cartesian} if any cell
  \[\xymatrix{ \ar[r]|-@{|}^-{r_1}\ar[d]_{fh} \ar@{}[drrr]|{{\textstyle\Downarrow}} &
    \ar[r]|-@{|}^-{r_2} & \dots \ar[r]|-@{|}^-{r_n} &
    \ar[d]^-{gk}\\
    \ar[rrr]|-@{|}_q &&& }\]
  factors through it uniquely as follows:
  \[\xymatrix{ \ar[r]|-@{|}^-{r_1}\ar[d]_{h} \ar@{}[drrr]|{{\textstyle\Downarrow}} &
    \ar[r]|-@{|}^-{r_2} & \dots \ar[r]|-@{|}^-{r_n} &
    \ar[d]^-{k}\\
    \ar[rrr]|p \ar[d]_f \ar@{}[drrr]|{{\textstyle\Downarrow} \mathrm{cart}} &&&
    \ar[d]^g\\
    \ar[rrr]|-@{|}_q &&& .}\]
\end{definition}

If there exists a cartesian cell~\eqref{eq:cart-2cell}, we say that
the $p$ is the \textbf{restriction $q(g,f)$}.  The notation is
intended to suggest precomposition of a profunctor $q(-,-)$ with $f$
and $g$.  When $f$ or $g$ is an identity, we write $q(g,1)$ or
$q(1,f)$, respectively.  It is evident from the universal property
that restrictions are unique up to isomorphism, and pseudofunctorial;
that is, we have $q(1,1)\iso q$ and $q(1,g)(1,f) \iso q(1,g f)$
coherently.

We say that $\mathbb{X}$ \textbf{has restrictions} if $q(g,f)$ exists
for all $q$, $f$, and $g$, and that a functor \textbf{preserves
restrictions} if it takes cartesian cells to cartesian cells.  We
write $\vDblr$ for the sub-2-category of $\vDbl$ determined by the
virtual double categories with restrictions and the
restriction-preserving functors.

\begin{examples}
  The virtual double categories $\MAT{V}$, $\PROF{V}$, $\SPAN{C}$, and
  $\iPROF{C}$ have all restrictions.  Restrictions in $\MAT{V}$ are
  given by reindexing matrices, restrictions in $\SPAN{C}$ are given
  by pullback, and restrictions in $\PROF{V}$ and $\iPROF{C}$ are given
  by precomposing with functors.
\end{examples}

Note that the restrictions in $\PROF{V}$ and $\iPROF{C}$ are induced by
those in $\MAT{V}$ and $\SPAN{C}$, in the following general way.

\begin{proposition}\label{prop:mod-restr}
  If $\mathbb{X}$ is a virtual double category with restrictions, then
  $\Mod(\mathbb{X})$ also has restrictions.
\end{proposition}
\begin{proof}
  If $B\hto^p D$ is a bimodule in $\mathbb{X}$ and
  \[\vcenter{\xymatrix{A_0 \ar[r]|-@{|}^-{A} \ar[d]_{f_0} \ar@{}[dr]|{{\textstyle\Downarrow}f} & A_0 \ar[d]^{f_0}\\
      B_0\ar[r]|-@{|}_-{B} & B_0}}
  \qquad\text{and}\qquad
  \vcenter{\xymatrix{C_0 \ar[r]|-@{|}^-{C} \ar[d]_{g_0} \ar@{}[dr]|{{\textstyle\Downarrow}g} & C_0 \ar[d]^{g_0}\\
    D_0\ar[r]|-@{|}_-{D} & D_0}}
  \]
  are monoid homomorphisms, then the restriction $p(g_0,f_0)$ in
  $\mathbb{X}$ becomes an $(A,C)$-bimodule in an obvious way,
  making it into the restriction $p(g,f)$ in
  $\Mod(\mathbb{X})$.
\end{proof}

The other ingredient in the construction of generalized
multicategories also preserves restrictions.

\begin{proposition}\label{prop:hkl-restr}
  If $\mathbb{X}$ is a virtual double category with restrictions and
  $T$ is a monad on $\mathbb{X}$, then $\Hkl{X}{T}$ also has
  restrictions.
\end{proposition}
\begin{proof}
  Let $A\hto^p T B$ be a horizontal arrow in $\mathbb{X}$,
  regarded as a horizontal arrow $A\hto B$ in $\Hkl{X}{T}$.  It is easy to
  verify that the restriction of $p$ along $C\to^f A$ and $D\to^g B$
  in $\Hkl{X}{T}$ is given by the restriction $p(T g, f)$ in
  $\mathbb{X}$.
\end{proof}

Therefore, if $\mathbb{X}$ has restrictions, so does
$\KMod(\mathbb{X},T)$ for any monad $T$ on $\mathbb{X}$.  Moreover, by
\autoref{prop:mod-units}, $\KMod(\mathbb{X},T)$ always also has units.
As suggested in the introduction to this section, units and restrictions
together are an especially important combination, so we give a special
name to this situation.

\begin{definition}
  A \textbf{virtual equipment} is a virtual double category in which
  all units and all restrictions exist.
\end{definition}

We write \vEquip\ for the locally full sub-2-category of $\vDbl$
determined by the virtual equipments and the normal
restriction-preserving functors between them (however, see
\autoref{thm:func-pres-restr}).  We can now observe that $\Mod$ is a
2-functor from $\vDblr$ to $\vEquip$.

\begin{examples}
  $\SPAN{C}$, $\PROF{V}$, and $\iPROF{C}$ are always virtual
  equipments, and $\MAT{V}$ is a virtual equipment whenever \V\ has an
  initial object preserved by $\otimes$.  More generally,
  $\Mod(\mathbb{X})$ and $\KMod(\mathbb{X},T)$ are virtual equipments
  whenever $\mathbb{X}$ has restrictions.
\end{examples}

If $A\to^f B$ is a vertical arrow in a virtual equipment,
we define its \textbf{base change objects} to be
\[B(1,f) = U_B(1,f) \qquad\text{and}\qquad B(f,1) = U_B(f,1).\]
These come with cartesian cells
\begin{equation}
  \begin{array}{c}
    \xymatrix{
      \ar[r]|-@{|}^-{B(1,f)} \ar[d]_f \ar@{}[dr]|{\twoar(0,-1)}
      & \ar@{=}[d]\\
      \ar[r]|-@{|}_-{U_B} & }
  \end{array}\quad\text{and}\quad
  \begin{array}{c}
    \xymatrix{
      \ar[r]|-@{|}^-{B(f,1)} \ar@{=}[d] \ar@{}[dr]|{\twoar(0,-1)}
      & \ar[d]^f\\
      \ar[r]|-@{|}_-{U_B} & }
  \end{array}.\label{eq:bcobj-cart}
\end{equation}
By factoring $U_f$ through these cartesian cells, we obtain two
further cells
\begin{equation}
  \begin{array}{c}
    \xymatrix{
      \ar[r]|-@{|}^-{U_A} \ar[d]_f \ar@{}[dr]|{\twoar(0,-1)}
      & \ar@{=}[d]\\
      \ar[r]|-@{|}_-{B(f,1)} & }
  \end{array}\quad\text{and}\quad
  \begin{array}{c}
    \xymatrix{
      \ar[r]|-@{|}^-{U_A} \ar@{=}[d] \ar@{}[dr]|{\twoar(0,-1)}
      & \ar[d]^f\\
      \ar[r]|-@{|}_-{B(1,f)} & }
  \end{array}\label{eq:bcobj-opcart}
\end{equation}
such that the following equations hold:
\begin{align}
  \begin{array}{c}
    \xymatrix@C=3.5pc{
      \ar[r]|-@{|}^-{U_A} \ar@{=}[d] \ar@{}[dr]|{\twoar(0,-1)}
      & \ar[d]^f\\
      \ar[r]|{B(1,f)} \ar[d]_f \ar@{}[dr]|{\twoar(0,-1)}
      & \ar@{=}[d]\\
      \ar[r]|-@{|}_-{U_B} & }
  \end{array} &= 
  \begin{array}{c}
    \xymatrix{ \ar[r]|-@{|}^-{U_A} \ar[d]_f
      \ar@{}[dr]|{{\textstyle\Downarrow} U_f} &  \ar[d]^f\\
      \ar[r]|-@{|}_-{U_B} & }
  \end{array}
  &
  \begin{array}{c}
    \xymatrix@C=3.5pc{
      \ar[r]|-@{|}^-{U_A} \ar[d]_f \ar@{}[dr]|{\twoar(0,-1)}
      & \ar@{=}[d]\\
      \ar[r]|{B(f,1)} \ar@{=}[d] \ar@{}[dr]|{\twoar(0,-1)}
      & \ar[d]^f\\
      \ar[r]|-@{|}_-{U_B} & }
  \end{array} &= 
  \begin{array}{c}
    \xymatrix{ \ar[r]|-@{|}^-{U_A} \ar[d]_f
      \ar@{}[dr]|{{\textstyle\Downarrow} U_f} &  \ar[d]^f\\
      \ar[r]|-@{|}_-{U_B} & }
  \end{array}\label{eq:frbi-compconj-1}
\end{align}
Moreover, the following equations also hold:
\begin{align}
  \begin{array}{c}
    \xymatrix@C=3.5pc{
      \ar[r]|-@{|}^-{U_A} \ar@{=}[d] \ar@{}[dr]|{\twoar(0,-1)} &
      \ar[r]|-@{|}^-{B(1,f)} \ar[d]|f \ar@{}[dr]|{\twoar(0,-1)}
      & \ar@{=}[d]\\
      \ar[r]|-{B(1,f)} &
      \ar[r]|-{U_B} &\\
      \ar[rr]|-@{|}_-{B(1,f)} \ar@{}[urr]|{{\textstyle\Downarrow}\mathrm{opcart}} \ar@{=}[u] &&
      \ar@{=}[u]}
  \end{array} &=
  \vcenter{\xymatrix{
      \ar[r]|-@{|}^{U_A}\ar@{=}[d] \ar@{}[drr]|{{\textstyle\Downarrow}\mathrm{opcart}} &
      \ar[r]|-@{|}^{B(1,f)} &
      \ar@{=}[d]\\
      \ar[rr]|-@{|}_{B(1,f)} &&
    }}
  \label{eq:frbi-compconj-2}
  \\
  \begin{array}{c}
    \xymatrix@C=3.5pc{
      \ar[r]|-@{|}^-{B(f,1)} \ar@{=}[d] \ar@{}[dr]|{\twoar(0,-1)} &
      \ar[r]|-@{|}^-{U_A} \ar[d]|f \ar@{}[dr]|{\twoar(0,-1)}
      & \ar@{=}[d]\\
      \ar[r]|-{U_B} &
      \ar[r]|-{B(f,1)} &\\
      \ar[rr]|-@{|}_-{B(f,1)} \ar@{}[urr]|{{\textstyle\Downarrow}\mathrm{opcart}} \ar@{=}[u] &&
      \ar@{=}[u]}
  \end{array} &=
  \vcenter{\xymatrix{
      \ar[r]|-@{|}^{B(f,1)}\ar@{=}[d] \ar@{}[drr]|{{\textstyle\Downarrow}\mathrm{opcart}} &
      \ar[r]|-@{|}^{U_A} &
      \ar@{=}[d]\\
      \ar[rr]|-@{|}_{B(f,1)} &&
    }}
  \label{eq:frbi-compconj-3}
\end{align}
We can verify these by postcomposing each side with the appropriate
cartesian cell, using the equations~\eqref{eq:frbi-compconj-1}, and
invoking the uniqueness of factorizations through cartesian cells.
In the terminology of~\cite{dpp:spans},
equations~\eqref{eq:frbi-compconj-1}--\eqref{eq:frbi-compconj-3}
are said to make $B(1,f)$ and $B(f,1)$ into a \textbf{companion} and a
\textbf{conjoint} of $f$, respectively.

\begin{example}
  In $\SPAN{C}$, the base change objects $B(1,f)$ and $B(f,1)$ are the
  spans
\[
          	\bfig
          	\morphism<-225,-250>[A`A;1_A]
          	\morphism<225,-250>[A`B;f]
          	
          	\place(525,-125)[\mbox{and}]
          	
          	\morphism(1000,0)<-225,-250>[A`B;f]
          	\morphism(1000,0)<225,-250>[A`A;1_A]
          	\efig
\]
respectively.  These are often called the
  \emph{graph} of $f$.
\end{example}

\begin{example}
  In $\MAT{V}$, for a function $f\maps X\to Y$  the base change object
  $Y(1,f)$ is the matrix
  \[Y(1,f)(x,y) =
  \begin{cases}
    I & \text{if } f(x)=y \\
    \emptyset & \text{otherwise}.
  \end{cases}
  \]
\end{example}

\begin{example}
  In $\PROF{V}$, the base change objects $B(1,f)$ and $B(f,1)$ are the
  representable distributors defined by $B(1,f)(b,a) = B(b,fa)$ and $B(f,1)(a,b) =
  B(fa,b)$.  Base change objects in $\iPROF{C}$ are analogous.
\end{example}

At first glance, base change objects may seem only to be a particular
special case of restrictions.  However, it turns out that \emph{all}
restrictions can be recovered by composition with the base change
objects (hence the name).

\begin{theorem}\label{thm:comp-bc-gives-restr}
  Let $B\hto^p D$ be a horizontal arrow and $f\maps A\to B$ and $g\maps
  C\to D$ be vertical arrows in a virtual equipment.  Then $B(1,f) \odot p\odot
  D(g,1)$ exists and is isomorphic to $p(g,f)$.
\end{theorem}
\begin{proof}
  Consider the composite
  \begin{equation}
    \xymatrix{
      \ar[r]|-@{|}^{B(1,f)} \ar[d]_f \ar@{}[dr]|{{\twoar(0,-1)}} &
      \ar[r]|-@{|}^p \ar@{=}[d] \ar@{}[dr]|{1_p} &
      \ar@{=}[d] \ar[r]|-@{|}^{D_g} \ar@{}[dr]|{{\twoar(0,-1)}} &
      \ar[d]^g\\
      \ar[r]|-@{|}_{U_B} & \ar[r]|p & \ar[r]|-@{|}_{U_D} & \\
      \ar@{=}[u] \ar@{}[urrr]|{{\textstyle\Downarrow}\mathrm{opcart}} \ar[rrr]|-@{|}_p &&&  \ar@{=}[u] .}\label{eq:frbi-compconj-cart}
  \end{equation}
  By the universal property of restriction, this factors through the
  cartesian cell defining $p(g,f)$ to give a canonical cell
  \begin{equation}
    \xymatrix{
      \ar[r]|-@{|}^{B(1,f)} \ar@{=}[d] \ar@{}[drrr]|{{\twoar(0,-1)}} &
      \ar[r]|-@{|}^p &
      \ar[r]|-@{|}^{D(g,1)}  &
      \ar@{=}[d]\\
      \ar[rrr]|-@{|}_{p(g,f)} &&&.}\label{eq:frbi-compconj-opcart}
  \end{equation}
  We claim that this cell is opcartesian.  To show this, suppose
  given a cell
  \[\xymatrix{
    \ar@{-->}[r] \ar[d]_h &
    \ar[r]|-@{|}^{B(1,f)}  \ar@{}[drrr]|{{\twoar(0,-1)}} &
    \ar[r]|-@{|}^p & \ar[r]|-@{|}^{D(g,1)}  &
    \ar@{-->}[r] &
    \ar[d]^k\\
    \ar[rrrrr]|-@{|}_{q} &&&&&.}\]
  We need to factor it uniquely
  through~\eqref{eq:frbi-compconj-opcart}.  A factorization is given
  by the composite
  \[\xymatrix@C=3.5pc{
    & \ar@{-->}[r] \ar@{=}[dl] &
    \ar[r]|-@{|}^{p(g,f)}\ar@{=}[dl] \ar@{=}[d] &
    \ar@{-->}[r] \ar@{=}[d] \ar@{=}[dr]  &
    \ar@{=}[dr]  &\\
    \ar@{-->}[r] \ar@{=}[d] &
    \ar@{=}[d]  \ar[r]|{U_A} \ar@{}[dr]|{\twoar(0,-1)} &
    \ar[d]|f \ar[r]|-@{|}^{p(g,f)} \ar@{}[dr]|{\twoar(0,-1)} &
    \ar[d]|g \ar[r]|-{U_C} \ar@{}[dr]|{\twoar(0,-1)} &
    \ar@{=}[d] \ar@{-->}[r] &
    \ar@{=}[d]\\
    \ar[d]_h \ar@{-->}[r] &
    \ar[r]|{B(1,f)} \ar@{}[drrr]|{{\twoar(0,-1)}} &
    \ar[r]|p &
    \ar[r]|{D(g,1)}  &
    \ar@{-->}[r] &
    \ar[d]^k\\
    \ar[rrrrr]|-@{|}_{q} &&&&&.}\]
  To verify that this is a factorization, and that it is unique, we
  use the
  equations~\eqref{eq:frbi-compconj-1}--\eqref{eq:frbi-compconj-3}.  The
  details are similar to the
  proof of~\cite[Theorem 4.1]{FramedBicats}.
\end{proof}

\begin{corollary}\label{thm:bcobj-psfr}
  For vertical arrows $f\maps A\to B$ and $g\maps B\to C$, the
  composite $B(1,f) \odot C(1,g)$ always exists and is isomorphic to
  $C(1,g f)$, and dually.\qed
\end{corollary}

We also have the following dual result.

\begin{theorem}\label{thm:comp-bc-gives-extn}
  For arrows $f\maps A\to B$ and $g\maps C\to D$ in a virtual
  equipment, we have a bijection between cells of the form
  \[\vcenter{\xymatrix{A \ar[r]|-@{|}^-{p} \ar[d]_f \ar@{}[dr]|{\twoar(0,-1)} & C \ar[d]^g\\
      B\ar[r]|-@{|}_-{q} & D}}
  \qquad\text{and}\qquad
  \vcenter{\xymatrix{B \ar[r]|-@{|}^-{B(f,1)} \ar@{=}[d] \ar@{}[drrr]|{\twoar(0,-1)} &
    A \ar[r]|-@{|}^-{p} &
    C \ar[r]|-@{|}^-{D(1,g)} & D \ar@{=}[d]\\
    B\ar[rrr]|-@{|}_-{q} &&& D}}\]
\end{theorem}
\begin{proof}
  The inverse bijections are given by composing with the
  cells~\eqref{eq:bcobj-cart} and~\eqref{eq:bcobj-opcart}.  (Recall
  that all composites with units exist in any virtual double
  category.)  The fact that they are inverses follows
  from~\eqref{eq:frbi-compconj-1}--\eqref{eq:frbi-compconj-3}.
\end{proof}

It follows that in the situation of \autoref{thm:comp-bc-gives-extn}, if
the composite $B(f,1) \odot p\odot D(1,g)$ exists, then it is a
``corestriction'' or ``extension'' of $p$ along $f$ and $g$---that is,
it satisfies a universal property dual to that of a restriction.


Combining Theorems \ref{thm:comp-bc-gives-restr} and
\ref{thm:comp-bc-gives-extn}, we obtain the following.

\begin{corollary}\label{thm:companion-mates}
  In a virtual equipment, there is a bijection between cells of the
  form
  \[\vcenter{\xymatrix{
      \ar[r]|-@{|}^{p}\ar[d]_h \ar@{}[drr]|{{\twoar(0,-1)}}&
      \ar[r]|-@{|}^{B(1,f)}  &
      \ar[d]^k\\
      \ar[r]|-@{|}_{D(1,g)} &
      \ar[r]|-@{|}_{q} &
    }}\quad\text{and}\quad
  \vcenter{\xymatrix{
      \ar[r]|-@{|}^{p}\ar[d]_h \ar@{}[ddr]|{{\twoar(0,-1)}} &
      \ar[d]^f\\
      \ar[d]_g &
      \ar[d]^k\\
      \ar[r]|-@{|}_{q} &
    }}
  \]
\end{corollary}

Taking $p$ and $q$ to be units and $h$ and $k$ to be identities, we
obtain:

\begin{corollary}\label{thm:bco-lff}
  For vertical arrows $f,g\maps A\to B$ in a virtual equipment
  $\mathbb{X}$, there is a bijection between cells $f\to g$ in
  $\calV\mathbb{X}$ and cells
  \[\xymatrix@-.5pc{ \ar[r]|-@{|}^-{B(1,f)} \ar@{=}[d] \ar@{}[dr]|{\Downarrow} &  \ar@{=}[d]\\
    \ar[r]|-@{|}_-{B(1,g)} & ,}\]
  which respects
  composition.  Similarly, we have a bijection between cells $f\to
  g$ and cells $B(g,1) \to B(f,1)$.\qed
\end{corollary}

Now suppose that $\mathbb{X}$ is a virtual equipment which moreover
has all composites.  Then it has a horizontal bicategory
$\calH\mathbb{X}$, and Corollaries \ref{thm:bcobj-psfr} and
\ref{thm:bco-lff} imply that $f\mapsto B(1,f)$ defines a pseudofunctor
$\calV\mathbb{X}\to\calH\mathbb{X}$ which is locally full and
faithful.  Furthermore, it is easy to verify that $B(f,1)$ is right
adjoint to $B(1,f)$ in $\calH\mathbb{X}$.

This structure---a pseudofunctor which is bijective on objects,
locally full and faithful, and which takes each 1-cell to one having a
right adjoint---was defined in~\cite{wood:proarrows-i} to constitute
an \textbf{equipment}.  This is the structure we referred to in the
introduction to this section, which many authors have used where we
would find double categories more natural.  In fact, it is not hard to
show (see \cite[Appendix C]{FramedBicats}) that an equipment in the
sense of \cite{wood:proarrows-i} is \emph{equivalent} to a virtual
equipment which has all composites (this was called a \emph{framed
  bicategory} in \cite{FramedBicats}).  Therefore, from now on we use
\emph{equipment} to mean a virtual equipment having all composites
(thereby justifying the terminology ``virtual equipment'').

\begin{remark}
  It was shown in~\cite{FramedBicats} that an equipment can equally
  well be defined as a pseudo double category in which all
  restrictions exist, or in which all ``extensions'' exist (in the sense
  mentioned after \autoref{thm:comp-bc-gives-extn}), or in which
  there exist base change objects with cells~\eqref{eq:bcobj-cart}
  and~\eqref{eq:bcobj-opcart}
  satisfying~\eqref{eq:frbi-compconj-1}--\eqref{eq:frbi-compconj-3}.
  In the virtual case, we have a
  Goldilocks trifurcation: merely having base change objects is too
  weak, and having all extensions is too strong, but having
  restrictions (together with units) is just right.
\end{remark}

We now consider how functors and transformations interact with
restrictions.

\begin{theorem}\label{thm:func-pres-restr}
  Any functor $F$ between virtual equipments preserves restriction.
\end{theorem}
\begin{proof}
  The proof given for equipments in~\cite[Theorem 6.4]{FramedBicats}
  applies basically verbatim to virtual equipments.
\end{proof}

In particular, \vEquip\ is in fact a \emph{full} sub-2-category of \vDbln.
Note, though, that an arbitrary functor $F$ between virtual equipments
still may not preserve units, so that while we
have $F(B(1,f)) \iso FU_B(1,Ff)$, neither need be the same as
$FB(1,Ff)$.  Of course, they are the same if $F$ is normal.

Now, recall that any transformation $\xymatrix{\mathbb{X} \rtwocell^F_G{\alpha}&\mathbb{Y}}$
of functors between virtual equipments induces a
strictly 2-natural transformation $\calV(\alpha)$ of 2-functors
between vertical 2-categories.  In particular, we have $\alpha_B \circ
F(f) = G(f) \circ \alpha_A$ for any vertical arrow $f\maps A\to B$ in
$\mathbb{X}$.  However, we also have the cell component
\[\vcenter{\xymatrix@C=3pc{ \ar[r]|-@{|}^-{F(B(1,f))} \ar[d]_{\alpha_A}
    \ar@{}[dr]|{{\textstyle\Downarrow}\alpha_{B(1,f)}} &  \ar[d]^{\alpha_B}\\
    \ar[r]|-@{|}_-{G(B(1,f))} & }}
\]
of $\alpha$.  If $F$ and $G$ are normal, so that $F(B(1,f))\iso
FB(1,{Ff})$ and $G(B(1,f))\iso GB(1,{Gf})$, then by
\autoref{thm:companion-mates}, $\alpha_{B(1,f)}$ induces a 2-cell $\alpha_B
\circ F(f) \to G(f) \circ \alpha_A$, which seems to be trying to make
$\calV(\alpha)$ into an \emph{oplax} natural transformation.
Fortunately, however, this is an illusion.

\begin{proposition}\label{thm:vert-oplax-is-id}
  In the above situation, the 2-cell $\alpha_B \circ F(f) \to G(f)
  \circ \alpha_A$ induced by $\alpha_{B(1,f)}$ is an identity.
\end{proposition}
\begin{proof}
  This follows by inspection of how this 2-cell is constructed, and
  use of the cell naturality of $\alpha$.
\end{proof}

\begin{remark}\label{rmk:oplax-htrans}
  Recall from \autoref{rmk:horiz-bicat} that any functor
  $\mathbb{X}\to^F \mathbb{Y}$ between pseudo double categories induces
  a lax functor $\calH(\mathbb{X})\to^{\calH(F)} \calH(\mathbb{Y})$
  between horizontal bicategories, but not every transformation $F\to G$
  induces a transformation $\calH(F)\to\calH(G)$.
  It is true, however, that if $\mathbb{X}$ and $\mathbb{Y}$ are
  \emph{equipments}, then any transformation $F\to^\alpha G$
  induces an \emph{oplax} transformation
  $\calH(F)\to^{\calH(\alpha)}\calH(G)$ whose component at $X$ is
  $GX(1,{\alpha_X})$.  Likewise, the components $(GX)({\alpha_X},1)$
  form a lax transformation $\calH(G)\to\calH(F)$.  See also
  \autoref{rmk:hstrong-htrans}.
\end{remark}


\section{Normalization}
\label{sec:normalization}

With the notion of virtual equipment under our belt, we now return to
the general theory of generalized multicategories.
We observed in \S\ref{sec:genmulti} that for virtual double categories
whose objects are ``category-like,'' such as $\PROF{V}$ and $\iPROF{C}$
(as opposed to those such as $\MAT{V}$ and $\SPAN{C}$, whose objects
are ``set-like''), general $T$-monoids often contain too much structure.
For instance, if $T$ is the ``free strict monoidal category'' monad on
$\PROF{\Set}$, then a $T$-monoid consists of a category $A$, a
multicategory $M$ and a bijective-on-objects functor from $A$ to the
underlying category of $M$.  Usually, the morphisms of $A$ constitute
superfluous data which we would like to eliminate.  (This is not
always true, though: in \cite{cheng:opetopic} these extra
morphisms played an important role.)

The obvious way to eliminate this extra data, which we adopted in
describing examples of this sort in \S\ref{sec:genmulti}, is to
require $A$ to be a discrete category; this way the extra morphisms
simply do not exist.  However, a different way to eliminate it is to
require the given functor $A\to M$ to induce an \emph{isomorphism}
between $A$ and the underlying category of $M$; this way the extra
morphisms exist, but are determined uniquely by the rest of the
structure.  In this section we define general analogues of both
approaches, show their equivalence under general hypotheses, and argue
that when they are not equivalent it is usually the \emph{second}
approach that is more useful.  (This second approach was also the one
taken in~\cite{HermidaCoherent}.)

\begin{definition}
  Let $\mathbb{X}$ be a virtual equipment and let $T$ be a monad on
  $\Mod(\mathbb{X})$.  A $T$-monoid $A\hto^M TA$ is called
  \textbf{object-discrete} if $A$ is a monoid in $\mathbb{X}$ of the
  form $U_X$.
\end{definition}

We write $\dMod(\Mod(\mathbb{X}),T)$ for the full sub-virtual-equipment of
$\KMod(\Mod(\mathbb{X}),T)$ determined by the object-discrete $T$-monoids, and
$\dKMon(\mathbb{X},T)$ for its vertical 2-category.  Note that
object-discreteness is only defined for a monad on a virtual equipment
of the form $\Mod(\mathbb{X})$.

\begin{definition}
  Let $T$ be a monad on a virtual equipment $\mathbb{X}$.  A
  $T$-monoid $A\hto^M TA$ is \textbf{normalized} if its unit cell
  \[\xymatrix{A \ar[r]|-@{|}^-{U_A} \ar@{=}[d] \ar@{}[dr]|{\twoar(0,-1)} & A \ar[d]^{\eta}\\
    A\ar[r]|-@{|}_-{M} & TA}\]
  is cartesian in $\mathbb{X}$.
\end{definition}

We write $\nMod(\mathbb{X},T)$ for the full sub-virtual-equipment of
$\KMod(\mathbb{X},T)$ determined by the normalized $T$-monoids, and
$\nKMon(\mathbb{X},T)$ for its vertical 2-category.  Unlike
object-discreteness, normalization is defined for monads on any
virtual equipment.

Now, to prove an equivalence between normalization and
object-discreteness, we need to introduce the following definitions.

\begin{definition}
  A monoid homomorphism
  \[\bfig\scalefactor{.8}
  \square/{@{>}|-*@{|}}`>`>`{@{>}|-*@{|}}/[A_0`A_0`B_0`B_0;A`f_0`f_0`f]
  \place(225,300)[\twoar(0,-1)]
  \place(300,250)[\scriptstyle f]
  \efig
  \]
  in a virtual double category $\mathbb{X}$ is called
  \textbf{bijective on objects} (or \textbf{b.o.})\ if $f_0$ is an
  isomorphism.  It is called \textbf{fully faithful} (or
  \textbf{f.f.})\ if the cell $f$ is cartesian.
\end{definition}

\begin{lemma}\label{lem:monoid-restr}
  If $A_0\hto^A A_0$ is a monoid in a virtual double category
  $\mathbb{X}$ with restrictions and $X\to^f A_0$ is any vertical
  arrow, then $X\hto^{A(f,f)} X$ is also a monoid, and its defining
  cartesian cell is a monoid homomorphism with $f$ as its vertical
  part.
\end{lemma}
\begin{proof}
  To obtain a multiplication for $A(f,f)$, we compose two copies of
  the defining cartesian cell with the multiplication of $A$, then
  factor the result through the defining cartesian cell.  The unit is
  similar.
\end{proof}

\begin{lemma}\label{lem:bo-ff}
  If $\mathbb{X}$ has restrictions, then (b.o., f.f.) is a
  factorization system on the category $\mathbf{Mon}(\mathbb{X})$ of
  monoids and monoid homomorphisms in $\mathbb{X}$.
\end{lemma}
\begin{proof}
  Orthogonality is supplied by the universal property of cartesian
  cells, together with the fact that isomorphisms are orthogonal to
  anything.  Factorizations are given by restriction along the
  vertical arrow component of a monoid homomorphism, using the
  previous lemma.
\end{proof}

\begin{theorem}\label{thm:norm-disc-eqv}
  Let $\mathbb{X}$ be a virtual equipment, and let $T$ be a monad on
  $\Mod(\mathbb{X})$ which preserves b.o.\ morphisms.  Then
  \begin{enumerate}
  \item $\dMod(\Mod(\mathbb{X}),T)$ is coreflective in
    $\KMod(\Mod(\mathbb{X}),T)$ (that
    is, its inclusion has a right adjoint in \vEquip),\label{item:nd1}
  \item $\nMod(\Mod(\mathbb{X}),T)$ is reflective in
    $\KMod(\Mod(\mathbb{X}),T)$, and\label{item:nd2}
  \item the induced adjunction
    \[\dMod(\Mod(\mathbb{X}),T)\two/>`<-/ \nMod(\Mod(\mathbb{X}),T)\]
    is an adjoint equivalence.\label{item:nd3}
  \end{enumerate}
\end{theorem}
\begin{proof}
  We first prove~\ref{item:nd1}.  Let $A\hto^M TA$ be a $T$-monoid in
  $\Mod(\mathbb{X})$, where $A_0\hto^A A_0$ is a monoid in
  $\mathbb{X}$.  The unit of $A$ is a monoid homomorphism $e\maps
  U_{A_0}\to A$, so by \autoref{lem:monoid-restr}, $M(Te,e)\maps
  U_{A_0}\hto T(U_{A_0})$ is an (object-discrete) $T$-monoid.
  Likewise, if $B\hto^N TB$ is another $T$-monoid and $M\hto^p TN$ is
  a horizontal arrow in $\KMod(\Mod(\mathbb{X}),T)$, then
  $p(Te_B,e_A)$ is a horizontal arrow from $M(Te_A,e_A)$ to
  $N(Te_B,e_B)$.  It is straightforward to extend these
  constructions to a functor $\KMod(\Mod(\mathbb{X}),T)\to
  \dMod(\Mod(\mathbb{X}),T)$.

  Note that $M(Te,e)$ comes with a natural map to $M$, whose
  vertical arrow component is $e$, and likewise for $p(Te_B,e_A)$.
  This supplies the counit of the desired coreflection.  We obtain the
  unit by observing that if $M$ were already object-discrete, then $e$
  would be the identity, so we would have $M(Te,e)\iso M$.  The
  triangle identities are easy to check.
  
  Note that the horizontal arrow in $\mathbb{X}$ underlying
  $M(Te,e)$ is the restriction of $A_0\hto^M (TA)_0$ along the
  identity $e_0\maps A_0 = A_0$ and the map $(Te)_0\maps T(U_{A_0})_0
  \to (TA)_0$.  Since $e$ is b.o.\ and $T$ preserves b.o.\ morphisms,
  $(Te)_0$ is an isomorphism; thus the coreflection of $M$ leaves its
  underlying horizontal arrow in $\mathbb{X}$ essentially unmodified.

  We now prove~\ref{item:nd2}; let $A$ and $M$ be as before.  We first
  observe that the $T$-monoid $M$ in $\Mod(\mathbb{X})$ has an
  underlying monoid in $\Mod(\mathbb{X})$, namely $M(\eta,1)$.  (This is
  a special case of a general functoriality result we will prove
  in~\cite{cs:functoriality}.)  As noted in \autoref{rmk:mod-monad}, a
  monoid in $\Mod(\mathbb{X})$ consists of two monoids in $\mathbb{X}$
  and a monoid homomorphism between them whose vertical arrow
  components are identities.  In this case the first monoid is of
  course $A$.  We denote the second by $A'$ and the monoid homomorphism
  by $c\colon A\to A'$.
  Note that the underlying horizontal arrow of $A'$ in $\mathbb{X}$ is
  just $M(\eta,1)$.

  Now since $T$ preserves b.o.\ morphisms, $TA \to^{Tc} TA'$ is b.o.,
  hence $(TA)_0 \to^{(Tc)_0} (TA')_0$ is an isomorphism.  By
  restricting along its inverse and using the identity $(A')_0 =
  A_0$, from $A_0\hto^M (TA)_0$ we obtain a horizontal arrow $(A')_0
  \hto (TA')_0$.  We abuse notation by continuing to denote this $M$
  (since restriction along isomorphisms leaves an arrow essentially
  unchanged).
  Now $A'$ is a restriction of $M$, so it acts on $M$
  from the left via the multiplication of $M$.  And since $T$
  preserves restrictions, $TA'$ is a restriction of $TM$, so it also
  acts on $M$ from the right via the multiplication of $M$.  Thus, the
  horizontal arrow in $\mathbb{X}$ underlying $M$ also admits the
  structure of a horizontal arrow $A'\hto TA'$ in $\Mod(\mathbb{X})$,
  which we denote $M'$.  Likewise, the multiplication and unit of the
  $T$-monoid $M$ induce a multiplication and unit on $M'$, making it
  also into a $T$-monoid, and we have a canonical $T$-monoid
  homomorphism $M\to M'$ which is an isomorphism in $\mathbb{X}$.
  By definition of $A'$, $M'$ is normalized.
  There is an analogous construction on horizontal arrows, and
  together they extend straightforwardly to a functor
  $\KMod(\Mod(\mathbb{X}),T)\to \nMod(\Mod(\mathbb{X}),T)$.

  Now recall that $A'$ came equipped with a b.o.\ monoid homomorphism
  $A\to A'$.  It is straightforward to check that this homomorphism
  underlies a $T$-monoid homomorphism $M\to M'$; in this way we obtain
  the unit of the desired reflection.  We obtain its counit by
  observing that if $M$ is already normalized, then $A\iso A'$ and
  hence $M\iso M'$.  The triangle identities are again easy to check.

  Finally, to show~\ref{item:nd3}, we observe that the unit of the
  reflection and the counit of the coreflection are isomorphisms on
  the underlying horizontal arrow in $\mathbb{X}$ (since they are
  restrictions along an isomorphism).  Moreover, the reflection and
  coreflection functors both \emph{invert} morphisms with this
  property.  Statement~\ref{item:nd3} then follows formally.
\end{proof}

Recall that we began this section by observing that ordinary multicategories can be
recovered as either object-discrete or normalized $T$-monoids, when
$T$ is the ``free strict monoidal category'' monad on $\PROF{Set}$.
Since this $T$ preserves b.o.\ morphisms, this statement is indeed an
instance of \autoref{thm:norm-disc-eqv}.
However, ordinary multicategories can also be obtained as
\emph{arbitrary} $S$-monoids, when $S$ is the free monoid monad on
$\MAT{Set} = \SPAN{Set}$.  Noting that in this case $T=\Mod(S)$, we
generalize this statement to the following.

\begin{theorem}\label{thm:nd-eqv-mon}
  Let $S$ be a monad on a virtual equipment $\mathbb{X}$.  Then the
  monad $\Mod(S)$ on $\Mod(\mathbb{X})$ preserves b.o.\ morphisms, and we
  have a diagram
  \[\xymatrix{\dMod(\Mod(\mathbb{X}),\Mod(S)) \ar@<1mm>[r]
    \ar@{<-}@<-1mm>[r] \ar[dr] &
    \KMod(\Mod(\mathbb{X}),\Mod(S)) \ar@<1mm>[r]
    \ar@{<-}@<-1mm>[r] \ar[d] &
    \nMod(\Mod(\mathbb{X}),\Mod(S)) \ar[dl]\\
    &\KMod(\mathbb{X},S)
  }\]
  which serially commutes (up to isomorphism).  Moreover, the two
  diagonal functors
  \begin{align}
    \dMod(\Mod(\mathbb{X}),\Mod(S)) &\to \KMod(\mathbb{X},S)\label{eq:dmod-eqv}\\
    \nMod(\Mod(\mathbb{X}),\Mod(S)) &\to \KMod(\mathbb{X},S)\label{eq:nmod-eqv}
  \end{align}
  are equivalences.
\end{theorem}
\begin{proof}
  By definition, $\Mod(S)$ takes a monoid $A_0\hto^A A_0$ to
  $SA_0\hto^{SA} SA_0$, so it preserves b.o.\ morphisms since $S$
  preserves isomorphisms.  We define the middle vertical arrow 
  \begin{equation}
    \KMod(\Mod(\mathbb{X}),\Mod(S)) \to \KMod(\mathbb{X},S)\label{eq:kmodmod-forget}
  \end{equation}
  by applying the 2-functor $\Mod$ to the functor
  \[\HKl(\Mod(\mathbb{X}),\Mod(S)) \to \HKl(\mathbb{X},S)\]
  which takes a monoid $A_0\hto^A A_0$ to its underlying object $A_0$,
  and similarly for horizontal arrows.  (This is again a special case
  of the general functorial result of~\cite{cs:functoriality}.)
  Thus,~\eqref{eq:kmodmod-forget} takes a
  $\Mod(S)$-monoid $A\hto^{M} \Mod(S)(A)$ to the $S$-monoid $A_0
  \hto^M SA_0$.

  We define the diagonal functors by composition with this, so that
  the triangles $\vcenter{\xymatrix@R=1pc@C=1pc{ \ar[r]\ar[dr] &
      \ar[d]\\ & }}$ and $\vcenter{\xymatrix@R=1pc@C=1pc{ \ar[d] &
      \ar[l]\ar[dl]\\ & }}$ commute by definition.  The other
  triangles commute up to isomorphism because the reflection and
  coreflection were defined to fix $A_0$, replace $A$, and restrict
  $M$ along an isomorphism, whereas~\eqref{eq:kmodmod-forget} simply
  forgets about $A$.

  Now, by the 2-out-of-3 property for equivalences, it suffices to
  show that~\eqref{eq:dmod-eqv} is an equivalence.  We will construct
  an explicit inverse to it.  By \autoref{prop:mod-adjt}, to construct
  a normal functor
  \begin{equation}
    \KMod(\mathbb{X},S) \to \KMod(\Mod(\mathbb{X}),\Mod(S))\label{eq:dmod-eqv-inv}
  \end{equation}
  it suffices to construct a not-necessarily-normal functor
  \begin{equation}
    \KMod(\mathbb{X},S) \to \HKl(\Mod(\mathbb{X}),\Mod(S)).\label{eq:dmod-eqv-inv-adjunct}
  \end{equation}
  We define~\eqref{eq:dmod-eqv-inv-adjunct} on objects by sending an
  $S$-monoid $A_0\hto^A T A_0$ to $U_{A_0}$, and likewise on vertical
  arrows.  A horizontal arrow $A\hto^p B$ in $\KMod(\mathbb{X},S)$ has
  an underlying horizontal arrow $A_0 \hto^p T B_0$ in $\mathbb{X}$,
  which acquires a $T(U_{B_0})$-module structure from the following
  composite:
  \[\bfig
  \square/`=`>`{@{>}|-*@{|}}/<1000,375>[A_0`T^2B_0`A_0`TB_0;``\mu`p]
  \square(500,375)/{@{>}|-*@{|}}`=`>`{@{>}|-*@{|}}/<500,375>[TB_0`TB_0`TB_0`T^2B_0;T(U_{B_0})``T\eta`TB]
  \square(0,375)/{@{>}|-*@{|}}`=`=`{@{>}|-*@{|}}/<500,375>[A_0`TB_0`A_0`TB_0;p```p]
  \place(250,550)[=]
  \place(700,600)[\twoar(0,-1)]
  \place(775,550)[\scriptstyle T\hat{b}]
  \place(450,225)[\twoar(0,-1)]
  \place(525,200)[\scriptstyle \bar{p}_r]
  \efig
  \]
  This defines~\eqref{eq:dmod-eqv-inv-adjunct} on horizontal arrows;
  its action on cells is straightforward.  The induced normal
  functor~\eqref{eq:dmod-eqv-inv} takes an $S$-monoid $A_0\hto^A TA_0$
  to itself, regarded as an $(U_{A_0},T(U_{A_0}))$-bimodule.
  Clearly~\eqref{eq:dmod-eqv-inv} followed by the forgetful functor is
  the identity, while the composite in the other direction is
  precisely the coreflection functor into $\dMod(\Mod(\mathbb{X}),\Mod(S))$;
  this completes the proof.
\end{proof}

\begin{examples}
  As remarked above, when $S$ is the ``free monoid'' monad on
  $\SPAN{Set}$, this shows that ordinary multicategories (i.e.\
  $S$-monoids) can also be identified with object-discrete or
  normalized $\Mod(S)$-monoids.  Likewise, virtual double categories
  are $S$-monoids for the ``free category'' monad, and thus can also
  be identified with object-discrete or normalized $\Mod(S)$-monoids.
\end{examples}

\begin{example}\label{eg:modular}
  Since topological spaces can be identified with $U$-monoids in
  $\Rel$, they can also be identified with object-discrete or
  normalized $\Mod(U)$-monoids in $\PROF{2}$.  In the terminology
  of~\cite{tholen:ordered-topstruct}, a $\Mod(U)$-monoid is a
  \emph{modular topological space}.  It is normalized precisely when
  its order is the specialization order, so that it is equivalent to
  an ordinary topological space---i.e.\ a $U$-monoid, as required by
  \autoref{thm:nd-eqv-mon}.
\end{example}

By no means are all interesting monads on $\Mod(\mathbb{X})$ of the
form $\Mod(S)$.  However, many of them do preserve b.o.\ morphisms, so
that \autoref{thm:norm-disc-eqv} at least applies.

\begin{example}
  The ``free symmetric strict monoidal category'' monad on $\PROF{Set}
  = \Mod(\MAT{Set})$ preserves b.o.\ morphisms but is not of the form
  $\Mod(S)$ for any monad $S$ on $\MAT{Set}$.  We have seen in
  \autoref{eg:sym-multicat} that object-discrete $T$-monoids are symmetric
  multicategories; hence so are normalized $T$-monoids.
\end{example}

\begin{example}\label{eg:lawvere-theories-2}
  The ``free category with strictly associative finite products''
  monad on $\PROF{Set}$ also preserves b.o.\ morphisms but is not of
  the form $\Mod(S)$.  We have seen in
  \autoref{eg:lawvere-theories} that object-discrete $T$-monoids are
  multi-sorted Lawvere theories; hence so are normalized $T$-monoids.
\end{example}

\begin{nonexample}
  Recall from \autoref{eg:clubs} that \emph{clubs} are $T$-monoids
  in $\SPAN{Cat}$ with a discrete category of objects, where $T$ is a
  monad like the previous two.  However, since $\SPAN{Cat}$ is not
  of the form $\Mod(\mathbb{X})$, the theory of this section does not
  apply to clubs.  In particular, their ``object-discreteness'' is not
  an instance of our definition, and is not the same as normalization.
\end{nonexample}

When $T$ does not preserve b.o.\ morphisms, however, normalized and
discrete $T$-monoids can be quite different, even on a virtual
equipment of the form $\Mod(\mathbb{X})$.  Intuitively, saying that
$T$ preserves b.o.\ morphisms says that the possible domains of
multimorphisms in a $T$-multicategory depend only on its objects.  If
this fails to be true, then how many morphisms are included in the
underlying monoid can change what these possible domains are.

\begin{example}
  Let $T$ be the ``free category with equalizers'' monad on $\PROF{Set}$.
  Then $T$ evidently does not preserve b.o.\ morphisms, but it is the
  identity on discrete categories.  Therefore, an object-discrete
  $T$-monoid is just a category, whereas a normalized $T$-monoid can
  have morphisms whose domain is a ``formal equalizer'' of ordinary
  morphisms.

  More interestingly, normalized $T$-monoids for the ``free category
  with finite limits'' monad (which also does not preserve b.o.\
  morphisms) can be considered a generalization of Lawvere theories to
  finite-limit logics.  We can also continue to generalize to more
  powerful logics (or ``doctrines'').
\end{example}

These examples suggest that when $T$ does not
preserve b.o.\ morphisms, it is often the normalized, rather than the
object-discrete, $T$-monoids that better capture the desired notion of
$T$-multicategory.  Note also that normalization makes sense for any
monad on a virtual equipment, while object-discreteness only makes
sense for monads on virtual equipments of the form $\Mod(\mathbb{X})$.
Finally, we will see in the next section that normalized $T$-monoids
are the most natural notion to compare with pseudo $T$-algebras.
This inspires us to take normalized $T$-monoids as our preferred
definition of ``generalized multicategory,'' and to make the following
informal definition.

\begin{definition}
  If $T$ is a monad on a virtual equipment for which (possibly pseudo)
  $T$-algebras are called \emph{widgets}, then normalized $T$-monoids
  are called \textbf{virtual widgets}.
\end{definition}

The reasons for this definition were summarized in the introduction.  In
\S\ref{sec:vert-alg} we will prove that any widget has an underlying
virtual widget, further justifying the terminology.  Of course, we have
seen that a number of types of virtual algebras already have their own
names, such as ``multicategory'' and ``Lawvere theory.''  When such
common names exist, we of course use them in preference to terms such as
``virtual monoidal category'' or ``category with virtual finite
products.''

Note that ``virtual widget'' is, strictly speaking, ambiguous: knowing
the notion of widget determines at most the vertical 2-category
$\calV\mathbb{X}$ and the 2-monad $\calV T$, rather than $\mathbb{X}$
and $T$ themselves.  However, many 2-categories that arise in practice
come with an obvious ``natural'' extension to a virtual equipment, so in
practice there is little ambiguity.  (In fact, there is a general
construction of an equipment from a well-behaved 2-category;
see~\cite{cjsv:modulated-bicats}.)  One case of ambiguity is if
``widget'' is the name for $T$-algebras in \Set\ or \Cat, but we
consider $T$-monoids in $\MAT{V}$ or $\PROF{V}$; in this case we may
speak of \textbf{\V-enriched virtual widgets}.

\begin{remark}
  The discussion above suggests that when the objects of $\mathbb{X}$
  are category-like, it is the normalized $T$-monoids (i.e.\ virtual
  $T$-algebras) that are more important, while when the objects of
  $\mathbb{X}$ are set-like, it is the non-normalized $T$-monoids
  (i.e. virtual $\MOD(T)$-algebras) that are more important.
  This does seem to \emph{usually} be the case, but there are exceptions
  on both sides, such as the following.
  \begin{itemize}
  \item As we have already remarked, the multicategories
    of~\cite{bd:hda3} and~\cite{cheng:opetopic} are non-normalized
    $T$-monoids, when $T$ is the ``free symmetric strict monoidal
    category'' monad on $\PROF{Set}$ (whose objects are obviously
    category-like).
  \item Let $U$ be the ultrafilter monad on $\Rel$, whose objects are
    set-like.  We have seen that a $U$-monoid is just a topological
    space, but it is easy to verify that a $U$-monoid is normalized
    just when it is a $T_1$-space---certainly also an important
    concept.
  \end{itemize}
\end{remark}


\section{Representability}
\label{sec:vert-alg}

We now turn to a general version of the comparison between monoidal
categories and multicategories.  Of course, we first need to identify
the analogue of a monoidal category in the general case.  We saw in
\S\ref{sec:normalization} that ordinary multicategories have two
different faces in our setup: they are the $S$-monoids where $S$ is
the ``free monoid'' monad on $\SPAN{Set}$, and also the normalized
$T$-monoids, where $T=\Mod(S)$ is the ``free strict monoidal
category'' monad on $\PROF{Set}$.  Monoidal categories, however, are
more visible from the second point of view: they are the \emph{pseudo
  $\calV(T)$-algebras} in $\calV(\PROF{Set})=\Cat$.

Accordingly, in this section we will assume that $T$ is a monad on a
virtual equipment $\mathbb{X}$ whose objects are ``category-like,'' and seek to
compare (pseudo) $\calV(T)$-algebras with (normalized) $T$-monoids.
We will additionally have to assume that $T$ is a \emph{normal} monad
as defined in \S\ref{sec:2cat-tmon}, since otherwise it doesn't even
induce a 2-monad on $\calV(\mathbb{X})$.  If we are given instead a
monad $S$ on a virtual double category whose objects are ``set-like,'' then in order to
apply the theory of this section we simply consider $\Mod(S)$ instead;
some examples of this can be found later on.  Generalizing the
terminology of \cite[p.~165]{leinster:higher-opds}, we may call a
(pseudo) $\calV(\Mod(S))$-algebra an \textbf{$S$-structured monoid}.

Actually, the most natural approach to the comparison turns out to be
via \emph{oplax} $T$-algebras.  Recall that for a 2-monad $T$ on a
2-category, an \textbf{oplax $T$-algebra} is an object $A$ with a map
$a\maps TA\to A$ and 2-cells
\begin{equation}
  \vcenter{\xymatrix{T^2A \ar[r]^\mu\ar[d]_{Ta} \drtwocell\omit{\overline{a}} & TA \ar[d]^a\\
      TA \ar[r]_a & A}}\qquad\text{and}\qquad
  \vcenter{\xymatrix{A \ar[r]^\eta \ar@{=}[dr]^{}="mid" &
      TA \ar[d]^a \ar@{=>}[0,0];"mid"^{\widehat{a}} \\ & A }}\label{eq:oplaxalg-2cells}
\end{equation}
satisfying certain straightforward axioms.  We call it \emph{normal}
if $\widehat{a}$ is an isomorphism, and a \emph{pseudo $T$-algebra} if
both $\widehat{a}$ and $\overline{a}$ are isomorphisms.
Finally, if $T$ is a monad on a virtual equipment, we will always abuse
terminology by saying ``$T$-algebra'' (with appropriate prefixes) to
mean $\calV(T)$-algebra.

\begin{theorem}\label{thm:oplax-alg}
  Let $T$ be a normal monad on a virtual equipment $\mathbb{X}$.  Then:
  \begin{enumerate}
  \item Any oplax $T$-algebra $TA\to^a A$ in $\calV\mathbb{X}$
    gives rise to a $T$-monoid $A\hto^{A(a,1)} TA$, which is
    normalized if and only if $A$ is normal.\label{item:oplax1}
  \item A $T$-monoid $A\hto^M TA$ arises from an oplax
    $T$-algebra if and only if $M\iso A(a,1)$ for some vertical arrow
    $TA \to^a A$.\label{item:oplax2}
  \end{enumerate}
\end{theorem}
\begin{proof}
  If $a\maps TA\to A$ is an oplax $T$-algebra, then by definition of
  $A(a,1)$, the 2-cell $\widehat{a}$ induces a unit cell
  \[\xymatrix{A \ar[r]|-@{|}^-{U_A} \ar@{=}[d] \ar@{}[dr]|{\twoar(0,-1)} & A \ar[d]^\eta\\
    A\ar[r]|-@{|}_-{A(a,1)} & TA.}\]
  Likewise, by the dual of \autoref{thm:companion-mates},
  $\overline{a}$ induces a multiplication
  \[\xymatrix@C=3pc{
    A \ar[r]|-@{|}^-{A(a,1)} \ar@{=}[d] \ar@{}[drr]|{\twoar(0,-1)} & TA
    \ar[r]|-@{|}^-{(TA)({Ta},1)} &
    T^2 A\ar[d]^\mu\\
    A \ar[rr]|-@{|}_-{A(a,1)} && TA}
  \]
  and using the isomorphism $(TA)({Ta},1)\iso T(A(a,1))$ (since $T$ is
  normal) we obtain a
  multiplication cell.  The axioms to make $A(a,1)$ into a $T$-monoid
  follow directly from the axioms for an oplax $T$-algebra.  To
  complete~\ref{item:oplax1}, we observe that $\widehat{a}$ is an
  isomorphism if and only if the induced cell $U_A \to
  A(a\eta,1) \iso A(a,1)(\eta,1)$ is an isomorphism, which says precisely
  that the unit defined above is cartesian.  Conversely, if
  $A\hto^M TA$ is a $T$-monoid and $M\iso A(a,1)$, then the same bijections
  supply 2-cells $\overline{a}$ and $\widehat{a}$ satisfying the same
  axioms making $a\maps TA\to A$ into an oplax $T$-algebra; this
  shows~\ref{item:oplax2}.
\end{proof}

The following example may serve to clarify the connection between
normality of oplax $T$-algebras and normalization of $T$-monoids.

\begin{example}\label{eg:oplax-moncat}
  Let $T$ be the ``free strict monoidal category'' monad on
  $\PROF{Set}$.  Then an oplax $T$-algebra is an \emph{oplax monoidal
    category}: a category $A$ equipped with tensor product functors
  \begin{align*}
    A\times \dots\times A &\to A\\
    (x_1,\dots, x_n) &\mapsto \langle x_1\ten\dots\ten x_n\rangle
  \end{align*}
  for $n\ge 0$, and transformations
  \begin{align*}
    \langle x\rangle &\to x\\
    \langle x_{11}\ten \dots \ten x_{nk_n}\rangle &\to
    \langle\langle x_{11}\ten\dots\ten x_{1k_1}\rangle
    \ten\dots\ten
    \langle x_{n1}\ten\dots\ten x_{nk_n}\rangle\rangle
  \end{align*}
  satisfying certain evident axioms.  Note that the 0-ary tensor
  product $\langle\rangle$ is a ``lax unit'' and the 1-ary tensor
  product $\langle x\rangle$ is not necessarily isomorphic to $x$,
  only related by the given unit transformation $\langle x\rangle\to
  x$.

  As mentioned previously, a $T$-monoid consists of a category $A$, a
  multicategory $M$ with the same objects, and an
  identity-on-objects functor from $A$ to the underlying ordinary
  category of $M$.  Now \autoref{thm:oplax-alg} says
  that we can make an oplax monoidal category into a $T$-algebra by
  defining the multimorphisms in $M$ from $(x_1,\dots,
  x_n)$ to $y$ to be the morphisms $\langle x_1\ten\dots\ten
  x_n\rangle\to y$ in $A$.

  Note that the morphisms from $x$ to $y$ in the underlying ordinary
  category of $M$ are the morphisms from $\langle
  x\rangle$ to $y$ in $A$.  The functor $A\to M$ is
  defined by composing with the unit transformation $\langle
  x\rangle\to x$.  Clearly this is fully faithful (i.e.\ the
  $T$-monoid is normalized) just when $\langle x\rangle\to x$ is an
  isomorphism (i.e.\ the oplax $T$-algebra is normal).
\end{example}

The following characterization of pseudo $T$-algebras is now obvious.

\begin{corollary}\label{thm:psalg}
  A normalized $T$-monoid $A\hto^M TA$ arises from a pseudo
  $T$-algebra if and only if
  \begin{enumerate}
  \item $M\iso A(a,1)$ for some $TA\to^a A$, and\label{item:psalg1}
  \item the induced 2-cell $\overline{a}$ is an isomorphism.\label{item:psalg2}
  \end{enumerate}
\end{corollary}

We say that a normalized $T$-monoid is \textbf{weakly representable}
if it satisfies~\ref{item:psalg1}, and \textbf{representable} if it
satisfies both~\ref{item:psalg1} and~\ref{item:psalg2} (hence is
equivalent to a pseudo $T$-algebra).

\begin{example}\label{eg:repr-moncat}
  When $T$ is the ``free strict monoidal category'' monad on
  $\PROF{Set}$, \autoref{thm:psalg} specializes to the
  characterization of monoidal categories as representable
  multicategories, as in~\cite{HermidaRepresentable}
  and~\cite[\S3.3]{leinster:higher-opds}. We will see in
  \S\ref{sec:cartesian-2-monads} that it also includes the general
  representability notion of~\cite{HermidaCoherent}. The analogue of
  \autoref{thm:oplax-alg} in the language of~\cite{burroni:t-cats} can
  be found in~\cite{penon:tcats-rep}, which uses ``representable'' for
  what we call ``weakly representable'' and ``lax algebra'' for what
  we call an ``oplax algebra.''
\end{example}

\begin{remark}\label{rmk:bias}
  Strictly speaking, the notion of monoidal category obtained in this
  way is the ``unbiased'' version, which is equipped with a specified
  $n$-ary tensor product for all $n\ge 0$, instead of the usual
  ``biased'' version having only a binary and nullary product
  (see~\cite[\S3.1]{leinster:higher-opds}).  This
  is generally what happens for pseudoalgebras: if $T$ is a monad
  whose strict algebras are some strict structure, then pseudo
  $T$-algebras are an ``unbiased'' sort of weak structure.  Generally
  the unbiased version is equivalent to the biased one, but there is
  real mathematical content in this statement; for instance, the
  equivalence of biased and unbiased monoidal categories is
  essentially equivalent to Mac Lane's coherence theorem.
\end{remark}

\begin{example}\label{eg:vdblcat}
  Recall that virtual double categories can be identified with
  $S$-monoids for the ``free category'' monad $S$ on directed graphs,
  and hence also with normalized $\Mod(S)$-monoids.  In this case it
  is easy to check that for a normalized $\Mod(S)$-monoid
  $A\hto^M TA$, we have $M\iso A(a,1)$ iff every composable string of
  horizontal arrows is the source of a \emph{weakly} opcartesian arrow
  (see \autoref{thm:weak-opcart}).  Thus, such virtual double
  categories can be identified with ``normal oplax double categories,''
  which are equipped with $n$-ary composites for all $n$ and
  comparison maps
  \[\langle p_{11}\odot \dots\odot p_{nk_n}\rangle \to
  \langle \langle p_{11} \odot \dots \odot p_{1k_1} \rangle \odot
  \dots\odot \langle p_{n1} \odot \dots\odot p_{nk_n}\rangle\rangle
  \]
  and invertible comparison maps
  \[\langle p \rangle \to^\iso p\]
  satisfying analogous axioms to an oplax monoidal category.
  Condition~\ref{item:psalg2} in
  \autoref{thm:psalg} is then equivalent to requiring weakly
  opcartesian cells to be closed under composition.  As observed in
  \autoref{thm:weak-opcart}, this suffices to ensure we have a
  pseudo double category, i.e.\ a pseudo $\Mod(S)$-algebra.
\end{example}

\begin{example}
  Let $U$ be the ultrafilter monad on $\Rel$.  We have seen that a
  $U$-monoid is a topological space, and a
  normalized $U$-monoid is a $T_1$-space.  A vertical
  $U$-algebra (which is automatically strict, since $\calV(\Rel)$ is
  locally discrete) is a compact Hausdorff space, and in
  this case \autoref{thm:oplax-alg} tells us what we already knew:
  any compact Hausdorff space is, in particular, a $T_1$ topological
  space.

  Now consider the induced monad $\Mod(U)$ on $\Mod(\Rel)$.  The
  objects of $\Mod(\Rel)=\PROF{2}$ are preorders.  In the language
  of~\cite{tholen:ordered-topstruct}, a strict $\Mod(U)$-algebra is an
  \emph{ordered compact Hausdorff space}, whereas by
  \autoref{thm:nd-eqv-mon} a normalized $\Mod(U)$-monoid is simply a
  topological space.  Thus, \autoref{thm:oplax-alg} tells us that any
  ordered compact Hausdorff space $X$ can be equipped with a topology
  in which an ultrafilter $\mathfrak{F}$ converges to a point $y$ if
  and only if the (unique) limit of $\mathfrak{F}$ in $X$ is $\le y$
  in the given preorder.
\end{example}

The next three examples can all be found in~\cite{HermidaCoherent} (see
\ref{sec:cartesian-2-monads} for more on the comparison between our
setting and Hermida's).

\begin{example}\label{eg:fibrations}
  Let $\mb{S}$ be a small category and $T$ the monad on
  $\C=\Set^{\ob(\mb{S})}$ whose algebras are functors $\mb{S}\to
  \Set$, as in \autoref{eg:cat-over}, and consider the monad
  $\Mod(T)$ on $\Mod(\SPAN{C}) = \iPROF{C}$.  A strict
  $\Mod(T)$-algebra is a functor $\mb{S}\to\Cat$, while a pseudo
  $\Mod(T)$-algebra is a pseudofunctor $\mb{S}\to\Cat$.  Now by
  \autoref{thm:nd-eqv-mon}, normalized $\Mod(T)$-monoids can be
  identified with $T$-monoids, which as we saw can be identified with
  functors $\mb{A}\to\mb{S}$.  It is then easy to verify that a
  normalized $\Mod(T)$-monoid satisfies
  \ref{thm:psalg}\ref{item:psalg1} iff the corresponding functor
  $\mb{A}\to\mb{S}$ admits all weakly opcartesian liftings, and
  \ref{thm:psalg}\ref{item:psalg2} iff weakly opcartesian arrows are closed
  under composition.  Thus, in this case \autoref{thm:psalg}
  specializes to the classical equivalence between pseudofunctors
  $\mb{S}\to\Cat$ and opfibrations over $\mb{S}$.
\end{example}

\begin{example}\label{thm:discfib}
  Let $T$ be as in \autoref{eg:fibrations}, but now consider
  $T$-monoids rather than $\Mod(T)$-monoids.  A $T$-monoid is
  normalized just when for each $x\in \mb{S}$, the induced span $A(x,1)
  \toleft M \to A(x,1)$ is the identity span; i.e.\ when the fibers of
  $\mb{A}\to\mb{S}$ are discrete categories.  Such a normalized
  $T$-monoid satisfies \ref{thm:psalg}\ref{item:psalg1} iff
  $\mb{A}\to\mb{S}$ admits all weakly opcartesian liftings, which in this
  case are automatically opcartesian by discreteness.  Thus,
  \autoref{thm:psalg} also specializes to the equivalence of functors
  $\mb{S}\to\Set$ and \emph{discrete} opfibrations over $\mb{S}$.
\end{example}

\begin{example}
  Let $T$ be the ``free strict $\omega$-category'' monad on
  $\SPAN{Globset}$, for which we saw in \autoref{eg:glob-opd} that
  $T$-monoids on $1$ are globular multicategories.  Pseudo $T$-algebras are
  an ``unbiased'' version of the \emph{monoidal
    globular categories} of~\cite{batanin:monglob}.  Thus any monoidal
  globular category has an underlying globular multicategory, and can
  be characterized among the latter by a representability property.
\end{example}

The requirement that $T$ be normal in \autoref{thm:oplax-alg} 
cannot be dispensed with.

\begin{example}\label{eg:powerset}
  Let $P$ be the extension of the powerset monad to $\Rel$
  described in \autoref{example:tv}.  Since $\calV(\Rel)=\Set$
  is locally discrete, oplax $P$-algebras are just $P$-algebras in
  \Set, which can be identified with complete meet-semilattices (the
  structure map $PA\to A$ takes a subset $X\subseteq A$ to its meet
  $\bigwedge X$).

  Now, we have observed in \autoref{eg:closure-spaces} that
  $P$-monoids can be identified with closure spaces.  If we attempt to
  follow the prescription of \autoref{thm:oplax-alg}, starting from a
  complete join-semilattice we would define the ``closure operation''
  $c(X) = \{\bigwedge X\}$; but this is neither extensive nor
  monotone.

  On the other hand, if we first apply $\Mod$, we obtain a monad
  $\Mod(P)$ on $\Mod(\Rel)=\PROF{2}$, which is normal.
  By \autoref{thm:nd-eqv-mon}, normalized
  $\Mod(P)$-monoids can be identified with $P$-monoids, i.e.\ closure
  spaces.  With a little effort, pseudo $\Mod(P)$-algebras can be
  identified with meet-complete preorders (that is, preorders that are
  complete as categories).  \autoref{thm:oplax-alg} then tells us
  that from a meet-complete preorder we can construct a closure space
  with $c(X) = \{x \mid \bigwedge X \le x\}$, which is certainly true.
\end{example}

We can also make the correspondence of \autoref{thm:oplax-alg}
functorial.  Recall that for any
$T$ we have a 2-category $\KMon(\mathbb{X},T)$ of $T$-monoids, defined
to be the vertical 2-category of $\KMod(\mathbb{X},T)$.  It turns out
that while $T$-monoids correspond to \emph{oplax} $T$-algebras,
morphisms of $T$-monoids correspond to \emph{lax} $T$-algebra
morphisms.  Recall that a \emph{lax $T$-morphism} between oplax
$T$-algebras consists of a map $f\maps A\to B$ and a 2-cell
\[\vcenter{\xymatrix{TA \ar[r]^{Tf}\ar[d]_{a} \drtwocell\omit{\overline{f}} & TB \ar[d]^b\\
    A\ar[r]_f & B}}\]
satisfying certain straightforward axioms.  And if $A\two^f_g B$
are two such, a \emph{$T$-transformation} is a 2-cell $f\to^\alpha g$ such that
\[\vcenter{\xymatrix@R=0.6pc{TA \ar[r]^{Tf}\ar[ddd]_a \ddrtwocell\omit{\overline{f}} & TB \ar[ddd]^b\\ &\\ & \\
    A\rtwocell^{f}_g{\alpha} & B}}\quad=\quad
\vcenter{\xymatrix@C=3pc@R=0.6pc{ TA \ar[ddd]_a \rtwocell^{Tf}_{Tg}{\;\;T\alpha}
    & TB \ar[ddd]^b\\
    \ddrtwocell\omit{\overline{f}}\\ & \\
    A\ar[r]_f & B.}}
\]
We write
$\mathcal{O}\mathit{plax}\text{-}T\text{-}\mathcal{A}\mathit{lg}_l$ for
the resulting 2-category.

\begin{theorem}\label{thm:oplaxalg-funct}
  Let $T$ be a normal monad on a virtual equipment $\mathbb{X}$.  Then
  there is a strict 2-functor
  $\mathcal{O}\mathit{plax}\text{-}T\text{-}\mathcal{A}\mathit{lg}_l \to
  \KMon(\mathbb{X},T)$, whose underlying 1-functor is fully faithful,
  and which becomes 2-fully-faithful (that is, an isomorphism on
  hom-categories) when restricted to \emph{normal} oplax $T$-algebras.
\end{theorem}
\begin{proof}
  For oplax $T$-algebras $(A,a)$ and $(B,b)$, regarded as $T$-monoids
  $A\hto^{A(a,1)} TA$ and $B\hto^{B(b,1)} TB$, a morphism of
  $T$-monoids consists of a vertical arrow $f\maps A\to B$ and a
  2-cell
  \[\vcenter{\xymatrix{A \ar[r]|-@{|}^-{A(a,1)} \ar[d]_f \ar@{}[dr]|{\twoar(0,-1)} &TA  \ar[d]^{Tf}\\
      B\ar[r]|-@{|}_-{B(b,1)} & TB}}\]  
  satisfying certain axioms.  But by definition of $B(b,1)$, and
  by \autoref{thm:comp-bc-gives-extn} applied to $A(a,1)$, this 2-cell is
  equivalent to one
  \[\vcenter{\xymatrix{
      TA\ar[r]|-@{|}^{U_{TA}}\ar[d]_a \ar@{}[ddr]|{\twoar(0,-1)} &
      TA\ar[d]^{Tf}\\
      A\ar[d]_f &
      TB\ar[d]^b\\
      \ar[r]|-@{|}_{U_B} &.
    }}\]
  This defines a 2-cell $b\circ Tf \to f\circ a$ in $\calV(\mathbb{X})$, which
  is precisely the additional data required to make $f$ into a lax
  $T$-algebra morphism.  It is easy to verify that the axioms of a
  $T$-monoid morphism are equivalent to the axioms of a lax
  $T$-algebra morphism under this translation, and that composition is
  preserved.

  Now let $f,g\maps A\to B$ be two such morphisms, and recall from
  \S\ref{sec:2cat-tmon} that a $T$-monoid transformation $\beta\maps
  f\to g$ consists of a cell
  \[\bfig\scalefactor{.7} \Atriangle[A`B`TB;g`\eta_{TB}\circ f`B(b,1)]
  \morphism(0,0)|m|<1000,0>[B`TB;|]
  \place(500,200)[\Downarrow\scriptstyle \beta]
  \efig
  \]
  satisfying a certain axiom.  Equivalently, $\beta$ is a 2-cell $b
  \eta f \to g$ in $\calV(\mathbb{X})$.  Thus, given a $T$-algebra
  transformation $\alpha\maps f\to g$, it is natural to define $\beta$
  to be the composite 
  \[ b \eta f \to^{\widehat{b} f} f \to^{\alpha} g,
  \]
  where $\widehat{b}\maps b\eta \to 1$ is the oplax unit map of
  $B$.  With this definition, the axiom that must be satisfied for
  $\beta$ to be a $T$-monoid transformation becomes the following
  equality of pasting diagrams in $\calV\mathbb{X}$.
  \begin{equation}
    \bfig\scalefactor{.8}
    \square(500,1000)|amrb|/<-`->`->`<-/[A`TA`B`TB;a`f`Tf`b]
    \place(750,1250)[_{\overline{f}}\twoar(-1,1)]
    \square(500,500)|bmrb|/<-`->`->`<-/[B`TB`TB`T^2B;b`\eta`\eta`Tb]
    \place(750,750)[=]
    \morphism(0,500)/=/<500,500>[B`B;]
    \morphism(500,1500)/{@{>}@/_10pt/}/<-500,-1000>[A`B;g]
    \place(300,1100)[\twoar(-2,1)]
    \place(300,1050)[\scriptstyle \alpha]
    \place(350,650)[_{\widehat{b}}\twoar(-1,1)]
    \morphism(0,500)|b|/<-/<500,0>[B`TB;b]
    \place(500,250)[_{\overline{b}}\twoar(-1,1)]
    \square/`=`->`<-/<1000,500>[B`T^2B`B`TB;``\mu_B`b]
    \efig \;=\;
    \bfig\scalefactor{.8}
    \square(0,500)|almb|/<-`->`->`<-/<500,1000>[A`TA`B`TB;a`g`Tg`b]
    \place(250,1000)[_{\overline{g}}\twoar(-1,1)]
    \morphism(500,1500)/{@{>}@/^10pt/}/<500,-500>[TA`TB;Tf]
    \morphism(1000,1000)|r|<0,-500>[TB`T^2B;T\eta_{B}]
    \morphism(500,500)/=/<500,500>[TB`TB;]
    \place(700,1100)[\twoar(-1,0)]
    \place(700,1025)[\scriptstyle T\alpha]
    \place(850,650)[_{T\widehat{b}}\twoar(-1,1)]
    \morphism(500,500)|b|/<-/<500,0>[TB`T^2B;Tb]
    \place(500,250)[_{\overline{b}}\twoar(-1,1)]
    \square/`=`->`<-/<1000,500>[B`T^2B`B`TB.;``\mu_B`b]
    \efig\label{eq:alg-transfeq}
  \end{equation}
  The cell marked ``$=$'' is an identity by
  \autoref{thm:vert-oplax-is-id} applied to a cell component of $\eta$.
  Now, two of the axioms for an oplax $T$-algebra say that
  \[\bfig\scalefactor{.8}
  \square(500,500)|amrb|/<-`->`->`<-/[B`TB`TB`T^2B;b`\eta`\eta`Tb]
  \place(750,750)[=]
  \morphism(0,500)/=/<500,500>[B`B;]
  \place(350,650)[_{\widehat{b}}\twoar(-1,1)]
  \morphism(0,500)|b|/<-/<500,0>[B`TB;b]
  \place(500,250)[_{\overline{b}}\twoar(-1,1)]
  \square/`=`->`<-/<1000,500>[B`T^2B`B`TB;``\mu_B`b]
  \efig \;=\;
  \vcenter{\xymatrix{B && TB \lltwocell^b_b{^\;1_B} }}
  \;=\;
  \bfig\scalefactor{.8}
  \morphism(1000,1000)|r|<0,-500>[TB`T^2B;T\eta_{B}]
  \morphism(500,500)/=/<500,500>[TB`TB;]
  \place(850,650)[_{T\widehat{b}}\twoar(-1,1)]
  \morphism(500,500)|b|/<-/<500,0>[TB`T^2B;Tb]
  \morphism(0,500)|a|/<-/<500,0>[B`TB;b]
  \place(500,250)[_{\overline{b}}\twoar(-1,1)]
  \square/`=`->`<-/<1000,500>[B`T^2B`B`TB.;``\mu_B`b]
  \efig
  \]
  After removing these composites from~\eqref{eq:alg-transfeq}, what
  is left is simply the equation for $\alpha$ to be a $T$-algebra
  transformation; thus $\alpha$ is such precisely when $\beta = \alpha
  . \widehat{b} f$ is a $T$-monoid transformation.  And, of course, if
  $B$ is normal, then $\widehat{b} f$ is an isomorphism, and so
  $\alpha$ can be recovered uniquely from $\beta$.  We leave it to the
  reader to verify that this association preserves both types of
  2-cell composition.
\end{proof}

The restriction to \emph{normal} oplax algebras in the final statement
of \autoref{thm:oplaxalg-funct} cannot be dispensed with either.

\begin{example}
  Let $A$ and $B$ be oplax monoidal categories, regarded as
  $T$-monoids for the ``free strict monoidal category'' monad $T$ as
  in \autoref{eg:oplax-moncat}, and let $A\two^f_g B$ be lax
  monoidal functors.  A $T$-algebra transformation $f\to g$ is, in
  particular, a natural transformation $f\to g$, and therefore has
  components $fx \to^{\alpha_x} gx$.  However, if we unravel the
  definition of a $T$-monoid transformation, we see that its
  components are of the form $\langle fx\rangle \to^{\beta_x} gx$.
  Thus, when $B$ is not normal, there can be no bijection in general.
\end{example}

\begin{remark}\label{rmk:free-models}
  Recall from \S\ref{sec:introduction} that generalized
  multicategories can be regarded as ``algebraic theories.''  For
  instance, ordinary multicategories correspond to strongly regular
  finitary theories, while Lawvere theories correspond to arbitrary
  finitary theories.  In language introduced by Jon Beck, one may say
  that the monad $T$ provides the ``doctrine'' in which the theories
  are written.  (Motivated by this, some authors use the word
  \emph{doctrine} to mean simply a 2-monad.)
  If $X$ is a $T$-monoid, regarded as a theory in the
  doctrine $T$, and $A$ is a pseudo $T$-algebra, then it is natural to
  define a \emph{model} of $X$ in $A$ to be a $T$-monoid homomorphism
  from $X$ to (the underlying $T$-monoid of) $A$.

Now, frequently the functor of \autoref{thm:oplaxalg-funct} has a left
adjoint $F_T$ when restricted to pseudo $T$-algebras and morphisms.
In such a case, a model of $X$ in $A$ can equally be defined as a
$T$-algebra morphism $F_TX\to A$.  That is, $F_T X$ is the ``free
$T$-algebra containing a model of $X$.''
Following~\cite{lawvere:functsem}, this is sometimes called the
``functorial semantics'' of $T$. 

  For example, when $X$ is a Lawvere theory, $F_T X$ is the category with finite
  products that incarnates it.  (In fact, as we remarked in
  \autoref{eg:lawvere-theories}, Lawvere theories are often
  \emph{defined} to be certain categories with finite products.)
  Likewise, when $X$ is a Lawvere \V-theory, then $F_T X$ is the
  \V-category with finite cotensors that incarnates it
  (see~\cite{power:enr-theories}) and when $X$ is an ordinary or
  non-symmetric operad, $F_T X$ is the ``PROP'' or ``PRO'' associated to
  it (see~\cite{bv:htpy-invar}).

  This adjunction can also be used to characterize representability.
  It turns out that the strict (resp.\ pseudo) algebras for the
  induced 2-monad $T_\dashv$ on $\KMon(\mathbb{X},T)$ can be identified
  with strict (resp.\ pseudo) $T$-algebras.  Moreover, $T_\dashv$ is a
  ``lax-idempotent'' 2-monad in the sense of~\cite{kl:property-like},
  so that $A$ is a pseudo $T_\dashv$-algebra precisely when the unit
  $A\to T_\dashv A$ has a left adjoint.  Thus the ``structure''
  imposed by the 2-monad $T$ has been transformed into ``property-like
  structure'' imposed by $T_\dashv$.  In particular cases, these
  observations can be found
  in~\cite{HermidaRepresentable,HermidaCoherent,penon:tcats-rep};
  in~\cite{cs:laxidem} we will study them in our general context.
\end{remark}

\appendix

\section{Composites in $\Mod$ and $\HKl$}
\label{sec:composites-mod-hkl}

In this appendix we consider the question of when $\Mod(\mathbb{X})$
and $\Hkl{X}{T}$ have composites and units, which will be needed for
our comparisons with existing theories in the next appendix.  The
first case is easy; the following was also observed in
\cite{FramedBicats}.

\begin{theorem}\label{thm:comp-mod}
  If $\mathbb{X}$ is an equipment in which each category
  $\calH\mathbb{X}(A,B)$ has coequalizers, which are preserved on both
  sides by $\odot$, then $\Mod(\mathbb{X})$ is also an equipment.
\end{theorem}
\begin{proof}[(Sketch)]
  We have seen already that $\Mod(\mathbb{X})$ always has units, and
  inherits restrictions from $\mathbb{X}$, so it remains only to
  construct composites.  We define the composite of bimodules $A\hto^p
  B \hto^q C$ to be the coequalizer of the two actions of $B$:
  \[ p\odot B\odot q \two p\odot q \epi p\odot_B q\]
  and similarly for longer composites.  Given a cell
  \[\vcenter{\xymatrix{
      \ar[r]|-@{|}^{p}\ar[d]_f \ar@{}[drr]|{{\twoar(0,-1)}} &
      \ar[r]|-@{|}^{q} &
      \ar[d]^g\\
      \ar[rr]|-@{|}_{r} &&
    }}
  \]
  in $\Mod(\mathbb{X})$, to factor it through $p\odot_B q$, we first
  factor it through a cartesian cell to obtain a cell
  $(p,q)\to r(g,f)$ in $\calH\mathbb{X}(A,C)$, then factor this through
  the coequalizer in
  $\calH\mathbb{X}(A,C)$, and finally compose again with the cartesian
  cell.  Thus $(p,q)\to p\odot_B q$ is weakly opcartesian; to show
  that these cells compose we use the fact that $\odot$ preserves
  coequalizers.
\end{proof}

Note that we require $\mathbb{X}$ to have restrictions, as well as
composites, in order to show that $\Mod(\mathbb{X})$ has composites.
We could instead assume explicitly that the coequalizers in
$\calH\mathbb{X}(A,B)$ satisfy a universal property relative to all
cells, but in practice this generally tends to hold only because of
the existence of restrictions.

\begin{example}
  If \V\ has small colimits preserved by $\otimes$, then $\MAT{V}$
  satisfies the hypotheses of \autoref{thm:comp-mod}, so (as we have
  seen) $\PROF{V}$ is an equipment.
\end{example}

\begin{example}
  If \C\ has pullbacks and coequalizers preserved by pullback, then
  $\SPAN{C}$ satisfies the hypotheses of \autoref{thm:comp-mod}, so
  (as we have also seen) $\iPROF{C}$ is an equipment.
\end{example}

We have also already seen that $\Hkl{X}{T}$ always inherits
restrictions from $\mathbb{X}$.  However, to show that it has
composites, we require fairly strong conditions not just on
$\mathbb{X}$ but on $T$ as well.  Recall that a functor between pseudo
double categories (and, therefore, also between equipments) is called
\emph{strong} if it preserves all composites.  We then make the
following definition:

\begin{definition}
  A transformation $F\to^\alpha G$ of functors $\mathbb{X}
  \two^{F}_G \mathbb{Y}$ between equipments is \textbf{horizontally
    strong} if for every horizontal arrow $X\to^p Y$ in $\mathbb{X}$,
  the cell induced by $\alpha_p$ under the bijection of
  \autoref{thm:companion-mates}:
  \[\vcenter{\xymatrix@C=3pc{
      FX \ar[r]|-@{|}^{Fp}\ar@{=}[d] \ar@{}[drr]|{{\twoar(0,-1)}} &
      FY \ar[r]|-@{|}^{GY(1,{\alpha_Y})}&
      GY \ar@{=}[d]\\
      FX \ar[r]|-@{|}_{GX(1,{\alpha_X})} &
      GX \ar[r]|-@{|}_{Gp} &
      GY
    }}
  \]
  is an isomorphism.  A monad $T$ on an equipment $\mathbb{X}$ is
  \textbf{horizontally strong} if $T$ is a strong functor and $\mu$ and
  $\eta$ are horizontally strong.
\end{definition}

\begin{remark}\label{rmk:hstrong-htrans}
  Recall from \autoref{rmk:oplax-htrans} that any transformation
  $\alpha$ of functors between equipments induces an oplax
  transformation $\calH(\alpha)$ of functors between horizontal
  bicategories.  Horizontal strength of $\alpha$ is equivalent to
  requiring $\calH(\alpha)$ to be a strong (aka pseudo) transformation
  (hence the name).
\end{remark}

\begin{example}\label{eg:cart-mnd}
  Let $\xymatrix{\mb{C} \rtwocell^F_G{\alpha} & \mb{D}}$ be a
  transformation between pullback-preserving functors between
  categories with pullbacks.  It is not hard to verify that the
  induced transformation $\SPAN{\alpha}$ is horizontally strong if and
  only if $\alpha$ is a \textbf{cartesian natural transformation},
  meaning that all its naturality squares are pullbacks.  Therefore,
  if $T$ is a pullback-preserving monad on \C, the monad $\SPAN{T}$ is
  horizontally strong if and only if $\mu$ and $\eta$ are cartesian
  natural transformations.  Such a $T$ is often called a
  \textbf{cartesian monad}; see \S\ref{section:cartesian}.
\end{example}

\begin{example}
  We have remarked that most of the monads from \autoref{example:tv} are
  not even strong functors in general, with the exception of the
  ultrafilter monad on $\Rel$.
  One can verify that for this monad, the multiplication transformation
  is horizontally strong, but the unit transformation is not; hence it
  is not a horizontally strong monad.
\end{example}

\begin{theorem}\label{thm:hklequip}
  If $T$ is a horizontally strong monad on an equipment $\mathbb{X}$,
  then $\Hkl{X}{T}$ is also an equipment.
\end{theorem}
\begin{proof}[(Sketch)]
  It suffices to show that $\Hkl{X}{T}$ has composites.  A composable
  string of horizontal arrows
  \[ X_0 \hto^{p_1} X_1 \hto^{p_2} \cdots \hto^{p_n} X_n\]
  in $\Hkl{X}{T}$ consists of horizontal arrows $X_i \hto^{p_{i+1}}
  TX_{i+1}$ in $\mathbb{X}$.  Since $\mathbb{X}$ is an equipment, we
  can form the composite
  \[ p_1\odot T(p_2)\odot \dots \odot T^{n-1}(p_n) \odot TX_n(1,{\mu^{n-1}})
  \]
  in $\mathbb{X}$, which clearly supplies a \emph{weakly} opcartesian
  cell in $\Hkl{X}{T}$.  Likewise, $TX(1,\eta)$ is a weak unit for
  $X$ in $\Hkl{X}{T}$.  The assumptions on $T$ are required to show
  that these weakly opcartesian cells compose, or equivalently that
  this composition is associative; rather than write this out in
  detail we merely compute the 3-fold associativity isomorphism for $X
  \hto^{p} Y \hto^{q} Z \hto^r W$.
  \begin{alignat*}{2}
    (p \odot Tq \odot TZ(1,\mu)) &\odot Tr \odot TW(1,\mu)\\
    &\iso p \odot Tq \odot T^2r \odot T^2W(1,{\mu T}) \odot TW(1,\mu)
    &&\quad\text{(strength of $\mu$)}\\
    &\iso p \odot Tq \odot T^2r \odot T^2W(1,{T\mu}) \odot TW(1,\mu)
    &&\quad\text{(associativity of $\mu$)}\\
    &\iso p \odot Tq \odot T^2r \odot T(TW(1,\mu)) \odot TW(1,\mu)
    &&\quad\text{(normality of $T$)}\\
    &\iso p \odot T(q \odot Tr \odot TW(1,\mu)) \odot TW(1,\mu)
    &&\quad\text{(strength of $T$)}.
  \end{alignat*}
  Of course, in the unit isomorphisms $\eta$ is used instead of $\mu$.
\end{proof}

\begin{corollary}\label{thm:kmod-equip}
  If $T$ is a horizontally strong monad on an equipment $\mathbb{X}$,
  such that each category $\calH\mathbb{X}(A,B)$ has coequalizers that
  are preserved by $\odot$ on both sides and also preserved by $T$,
  then $\KMod(\mathbb{X},T)$ is also an equipment.
\end{corollary}
\begin{proof}
  The hypotheses ensure that $\Hkl{X}{T}$
  satisfies the conditions of \autoref{thm:comp-mod}.
\end{proof}

\begin{example}
  If $T$ is a cartesian monad on a category \C\ with pullbacks, then
  we have seen that the induced monad on $\SPAN{C}$ is horizontally
  strong; thus $\HKl(\SPAN{C},T)$ is an equipment.  If furthermore \C\
  has coequalizers that are preserved by pullback and by $T$, then
  \autoref{thm:kmod-equip} implies that $\KMod(\SPAN{C},T)$ is also an
  equipment.

  For example, the ``free $M$-set'' monad $(M\times -)$ on \Set\
  preserves coequalizers, so we have an equipment of $M$-graded
  categories.
  However, the ``free monoid'' monad does not preserve coequalizers, and
  the virtual equipment of ordinary multicategories is not an equipment.
\end{example}

\begin{example}\label{eg:fmc-hstrong}
  Let \V\ be a symmetric monoidal category with small coproducts
  preserved by $\otimes$ on both sides, and let $T$ be the extension of
  the ``free monoid'' monad on \Set\ to a monad on $\MAT{V}$ defined in
  \autoref{eg:free-monoid-on-mat}.
  We have already remarked in \autoref{eg:fmc-strong} that $T$ is
  strong, and an easy calculation shows that it is in fact horizontally
  strong.
  Thus, $\Hkl{\MAT{V}}{T}$ is an equipment.
  However, even if \V\ is cocomplete, and in particular has coequalizers
  preserved by $\otimes$ on both sides, these coequalizers will not in
  general be preserved by $T$.
  Thus, the virtual equipment $\KMod(\MAT{V},T)$ of \V-enriched ordinary
  multicategories fails to be an equipment.
  (For $\V=\Set$ we have seen this already in the previous example.)
\end{example}

\begin{example}\label{eg:lawvth-hstrong}
  Let \V\ be a cocomplete symmetric monoidal category with small
  colimits preserved by $\otimes$ on both sides, and let $T$ be the
  ``free symmetric strict monoidal \V-category'' monad on $\PROF{V}$
  from Examples~\ref{eg:free-ssmc} and~\ref{eg:sym-multicat}.
  We remarked in \autoref{eg:fssmc-strong} that $T$ is strong.  In fact,
  it is easily seen to be horizontally strong, so that
  $\Hkl{\PROF{V}}{T}$ is an equipment.
  As in the previous example, however, $T$ fails to preserve
  coequalizers in $\calH(\PROF{V})(A,B)$, so that $\KMod(\PROF{V},T)$ is
  not an equipment even when $\V=\Set$.

  When the horizontal arrows in $\Hkl{\PROF{Set}}{T}$ are identified
  with the \emph{generalized structure types} of~\cite{fghw:gen-species}
  as in \autoref{rmk:structure-types}, their horizontal composites are
  identified with the \emph{substitution} operation on structure types.
  In~\cite{fghw:gen-species} the bicategory $\calH(\Hkl{\PROF{Set}}{T})$
  was constructed in this way from structure types and substitution.
\end{example}

\begin{example}\label{eg:finitary-monads}
  Let $T$ be the ``free category with strictly associative finite
  products'' monad on $\PROF{V}$ from Examples \ref{eg:free-safp} and
  \ref{eg:lawvere-theories}.
  We remarked in \autoref{eg:fssmc-strong} that $T$ is strong if \V\ is
  cartesian monoidal.
  In fact, in this case it can moreover be shown to be horizontally
  strong, so that $\HKl(\PROF{V},T)$ is an equipment.
  
  Now we specialize to $\V=\Set$.
  If $1$ denotes the terminal category, then $T1$ is equivalent to
  $\Set_f^{\op}$, the opposite of the category of finite sets.
  Thus, a profunctor $1\hto T1$ is equivalent to a functor
  $\Set_f\to\Set$, which (since \Set\ is locally finitely presentable)
  is equivalent to a finitary endofunctor of \Set.
  It is then not hard to verify that the equivalence
  \[\calH(\HKl(\PROF{Set},T))(1,1) \simeq [\Set_f,\Set]\]
  is actually an equivalence of monoidal categories, and thus induces an
  equivalence between categories of monoids.
  But a monoid in $\calH(\HKl(\PROF{Set},T))(1,1)$ is a $T$-monoid
  $1\hto T1$, i.e.\ a Lawvere theory; thus we recover the classical
  result of~\cite{lawvere:functsem} that Lawvere theories can be
  identified with \emph{finitary monads} on \Set.

  An analogous argument for the ``free \V-category with finite
  cotensors'' monad on $\PROF{V}$ from Examples \ref{eg:free-coten} and
  \ref{eg:enr-lawvth} reproduces the result of~\cite{power:enr-theories}
  that Lawvere \V-theories can be identified with finitary \V-monads on
  \V.
  In this case, $\HKl(\PROF{V},T)$ seemingly need not be an equipment,
  but at least the multicategory $\HKl(\PROF{V},T)(1,1)$ is a monoidal
  category, precisely because it can be identified with the monoidal
  category of finitary endofunctors of \V\ under composition.
\end{example}

\begin{example}
  Let $\PROF{Set}_{\mathrm{ls}}$ denote the virtual double category
  whose objects are \emph{locally small} categories (that is, large
  categories with small hom-sets), whose vertical arrows are functors,
  and whose horizontal arrows are profunctors taking values in
  \emph{small} sets.  Then there is a monad $T$ on
  $\PROF{Set}_{\mathrm{ls}}$ whose algebras are categories with all
  \emph{small} products, and we have $T1\simeq\Set^{\op}$.  Thus, by
  analogous reasoning to \autoref{eg:finitary-monads}, we see that
  $T$-monoids $1\hto T1$ for this $T$ can be identified with
  \emph{arbitrary} monads on $\Set$.  We can likewise obtain monads on
  any suitable \V\ by using the monad for arbitrary cotensors on
  $\PROF{V}_{\mathrm{ls}}$.  In~\cite{cs:functoriality} we will see
  that by regarding monads as particular generalized multicategories
  in this way, we can recover the monad associated to an operad (as
  originally defined in~\cite{may:goils}) as a particular case of the
  functoriality of generalized multicategories.
\end{example}


\section{Comparisons to previous theories}
\label{sec:comparisons}

We now describe the existing approaches to generalized
multicategories, and show how they compare to our theory.  Most
existing approaches turn out to be instances of our theory, applied to
a particular sort of monad on a particular sort of virtual equipment.
Unsurprisingly, however, often more can be said in such special cases
that is not true in general.  Thus, in each section below we briefly
mention some of the additional results that different authors have
obtained in their particular contexts.

\subsection{Cartesian Monads}\label{section:cartesian}

In order to study and define a type of $n$-category, Leinster
developed a theory of cartesian monads and their associated
multicategories.  This theory was developed in a series of papers,
eventually culminating in his book \cite{leinster:higher-opds}.  Recall
the following from \autoref{eg:cart-mnd}:

\begin{definition}
Let $\mb{E}$ be a category with pullbacks.  A monad $(T, \eta, \mu)$ on
$\mb{E}$ is \textbf{cartesian} if $T$ preserves pullbacks, and all
naturality squares of $\eta$ and $\mu$ are pullbacks.
\end{definition}

Given a cartesian monad, Leinster constructs a bicategory
$\mb{E}_{(T)}$ of \emph{$T$-spans} and defines a \emph{$(\mb{E},
  T)$-multicategory} (or simply a \emph{$T$-multicategory}) to be a
monad in $\mb{E}_{(T)}$.  To compare this to our context, recall that
whenever $\mb{E}$ has pullbacks, $\SPAN{E}$ is an equipment, and any
pullback-preserving monad $T$ on such a $\mb{E}$ extends to a strong
monad on $\SPAN{E}$, which is horizontally strong just when $T$ is
cartesian.  It is then easy to see:

\begin{proposition}
  For a cartesian monad $T$ on $\mb{E}$, Leinster's bicategory
  $\mb{E}_{(T)}$ is isomorphic to $\calH(\Hkl{\SPAN{E}}{T})$.
\end{proposition}

Therefore, Leinster's category of $(\mb{E}, T)$-multicategories
is the vertical category of our $\KMod(\SPAN{E}, T)$.  In particular,
what he calls a \emph{$T$-operad} is a $T$-monoid $1\hto T1$ in $\SPAN{E}$.

Actually, Leinster constructs the whole
virtual double category $\KMod(\SPAN{E}, T)$ (which he calls an
$\mb{fc}$-multicategory) in \cite[\S5.3]{leinster:higher-opds}, and
uses it to define transformations of his
generalized multicategories just as we did in \S\ref{sec:2cat-tmon}.  
(The notes at the end of his \S5.3 also point in the direction of
our \S\S\ref{sec:composites-units} and \ref{sec:pfrbi}.)

Leinster also proves most of the results of our \S\ref{sec:vert-alg}
in his context, as well as the functoriality of the construction
mentioned in \autoref{rmk:functoriality}.
Furthermore, he shows that if $T$ is a ``suitable'' cartesian monad on a
``suitable'' cartesian category $\mb{E}$, then the category of $(\mb{E},
T)$-multicategories is itself monadic over a cartesian category of
``$T$-graphs,'' and this monad is also cartesian.  Thus the process can
be iterated, leading to a definition of the ``opetopes'' used in the
Baez-Dolan definition of weak $n$-categories.
Finally, Leinster also studies \emph{algebras} for generalized operads, which are
closely related to the horizontal arrows in $\KMod(\mathbb{X},T)$.

\subsection{Clubs}
\label{sec:clubs}

Essentially the same theory was developed in~\cite{kelly:clubs} for
the case of generalized operads (generalized multicategories on $1$).
Observe that when $T$
is a cartesian monad on a category $\mb{E}$ with finite limits, so
that $\HKl(\SPAN{E},T)$ is an equipment, then in particular the
category
\[\mb{E}/T1 \simeq \calH(\HKl(\SPAN{E},T))(1,1)\]
inherits a monoidal structure.  It is easy to verify that this
monoidal structure on $\mb{E}/T1$ is the same as that constructed by
Kelly (see the explicit description in
\cite[p.~174--175]{kelly:clubs}).  Kelly defined a \emph{club over
  $T$} to be a monoid for this monoidal structure; thus such clubs can
be identified with $T$-monoids $1\hto T1$ in $\SPAN{E}$.  He also
proved that such clubs are essentially equivalent to cartesian monads
equipped with a ``cartesian map'' to $T$.

(Actually, in~\cite{kelly:clubs} it was assumed that $T$ preserves
\emph{certain} pullbacks rather than all pullbacks.  This suffices to
construct a monoidal structure on $\mb{E}/T1$, though not for $T$ to
define a monad on all of $\SPAN{E}$.)

\subsection{Pseudomonads on Prof}
\label{sec:pseudomonads-prof}

The existing theory which is probably closest to our approach involves
the construction of a Kleisli bicategory from a pseudomonad on a
bicategory such as $\PROF{V}$.  A general theory of multicategories
based on such pseudomonads does not appear to exist in the literature,
but it is implicit
in~\cite{bd:hda3,cheng:opetopic,fghw:gen-species,garner:polycats,ds:lmpo-conv}
among other places.

The general framework is, however, quite simple to state: from a
pseudomonad $T$ on a bicategory $\mathcal{B}$, one
can construct the Kleisli bicategory $\mathcal{B}_T$ and consider
monads in $\mathcal{B}_T$ as a notion of generalized multicategory.
(Leinster's approach is a special case of this for
$\mathcal{B}=\SPAN{C}$, as is Hermida's for a different
bicategory---see \S\ref{sec:cartesian-2-monads}, below.)
The relationship with our theory is that if $T$ is a horizontally
strong monad on an equipment $\mathbb{X}$, it gives rise to a
pseudomonad $\calH(T)$ on the bicategory $\calH(\mathbb{X})$, and
$\calH(\mathbb{X})_{\calH(T)} \simeq \calH(\Hkl{\mathbb{X}}{T})$ so
that the resulting notions of multicategory agree.

In the converse direction, of course not every pseudomonad on
$\calH(\mathbb{X})$ arises from a monad on $\mathbb{X}$ itself, but we
have seen that this is true for most monads relative to which one may
want to define generalized multicategories.  In a few cases, however,
the extension to $\mathbb{X}$ may not be vertically strict,
necessitating the extension of \vEquip\ to a tricategory.

Note that if $T$ is also \emph{co-horizontally strong}, in the sense
that its horizontal dual is horizontally strong, then it also induces
a pseudo-\emph{comonad} \(\hat{\calH}(T)\) on the bicategory
$\calH(\mathbb{X})$.  From this perspective, $T$-monoids can be
identified with \emph{lax $\hat{\calH}(T)$-coalgebras}.  Of course, if
we work in the horizontal dual, then \(\hat{\calH}(T)\) becomes a
pseudomonad and $T$-monoids are its lax algebras.  This is the
terminology used by several authors, including Hermida.  The authors
of~\cite{ds:lmpo-conv} consider the special case when the bicategory
$\mathcal{B}$ is monoidal and $T$ is its free monoid monad, so that
$T$-monoids can be called \emph{lax monoids}.  Such lax monoids can
also be described directly in terms of $\mathcal{B}$, without the need
for cocompleteness hypotheses to ensure that $T$ exists.

\subsection{$(T,\mb{V})$-algebras}
\label{sectiontv}

Following Barr \cite{RelationalAlgebras}, \cite{CommonApproach} started
a series of papers which described the ideas of a ``set-monad with lax
extension to $\mat{V}$'', as well as the $(T, \mb{V})$-algebras
associated to these monads.  Barr's original idea showed that the ``lax
algebras'' of the ultrafilter monad are topological spaces; Clementino
and Tholen's idea extended this further, developing a framework that
eventually included not only topological spaces, but also metric spaces,
approach spaces, and closure spaces.

In the work on $(T,\mb{V})$-algebras, two definitions of set-monad with
lax extension to $\mat{V}$ have been proposed.  The original version was
applicable to all monoidal $\mb{V}$.  However, it failed to capture all
relevant examples, and so a second, slightly different definition was
proposed in \cite{SealCanonical}, which captured further examples.
However, this definition was only applicable when $\mb{V}$ was a
preorder.  As we shall see, Seal's definition turns out to be equivalent
to asking for a monad on $\MAT{V}$ (in the case when $\mb{V}$ is
ordered), while the original definition is equivalent to asking for a
monad on $\MAT{V}$ which is normal.

The original definition, given in \cite{CommonApproach}, was as follows:

\begin{definition}
A \textbf{set monad with lax extension to $\mb{V}$} consists of a monad
$(T, \eta, \mu)$ on $\mb{Set}$, together with a lax functor $T_M$ on
$\mat{V}$ such that:
\begin{packeditem}
	\item $T_M$ is the same as $T$ on objects and functions (viewed as $\mb{V}$-matrices),
	\item the comparisons $(T_Ms)(T_Mr) \to T_M(sr)$ are isomorphisms when $r$ is a function, 
	\item when viewed as transformations on $T_M$, $\eta$ and $\mu$ have op-lax structure.
\end{packeditem}
\end{definition}

In general, however, this definition was found to be too restrictive, as
it didn't allow for examples such as extensions of the powerset monad,
whose algebras would be closure spaces.  To include this type of
example, the requirement that $T_M$ was the same as $T$ on functions
needed to be removed.  Seal's definition, given in \cite{SealCanonical},
was the following:

\begin{definition}
  Suppose that $\mb{V}$ is a monoidal preorder.  A \textbf{set monad
    with lax extension to $\mb{V}$} consists of a monad $(T, \eta, \mu)$
  on $\mb{Set}$, together with a lax functor $T_M$ on $\mat{V}$ which is
  the same as $T$ on objects, and satisfies
	\begin{enumerate}
	\item $Tf \leq T_M f$,
	\item $(Tf)^{\circ} \leq T_M f^{\circ}$.
	\end{enumerate}
\end{definition}
Here, $()^{\circ}$ denotes taking the opposite $\mb{V}$-matrix.  By
\cite[p.~225]{SealCanonical} the conditions imply that if $q$ is a
$\mb{V}$-matrix and $f$ a function, then $T_M(qf) = (T_Mq)(Tf)$ and
$T_M(f^{\circ}q) = (Tf)^{\circ}(T_Mq)$.

In \cite[p.~203]{SealCanonical}, he also shows that when $\mb{V}$ is
completely distributive, $\eta$ and $\mu$ have op-lax structure.
However, this is not a priori required in his definition.  If, however,
we include this axiom in his definition, then his notion of set monad
with lax extension is equivalent to giving our notion of a monad on the
virtual double category $\MAT{V}$.

\begin{proposition}
  Suppose that $\mb{V}$ is a monoidal preorder.  If $T$ is a set monad
  with lax extension $T_M$ (in the sense of Seal) for which $\eta$ and
  $\mu$ are op-lax, then we can define a monad $T$ on $\MAT{V}$ which is
  $T$ on vertical arrows, and $T_M$ on horizontal arrows.  Conversely,
  given a monad $T$ on $\MAT{V}$, we can define a set monad with lax
  extension which is $T$ on functions, and uses the horizontal action of
  $T$ to define $T_M$.
\end{proposition}
\begin{proof}
Suppose that we have a set monad with lax extension to $\mb{V}$, in the
second sense given above.  Define a functor $T$ on the double category
$\MAT{V}$, which is $T$ on vertical arrows, and $T_M$ on horizontal
arrows.  Using the $\eta$ and $\mu$, we get all of the necessary data
for a monad on $\MAT{V}$, with the exception of checking that 
\[
	\bfig
	\square/`>`>`/<400,350>[X`Y`W`Z;`f`g`]
	\morphism(0,350)/@{>}|-*@{|}/<400,0>[X`Y;p]
	\morphism(0,0)|b|/@{>}|-*@{|}/<400,0>[W`Z;q]
	\place(200,225)[\twoar(0,-1)]
	
	\place(800,175)[\mbox{implies}]
	
	\square(1200,0)/`>`>`/<400,350>[TX`TY`TW`TZ;`Tf`Tg`]
	\morphism(1200,350)/@{>}|-*@{|}/<400,0>[TX`TY;Tp]
	\morphism(1200,0)|b|/@{>}|-*@{|}/<400,0>[TW`TZ;Tq]
	\place(1400,225)[\twoar(0,-1)]
	
	\efig
\]
This is equivalent to checking that $p \leq g^{\circ} qf \Rightarrow
T_M(p) \leq (Tg)^{\circ} T_Mq (Tf)$.  But this is easy to check by using
the two results given after Seal's definition.
\[ T_M(p) \leq T_M(g^{\circ}qf) = (Tg)^{\circ} T_M(qf) = (Tg)^{\circ}(T_Mq)(Tf) \]

Conversely, suppose that we have a monad $T$ on $\MAT{V}$.  We would
like to define a lax extension of $T$ (considered as a $\Set$-monad) to
$\mat{V}$.  Define $T_M$ on matrices as for $T$.  The only conditions we
need to check are $Tf \leq T_M f$ and $(Tf)^{\circ} \leq T_M f^{\circ}$.
To show the first is equivalent to showing that $TB(1,{Tf}) \leq
T(B(1,f))$.  To show this, recall that we have a cartesian cell
\[
	\bfig
	\node a(0,350)[A]
	\node b(400,350)[B]
	\node c(0,0)[B]
	\node d(400,0)[B]
	\arrow/@{>}|-*@{|}/[a`b;B(1,f)]
	\arrow|b|/@{>}|-*@{|}/[c`d;U_B]
	\arrow[a`c;f]
	\arrow/=/[b`d;]
	\place(200,225)[\twoar(0,-1)]
	\efig
\]
Moreover, since $T$ is a functor, it preserves cartesian cells
(\autoref{thm:func-pres-restr}), and so
\[
	\bfig
	\node a(0,350)[TA]
	\node b(600,350)[TB]
	\node c(0,0)[TB]
	\node d(600,0)[TB]
	\arrow/@{>}|-*@{|}/[a`b;T(B(1,f))]
	\arrow|b|/@{>}|-*@{|}/[c`d;T(U_B)]
	\arrow[a`c;Tf]
	\arrow/=/[b`d;]
	\place(300,225)[\twoar(0,-1)]
	\efig
\]
is also cartesian.  We can thus factor the cell
\[
	\bfig
	\node a(0,350)[TA]
	\node b(600,350)[TB]
	\node c(0,0)[TB]
	\node d(600,0)[TB]
	\node e(0,-350)[TB]
	\node f(600,-350)[TB]
	
	\arrow/@{>}|-*@{|}/[a`b;TB(1,Tf)]
	\arrow|b|/@{>}|-*@{|}/[c`d;U_{TB}]
	\arrow[a`c;Tf]
	\arrow/=/[b`d;]
	\arrow/=/[c`e;]
	\arrow/=/[d`f;]
	\arrow|b|/@{>}|-*@{|}/[e`f;T(U_B)]
	
	\place(300,225)[\twoar(0,-1)]
	\place(300,-150)[\twoar(0,-1)]
	\efig
\]
through it to get a cell
\[
	\bfig
	\node a(0,350)[TA]
	\node b(600,350)[TB]
	\node c(0,0)[TA]
	\node d(600,0)[TB]
	\arrow/@{>}|-*@{|}/[a`b;TB(1,Tf)]
	\arrow|b|/@{>}|-*@{|}/[c`d;T(B(1,f))]
	\arrow/=/[a`c;]
	\arrow/=/[b`d;]
	\place(300,225)[\twoar(0,-1)]
	\efig
\]
as required.  The second inequality follows similarly, using the cartesian cell
\[
	\bfig
	\node a(0,350)[B]
	\node b(400,350)[A]
	\node c(0,0)[B]
	\node d(400,0)[B]
	\arrow/@{>}|-*@{|}/[a`b;B(f,1)]
	\arrow|b|/@{>}|-*@{|}/[c`d;U_B]
	\arrow|r|[b`d;f]
	\arrow/=/[a`c;]
	\place(200,225)[\twoar(0,-1)]
	\efig
\]
Thus, a monad on $\MAT{V}$ defines a set-monad with lax extension to $\mat{V}$.
\end{proof}

For general $\mb{V}$, we can also use the above correspondence to
recover the first notion of notion of set-monad with lax extension: they
are the monads which are normal.

\begin{proposition}
Using the above correspondence, we get a set-monad with lax extension in
the first sense if and only if the monad $T$ on $\MAT{V}$ is normal.
\end{proposition}
\begin{proof}
Suppose we have a set-monad $T$ with lax extension $T_M$ in the first
sense.  Then we have $T_Mf \cong Tf$ for all functions $f$, and we get a
monad on $\MAT{V}$.  Moreover,  we also have $T(B(1,f)) \cong TB(1,Tf)$.
In particular, we have $T(A(1,1)) \cong TA(1,1)$.  But $A(1,1) \cong
U_A$, so we have $T(U_A) \cong U_{TA}$.  Thus $T$ is normal.

Conversely, suppose that we have a monad $T$ on $\MAT{V}$ which is
normal.  We would like to show that $T(B(1,f)) \cong TB(1,Tf)$.  Using
the factoring as for the above proposition, we get a cell in one
direction.  To get the other direction, we factor
\[
	\bfig
	\square/`=`=`/<600,350>[TB`TB`TB`TB;```]
	\morphism(0,700)/@{>}|-*@{|}/<600,0>[TA`TB;T(B(1,f))]
	\morphism(0,350)/@{>}|-*@{|}/<600,0>[TB`TB;T(U_B)]
	\place(300,225)[\twoar(0,-1)]

	\square(0,350)/`>`=`/<600,350>[TA`TB`TB`TB;`Tf``]
	\morphism|b|/@{>}|-*@{|}/<600,0>[TB`TB;U_{TB}]
	\place(300,600)[\twoar(0,-1)]

	\place(1200,350)[\mbox{through}]

	\square(1700,175)/`>`=`/<600,350>[TA`TB`TB`TB;`Tf``]
	\morphism(1700,525)/@{>}|-*@{|}/<600,0>[TA`TB;TB(1,Tf)]
	\morphism(1700,175)|b|/@{>}|-*@{|}/<600,0>[TB`TB;U_{TB}]
	\place(2000,375)[\twoar(0,-1)]
		
	\efig
\]
(note that the bottom cell on the left exists by normality of $T$).  The
composites of the two cells are identities by the universal property of
the cartesian cells, and so we have $T(B(1,f)) \cong TB(1,Tf)$, as
required.  We also need to check the condition that the comparison cell
be an isomorphism when $r$ is a function.  However, this is equivalent
to asking that $T(r(1,s)) \cong Tr(1,Ts)$, and this follows from
\autoref{thm:func-pres-restr}.
\end{proof}

For a set-monad with lax extension $T$, the category of
$(T,\mb{V})$-algebras that Tholen, Clementino and Seal define is exactly
the vertical category of the virtual double category $\KMod(\MAT{V},
T)$.  They also describe $(T, \mb{V})$-modules, and these are the
horizontal arrows of $\KMod(\MAT{V}, T)$.

In addition to providing a more conceptual explanation of the notion of
lax extension, and a way to compare $(T,\mb V)$-algebras with other
notions of generalized multicategory, our general framework improves the
theory of $(T, \mb{V}$)-algebras in two ways.  Firstly, it gives a
context in which the horizontal Kleisli construction makes sense; such a
construction has been recognized as desirable (see, for example,
\cite[p.~7]{tholen:coursenotes}), but is impossible using only
bicategories since the monads used in this case are not horizontally
strong.  Secondly, it provides a general reason for the observation of
\cite[p.~15]{tholen:coursenotes} that any set-monad with lax extension
to $\mat{V}$ can also be extended to a monad on $\cat{V}$ with a lax
extension to $\prof{V}$ (we simply apply the 2-functor $\MOD$).
 
There are, however, other special aspects of the theory of
$(T,\mb{V})$-algebras which we have not discussed.  One is that the
category of $(T,\mb{V})$-algebras is generally \emph{topological} over
\Set, and has many other similar formal properties.  Another is the use
of the ``pro'' construction found in \cite{OneSetting}.  This is useful
to describe additional topological structures; for example, monoids in
$\pro\Rel$ are quasi-uniform spaces.  In general, given a virtual double
category $\mathbb{X}$, one can define a new virtual double category
$\pro\mathbb{X}$, and go on to describe ``pro-generalized
multicategories''.  Further discussion of this, however, awaits a future
paper.

\subsection{Non-cartesian monads}
\label{sec:non-cartesian-monads}

The earliest work on generalized multicategories was by Burroni in
\cite{burroni:t-cats}.  His framework is very similar to Leinster's (see
\S\ref{section:cartesian})
except that he requires nothing at all about the monad $T$, not even
that it preserve pullbacks (although the category $\mb{E}$ must still
have pullbacks).  This level of generality does not fit into our
existing framework, since if $T$ does not preserve pullbacks then it
does not induce a functor on $\SPAN{E}$ of the sort we have defined.
However, it does induce an \emph{oplax functor} between pseudo double
categories.

The simplest way for us to define an oplax functor is to say it is a
functor between pseudo double categories regarded as co-virtual double
categories.  Pseudo double categories, oplax functors, and
transformations form a 2-category, and we define an \emph{oplax monad}
on a pseudo double category to be a monad in this 2-category.  Since
\emph{any} functor on a category $\mb{E}$ with pullbacks induces an
oplax functor on $\SPAN{E}$, any monad on such an $\mb{E}$ induces an
oplax monad on $\SPAN{E}$.  Moreover, we can extend our framework to
deal with oplax functors as follows.

\begin{definition}
  If $T$ is an oplax monad on a pseudo double category $\mathbb{X}$,
  the \textbf{horizontal Kleisli virtual double category} $\Hkl{X}{T}$
  of $T$ is defined as follows.
  \begin{itemize}
  \item Its objects, vertical and horizontal arrows, and cells with
    nullary source are defined as when $T$ is a lax functor.
  \item A cell
    \[
    \bfig
    \square/`>`>`/<2000,250>[X_0`X_n`Y_0`Y_1;`f`g`]
    \morphism|b|/@{>}|-*@{|}/<2000,0>[Y_0`Y_1;q]
    \morphism(0,250)/@{>}|-*@{|}/<500,0>[X_0`X_1;p_1]
    \morphism(500,250)/@{>}|-*@{|}/<500,0>[X_1`X_2;p_2]
    \morphism(1000,250)/@{>}|-*@{|}/<500,0>[X_2`\cdots;p_3]
    \morphism(1500,250)/@{>}|-*@{|}/<500,0>[\cdots`X_n;p_n]
    \place(950,150)[\twoar(0,-1)]
    \place(1050,100)[\scriptstyle\alpha]
    \efig
    \]
    in $\Hkl{X}{T}$ is a cell
    \[\vcenter{\xymatrix@C=15pc{X_0
        \ar[r]|-@{|}^-{p_1 \odot T(p_2 \odot T(\cdots T(p_{n-1}\odot Tp_n)\cdots))} \ar[d]_f
        \ar@{}[dr]|{\twoar(0,-1)} & T^n X_n\ar[d]^{Tg\circ \mu^n}\\
        Y_0\ar[r]|-@{|}_-{q} & TY_1}}
    \]
    in $\mathbb{X}$.
  \item Composition is defined using the multiplication, unit, and
    oplax structure of $T$.  For example, the composite of
    \[\vcenter{\xymatrix{
        \ar[r]|-@{|}^{}\ar[d] \ar@{}[drr]|{{\textstyle\Downarrow} \alpha} &
        \ar[r]|-@{|}^{} &
        \ar[r]|-@{|}^{}\ar[d] \ar@{}[drr]|{{\textstyle\Downarrow} \beta} &
        \ar[r]|-@{|}^{} &
        \ar[r]|-@{|}^{}\ar[d] \ar@{}[drr]|{{\textstyle\Downarrow} \gamma} &
        \ar[r]|-@{|}^{} &
        \ar[d]\\
        \ar[rr]|-@{|}^{}\ar[d] \ar@{}[drrrrrr]|{{\textstyle\Downarrow}\delta} &&
        \ar[rr]|-@{|}^{}&&
        \ar[rr]|-@{|}^{}&&
        \ar[d]\\
        \ar[rrrrrr]|-@{|}_{} &&&&&&
        }}
      \]
      is given by the composite
      \[\delta \circ ( \alpha \odot T(\beta \odot T\gamma)) \circ
      (1 \odot \mu) \circ (1\odot T(1\odot \mu)) \circ T_\odot
      \]
      in $\mathbb{X}, $where $T_\odot$ is a composite of the oplax
      structure maps of $T$:
      \[p_1 \odot T(p_2 \odot T(p_3 \odot T(p_4\odot T(p_5\odot Tp_6))))
      \to
      p_1\odot Tp_2 \odot T^2(p_3\odot Tp_4 \odot T^2(p_5\odot Tp_6))
      \]
  \end{itemize}
\end{definition}

Note that when $T$ is a strong functor, both definitions of
$\Hkl{X}{T}$ make sense; however, in this case, they are equivalent.

\begin{definition}
  If $T$ is an oplax monad on a pseudo double category $\mathbb{X}$,
  then a \textbf{$T$-monoid} is a monoid in $\Hkl{X}{T}$, and we write
  $\KMod(\mathbb{X},T) = \Mod(\Hkl{X}{T})$.
\end{definition}

To compare this to Burroni's definition, note that $\Hkl{X}{T}$
clearly has weak composites, and hence is an oplax double category
(see \autoref{eg:vdblcat}).
Burroni works instead with what he calls a \emph{pseudo-category}, and
what we would probably call a \emph{lax-biased bicategory}: a
bicategory-like structure with units and binary composites and
noninvertible comparison maps
\begin{align*}
  p &\to p\odot U_A\\
  p &\to U_B \odot p\\
  p \odot(q\odot r) &\to (p\odot q)\odot r.
\end{align*}
satisfying suitable axioms.  In fact, of course, Burroni's
``pseudo-category of $T$-spans'' extends to a lax-biased \emph{double}
category.  Moreover, any lax-biased double category defines an oplax
double category (and hence a virtual double category), if we take the
$n$-ary composite to be
\[\langle p_1\odot \dots\odot p_n \rangle =
p_1 \odot (p_2\odot \cdots (p_{n-1}\odot p_n)\cdots).
\]
(Burroni points this out as well; see~\cite[p.~66,
example~3]{burroni:t-cats}.  He refers to virtual double categories as
simply ``multicat\'egories''.)
In this way, Burroni's $\mathrm{Sp}(T)$ is identified with our
$\calH(\Hkl{\SPAN{E}}{T})$, and thus his ``$T$-categories'' can be
identified with our $T$-monoids.  However, he must move
outside the bicategory (or ``pseudo-category'') framework to define
functors of $T$-categories,
whereas they emerge naturally from our setup.

Burroni also constructs the left adjoint $F_T$ from
\autoref{rmk:free-models}, and proves the functoriality of his
construction under lax morphisms of monads.  Working with his
language,~\cite{penon:tcats-rep} gives a version of the
representability results from \S\ref{sec:vert-alg}.

\subsection{Cartesian 2-monads}
\label{sec:cartesian-2-monads}

We have described the horizontal arrows $A\hto B$ in $\PROF{Set}$ as
profunctors, i.e.\ functors $B^{op}\times A\to \Set$, but it is
well-known that such functors can equivalently be described by
\emph{two-sided discrete fibrations}, and that this notion can be
internalized to a sufficiently well-behaved 2-category.
\cite{HermidaCoherent} develops a theory of generalized
multicategories in such a context.

Let $\calK$ be a finitely complete 2-category
(see~\cite{street:catval-2lim}).  Since its underlying ordinary
category $\calK_0$ has pullbacks, we can form the virtual equipment
$\iPROF{\calK_0}$.  Now, as observed in \cite{street:fib-yoneda-2cat},
for any object $A\in\calK$ we have an internal category $\Phi A$ in
$\calK$, defined by $(\Phi A)_0=A$ and $(\Phi A)_1 = A^\bbtwo$ (the
cotensor with the arrow category \bbtwo).  Similarly, every morphism
$f\maps A\to B$ defines an internal functor $\Phi f \maps \Phi
A\to\Phi B$, and the same for 2-cells; thus we have a 2-functor $\calK
\to \calV(\iPROF{\calK_0})$.  Moreover, this 2-functor is locally full
and faithful, i.e.\ 2-cells $f\to g$ are in bijection with internal
natural transformations $\Phi f \to\Phi g$.

We define an internal profunctor $\Phi A \hto^H \Phi B$ to be a
\textbf{discrete fibration} from $A$ to $B$ if the object $H\in
\calK/(A\times B)$ is internally discrete, i.e.\ $(\calK/(A\times
B))(C,H)$ is a discrete category for any $C\in \calK/(A\times B)$.  We
define $\DFib(\calK)$ to be the sub-virtual double category of
$\iPROF{\calK_0}$ determined by
\begin{packeditem}
\item The internal categories of the form $\Phi A$,
\item The internal functors of the form $\Phi f$,
\item The internal profunctors which are discrete fibrations, and
\item All cells between these.
\end{packeditem}
Since the pullback of a discrete fibration is a discrete fibration,
and $U_{\Phi A}$ is a discrete fibration, $\DFib(\calK)$ is a virtual
equipment.  Our remarks above show that $\calK \simeq
\calV(\DFib(\calK))$.  Under suitable conditions on \calK, discrete
fibrations can be composed, so that $\DFib(\calK)$ becomes an
equipment.  (Several authors have tried to isolate these conditions,
with varying degrees of success; in addition to \cite{HermidaCoherent}
see \cite{street:fibi} and \cite{cjsv:modulated-bicats}.
Our approach sidesteps this issue completely.)

Now suppose that $F\maps \calK\to\calL$ is a 2-functor that preserves
pullbacks and comma objects.  Then it preserves internal categories,
profunctors, the $\Phi$ construction, and discrete fibrations, so it
induces a normal functor $\DFib(F)\maps
\DFib(\calK)\to\DFib(\calL)$. Likewise, any 2-natural transformation
$F\to G$ induces a transformation $\DFib(F)\to\DFib(G)$, so we have a
2-functor $\DFib$ from finitely complete 2-categories to $\vEquip$.
In particular, any 2-monad $T$ on \calK\ whose functor part preserves
pullbacks and comma objects induces a normal monad on $\DFib(\calK)$,
so we can talk about $\DFib(T)$-monoids.

This is basically the context of \cite{HermidaCoherent}, except that,
like most other authors, he works only with bicategories.  Thus, he assumes that
\calK\ has the structure required to compose discrete fibrations, and
that moreover $T$ preserves this structure and that $\mu$ and $\eta$
are cartesian transformations.  This ensures that $\DFib(T)$ is
horizontally strong, so that $\Hkl{\DFib(\calK)}{\DFib(T)}$ is an
equipment.  Under these hypotheses, we have:

\begin{theorem}
  The 2-category $\nKMon(\DFib(\calK),\DFib(T))$ of normalized
  $\DFib(T)$-monoids is isomorphic to the 2-category
  $\mbox{Lax-Bimod}(T)\mbox{-alg}$ defined in \cite[4.3 and
  4.4]{HermidaCoherent}.
\end{theorem}

Our \autoref{thm:oplax-alg} is also a generalization of results of
\cite{HermidaCoherent}.  Hermida proves furthermore that under his
hypotheses, the left adjoint $F_T$ from \autoref{rmk:free-models}
exists, the adjunction is monadic when restricted to pseudo
$T$-algebras, and the induced monad $T_\dashv$ on
$\mbox{Lax-Bimod}(T)\mbox{-alg}$ is lax-idempotent
(see~\cite{kl:property-like}).  In~\cite{cs:laxidem} we will show that
an analogous result is true for any monad $T$ on an equipment
$\mathbb{X}$ satisfying suitable cocompleteness conditions.

\subsection{Monoidal pseudo algebras}
\label{sec:mon-psalg}

In~\cite{weber:operads}, Weber gives a definition of generalized
operads enriched in \emph{monoidal pseudo algebras}.  More precisely,
for any 2-monad $T$ on a 2-category \calK\ with finite products, and
any pseudo $T$-algebra $A$ which is also a pseudomonoid in a
compatible way, he defines a notion of \emph{$T$-operad in $A$}.
A general description of the relationship of this theory to ours would
take us too far afield, so we will remark only briefly on how such a
comparison should go.

For any pseudomonoid $A$ in a 2-category \calK\ with finite products,
there is a virtual equipment $A\Mat$ defined as follows.  Its objects
and vertical arrows are the objects and arrows of \calK.  A horizontal
arrow from $X\hto Y$ is a morphism $X\times Y \to^p A$ in \calK, and a
cell
\[\vcenter{\xymatrix{
    X_0\ar[r]|-@{|}^{p_1}\ar[d]_f \ar@{}[drrr]|{\twoar(0,-1)} &
    X_1\ar@{-->}[r] &
    X_{n-1}\ar[r]|-@{|}^-{p_n} &
    X_n\ar[d]^g\\
    Y_0\ar[rrr]|-@{|}_{q} &&& Y_1
  }}
\]
is a 2-cell
\[\bfig\scalefactor{.8}
\Atriangle[(X_0\times X_1)\times (X_1\times X_2) \times\dots \times
(X_{n-1}\times X_n)`Y_0\times Y_1`A.;(f,g)`(p_1\otimes\dots\otimes p_n)`q]
\place(500,250)[\twoar(-2,-1)]
\efig
\]
in \calK, where $\otimes\maps A\times A\to A$ is the multiplication of
the pseudomonoid $A$ (if $n=0$ we use its unit instead).

Now, if $A$ is additionally a pseudo $T$-algebra in a compatible way,
we might hope to be able to extend $T$ to a monad on $A\Mat$.
However, given a horizontal arrow $X\hto^p Y$, from $X\times Y \to^p
A$ we can form the composite
\[T(X\times Y) \to^{Tp} TA \to^a A,
\]
but this is not yet a horizontal arrow $TX \hto TY$.  If we assume
that $A$ admits well-behaved \emph{left (Kan) extensions}, then we can
define $TX\hto^{Tp} TY$ to be the extension of the above composite
along $T(X\times Y)\to TX \times TY$.  We can then construct
$\HKl(A\Mat,T)$ and $\KMod(A\Mat,T)$ as usual.  Moreover, we can give
an equivalent characterization of $\HKl(A\Mat,T)$ which is valid even
in the absence of left extensions: a horizontal arrow $X\hto Y$ is a
morphism $X\times TY \to^p A$ in \calK, and a cell
\[\vcenter{\xymatrix{
    X_0\ar[r]|-@{|}^{p_1}\ar[d]_f \ar@{}[drrr]|{{\twoar(0,-1)}} &
    X_1\ar@{-->}[r] &
    X_{n-1}\ar[r]|-@{|}^-{p_n} &
    X_n\ar[d]^g\\
    Y_0\ar[rrr]|-@{|}_{q} &&& Y_1
  }}
\]
is a 2-cell
\[\vcenter{\xymatrix{
    X_0\times T(X_1\times \dots T(X_{n-1}\times TX_n)\dots) \ar[d]_\pi
    \ar[dr]^(.7){p_1\ten\dots\ten p_n}\\
    X_0 \times T^nX_n \ar[d]_{1\times \mu^{n-1}} \rtwocell\omit &  A\\
    X_0\times TX_n \ar[d] _{f\times Tg}\\
    Y_0 \times TY_1 \ar[uur]_{q}
  }}
\]
in \calK.
Thus, we obtain a notion of $T$-monoid in $A$ for any monoidal pseudo
algebra $A$.  Weber only considers the case of \emph{operads}, rather
than more general \emph{multicategories}, but it is easy to verify
that his $T$-operads in $A$ coincide with those $T$-monoids in $A$
whose underlying object in \calK\ is the terminal object $1$.

Actually, there is a good reason that Weber considers only operads:
$T$-monoids $X_0\hto^X TX_0$ in this context for which $X_0$ is not
discrete, or at least a groupoid, are not very familiar objects.  In
familiar cases such as $\calK=\Cat$, we would obtain familiar types of
$A$-enriched multicategory only by taking the horizontal arrows $X\hto
Y$ in $A\Mat$ to be morphisms $X\times Y^{\op}\to A$, rather than
$X\times Y\to A$.  To put this in a general context, however, requires
a 2-category \calK\ in which ``opposites'' make sense, such as the
``2-toposes'' of~\cite{weber:2toposes}.

\bibliographystyle{alpha}

\begin{thebibliography}{FGHW08}

\bibitem[Bar70]{RelationalAlgebras}
Michael Barr.
\newblock Relational algebras.
\newblock {\em Lectures Notes in Mathematics 137: Reports of the Midwest
  Category Seminar IV}, pages 39--55, 1970.

\bibitem[Bat98]{batanin:monglob}
M.~A. Batanin.
\newblock Monoidal globular categories as a natural environment for the theory
  of weak {$n$}-categories.
\newblock {\em Adv. Math.}, 136(1):39--103, 1998.

\bibitem[BD98]{bd:hda3}
John~C. Baez and James Dolan.
\newblock Higher-dimensional algebra. {III}. {$n$}-categories and the algebra
  of opetopes.
\newblock {\em Adv. Math.}, 135(2):145--206, 1998.

\bibitem[BD01]{bd:fin-feyn}
John~C. Baez and James Dolan.
\newblock From finite sets to {F}eynman diagrams.
\newblock In {\em Mathematics unlimited---2001 and beyond}, pages 29--50.
  Springer, Berlin, 2001.

\bibitem[Bec69]{beck:dl}
Jon Beck.
\newblock Distributive laws.
\newblock In {\em Sem. on Triples and Categorical Homology Theory (ETH,
  Z\"urich, 1966/67)}, pages 119--140. Springer, Berlin, 1969.

\bibitem[B{\'{e}}n67]{benabou:bicats}
Jean B{\'{e}}nabou.
\newblock Introduction to bicategories.
\newblock In {\em Reports of the Midwest Category Seminar}, pages 1--77.
  Springer, 1967.

\bibitem[BM03]{brgmoer:operads}
Clemens Berger and Ieke Moerdijk.
\newblock Axiomatic homotopy theory for operads.
\newblock {\em Comment. Math. Helv.}, 78(4):805--831, 2003.

\bibitem[BS76]{bs:dblgpd-xedmod}
Ronald Brown and Christopher~B. Spencer.
\newblock Double groupoids and crossed modules.
\newblock {\em Cahiers Topologie G\'eom. Diff\'erentielle}, 17(4):343--362,
  1976.

\bibitem[Bur71]{burroni:t-cats}
Albert Burroni.
\newblock {$T$}-cat\'egories (cat\'egories dans un triple).
\newblock {\em Cahiers Topologie G\'eom. Diff\'erentielle}, 12:215--321, 1971.

\bibitem[BV73]{bv:htpy-invar}
J.~M. Boardman and R.~M. Vogt.
\newblock {\em Homotopy invariant algebraic structures on topological spaces}.
\newblock Lecture Notes in Mathematics, Vol. 347. Springer-Verlag, Berlin,
  1973.

\bibitem[Che04]{cheng:opetopic}
Eugenia Cheng.
\newblock Weak n-categories: opetopic and multitopic foundations.
\newblock {\em Journal of Pure and Applied Algebra}, 186(2):109 -- 137, 2004.

\bibitem[CHT04]{OneSetting}
M.~Clementino, D.~Hofmann, and W.~Tholen.
\newblock One setting for all: metric, topology, uniformity, approach
  structure.
\newblock {\em Applied Categorical Structures}, 12:127--154, 2004.

\bibitem[CJSV94]{cjsv:modulated-bicats}
Aurelio Carboni, Scott Johnson, Ross Street, and Dominic Verity.
\newblock Modulated bicategories.
\newblock {\em J. Pure Appl. Algebra}, 94(3):229--282, 1994.

\bibitem[CS10a]{cs:functoriality}
G.~S.~H. Cruttwell and Michael~A.\ Shulman.
\newblock On the functoriality of generalized multicategories.
\newblock In preparation, 2010.

\bibitem[CS10b]{cs:laxidem}
G.~S.~H. Cruttwell and Michael~A.\ Shulman.
\newblock Representability of generalized multicategories.
\newblock In preparation, 2010.

\bibitem[CT03]{CommonApproach}
M.~Clementino and W.~Tholen.
\newblock Metric, topology, and multicategory - a common approach.
\newblock {\em J. Pure and Applied Algebra}, 153:13--47, 2003.

\bibitem[DL07]{dl:lim-smallfr}
Brian~J. Day and Stephen Lack.
\newblock Limits of small functors.
\newblock {\em J. Pure Appl. Algebra}, 210(3):651--663, 2007.

\bibitem[DPP06]{dpp:paths-dbl}
R.~Dawson, R.~Par{\'e}, and D.~A. Pronk.
\newblock Paths in double categories.
\newblock {\em Theory Appl. Categ.}, 16:No. 18, 460--521 (electronic), 2006.

\bibitem[DPP07]{dpp:spans}
R.~Dawson, R.~Par{\'e}, and D.~A. Pronk.
\newblock Spans for 2-categories.
\newblock In preparation; slides accessed May 2007 at
  \url{http://www.mathstat.dal.ca/~selinger/ct2006/slides/CT06-Pare.pdf}, 2007.

\bibitem[DS97]{ds:monbi-hopfagbd}
Brian Day and Ross Street.
\newblock Monoidal bicategories and {H}opf algebroids.
\newblock {\em Adv. Math.}, 129(1):99--157, 1997.

\bibitem[DS03]{ds:lmpo-conv}
Brian Day and Ross Street.
\newblock Lax monoids, pseudo-operads, and convolution.
\newblock In {\em Diagrammatic morphisms and applications ({S}an {F}rancisco,
  {CA}, 2000)}, volume 318 of {\em Contemp. Math.}, pages 75--96. Amer. Math.
  Soc., Providence, RI, 2003.

\bibitem[EK66]{ek:closed-cats}
Samuel Eilenberg and G.~Max Kelly.
\newblock Closed categories.
\newblock In {\em Proc. Conf. Categorical Algebra (La Jolla, Calif., 1965)},
  pages 421--562. Springer, New York, 1966.

\bibitem[EM06]{em:rma-infloop}
A.~D. Elmendorf and M.~A. Mandell.
\newblock Rings, modules, and algebras in infinite loop space theory.
\newblock {\em Adv. Math.}, 205(1):163--228, 2006.

\bibitem[FGHW08]{fghw:gen-species}
M.~Fiore, N.~Gambino, M.~Hyland, and G.~Winskel.
\newblock The {C}artesian closed bicategory of generalised species of
  structures.
\newblock {\em J. Lond. Math. Soc. (2)}, 77(1):203--220, 2008.

\bibitem[Gar08]{garner:polycats}
Richard Garner.
\newblock Polycategories via pseudo-distributive laws.
\newblock {\em Adv. Math.}, 218(3):781--827, 2008.

\bibitem[GP99]{gp:double-limits}
Marco Grandis and Robert Par{\'e}.
\newblock Limits in double categories.
\newblock {\em Cahiers Topologie G\'eom. Diff\'erentielle Cat\'eg.},
  XL(3):162--220, 1999.

\bibitem[GP04]{gp:double-adjoints}
Marco Grandis and Robert Par{\'e}.
\newblock Adjoints for double categories.
\newblock {\em Cah. Topol. G\'eom. Diff\'er. Cat\'eg.}, 45(3):193--240, 2004.

\bibitem[Her00]{HermidaRepresentable}
Claudio Hermida.
\newblock Representable multicategories.
\newblock {\em Advances in Mathematics}, 151:164--225, 2000.

\bibitem[Her01]{HermidaCoherent}
Claudio Hermida.
\newblock From coherent structures to universal properties.
\newblock {\em J. Pure and Applied Algebra}, 165:7--61, 2001.

\bibitem[Joh02]{ptj:elephant1}
Peter~T. Johnstone.
\newblock {\em Sketches of an Elephant: A Topos Theory Compendium: Volume 1}.
\newblock Number~43 in Oxford Logic Guides. Oxford Science Publications, 2002.

\bibitem[Joy81]{joyal:species}
Andr{\'e} Joyal.
\newblock Une th\'eorie combinatoire des s\'eries formelles.
\newblock {\em Adv. in Math.}, 42(1):1--82, 1981.

\bibitem[Joy86]{joyal:funcan-species}
Andr{\'e} Joyal.
\newblock Foncteurs analytiques et esp\`eces de structures.
\newblock In {\em Combinatoire \'enum\'erative ({M}ontreal, {Q}ue.,
  1985/{Q}uebec, {Q}ue., 1985)}, volume 1234 of {\em Lecture Notes in Math.},
  pages 126--159. Springer, Berlin, 1986.

\bibitem[Kel72a]{kelly:abst-coh}
G.~M. Kelly.
\newblock An abstract approach to coherence.
\newblock In {\em Coherence in categories}, pages 106--147. Lecture Notes in
  Math., Vol. 281. Springer, Berlin, 1972.

\bibitem[Kel72b]{kelly:mv-funct-calc}
G.~M. Kelly.
\newblock Many-variable functorial calculus. {I}.
\newblock In {\em Coherence in categories}, pages 66--105. Lecture Notes in
  Math., Vol. 281. Springer, Berlin, 1972.

\bibitem[Kel82]{kelly:enr-lfp}
G.~M. Kelly.
\newblock Structures defined by finite limits in the enriched context. {I}.
\newblock {\em Cahiers Topologie G\'eom. Diff\'erentielle}, 23(1):3--42, 1982.
\newblock Third Colloquium on Categories, Part VI (Amiens, 1980).

\bibitem[Kel92]{kelly:clubs}
G.~M. Kelly.
\newblock On clubs and data-type constructors.
\newblock In {\em Applications of categories in computer science ({D}urham,
  1991)}, volume 177 of {\em London Math. Soc. Lecture Note Ser.}, pages
  163--190. Cambridge Univ. Press, Cambridge, 1992.

\bibitem[Kel05]{kelly:operads}
G.~M. Kelly.
\newblock On the operads of {J}. {P}. {M}ay.
\newblock {\em Repr. Theory Appl. Categ.}, (13):1--13 (electronic), 2005.

\bibitem[KL97]{kl:property-like}
G.~M. Kelly and Stephen Lack.
\newblock On property-like structures.
\newblock {\em Theory Appl. Categ.}, 3(9):213--250, 1997.

\bibitem[KM95]{km:oamm}
Igor K{\v{r}}{\'{\i}}{\v{z}} and J.~P. May.
\newblock Operads, algebras, modules and motives.
\newblock {\em Ast\'erisque}, (233):iv+145pp, 1995.

\bibitem[Lac]{lack:icons}
Stephen Lack.
\newblock Icons.
\newblock arXiv:0711.4657.

\bibitem[Lam72]{lambek:dedsys-iii}
Joachim Lambek.
\newblock Deductive systems and categories. {III}. {C}artesian closed
  categories, intuitionist propositional calculus, and combinatory logic.
\newblock In {\em Toposes, algebraic geometry and logic ({C}onf., {D}alhousie
  {U}niv., {H}alifax, {N}.{S}., 1971)}, pages 57--82. Lecture Notes in Math.,
  Vol. 274. Springer, Berlin, 1972.

\bibitem[Lam89]{lambek:multi-revis}
J.~Lambek.
\newblock Multicategories revisited.
\newblock In {\em Categories in computer science and logic ({B}oulder, {CO},
  1987)}, volume~92 of {\em Contemp. Math.}, pages 217--239. Amer. Math. Soc.,
  Providence, RI, 1989.

\bibitem[Law63]{lawvere:functsem}
F.~William Lawvere.
\newblock Functorial semantics of algebraic theories.
\newblock {\em Proc. Nat. Acad. Sci. U.S.A.}, 50:869--872, 1963.

\bibitem[Law02]{LawvereMetric}
F.~W. Lawvere.
\newblock Metric spaces, generalized logic and closed categories.
\newblock {\em Reprints in Theory and Applications of Categories}, 1:1--37,
  2002.

\bibitem[Lei02]{GeneralizedEnrichmentofCategories}
Tom Leinster.
\newblock Generalized enrichment of categories.
\newblock {\em J. Pure and Applied Algebra}, 168:391--406, 2002.

\bibitem[Lei04]{leinster:higher-opds}
Tom Leinster.
\newblock {\em Higher operads, higher categories}, volume 298 of {\em London
  Mathematical Society Lecture Note Series}.
\newblock Cambridge University Press, Cambridge, 2004.

\bibitem[Low88]{lowen:approach}
R.~Lowen.
\newblock {\em Approach Spaces: The Missing Link in the
  Topology-Uniformity-Metric Triad}.
\newblock Oxford Science Publications, 1988.

\bibitem[LS02]{ls:ftm2}
Stephen Lack and Ross Street.
\newblock The formal theory of monads. {II}.
\newblock {\em J. Pure Appl. Algebra}, 175(1-3):243--265, 2002.
\newblock Special volume celebrating the 70th birthday of Professor Max Kelly.

\bibitem[Man09]{manzyuk:closed}
Oleksandr Manzyuk.
\newblock Closed categories vs. closed multicategories.
\newblock arXiv:0904.3137, 2009.

\bibitem[May72]{may:goils}
J.~Peter May.
\newblock {\em The geometry of iterated loop spaces}.
\newblock Springer-Verlag, Berlin, 1972.
\newblock Lectures Notes in Mathematics, Vol. 271.

\bibitem[Mor06]{morton:catalg-qm}
Jeffrey Morton.
\newblock Categorified algebra and quantum mechanics.
\newblock {\em Theory Appl. Categ.}, 16:No. 29, 785--854 (electronic), 2006.

\bibitem[Pen09]{penon:tcats-rep}
Jacques Penon.
\newblock T-cat\'egories repr\'esentables.
\newblock {\em Theory Appl. Categ.}, 22(15):376--387, 2009.

\bibitem[Pow99]{power:enr-theories}
John Power.
\newblock Enriched {L}awvere theories.
\newblock {\em Theory Appl. Categ.}, 6:83--93 (electronic), 1999.
\newblock The Lambek Festschrift.

\bibitem[Sea05]{SealCanonical}
Gavin~J. Seal.
\newblock Canonical and op-canonical lax algebras.
\newblock {\em Theory and Applications of Categories}, 14:221--243, 2005.

\bibitem[Shu08]{FramedBicats}
Michael Shulman.
\newblock Framed bicategories and monoidal fibrations.
\newblock {\em Theory and Applications of Categories}, 20:650--738, 2008.

\bibitem[Shu10]{shulman:vequip}
Michael~A.\ Shulman.
\newblock Limits and completions in formal category theory.
\newblock In preparation, 2010.

\bibitem[SS08]{ss:extn-laxalg}
Christoph Schubert and Gavin~J. Seal.
\newblock Extensions in the theory of lax algebras.
\newblock {\em Theory Appl. Categ.}, 21:No. 7, 118--151, 2008.

\bibitem[Str72]{street:formalmonads}
Ross Street.
\newblock The formal theory of monads.
\newblock {\em Journal of Pure and Applied Algebra}, 2:149--168, 1972.

\bibitem[Str74a]{street:elem-cosmoi-i}
Ross Street.
\newblock Elementary cosmoi. {I}.
\newblock In {\em Category Seminar (Proc. Sem., Sydney, 1972/1973)}, pages
  134--180. Lecture Notes in Math., Vol. 420. Springer, Berlin, 1974.

\bibitem[Str74b]{street:fib-yoneda-2cat}
Ross Street.
\newblock Fibrations and {Y}oneda's lemma in a {$2$}-category.
\newblock In {\em Category Seminar (Proc. Sem., Sydney, 1972/1973)}, pages
  104--133. Lecture Notes in Math., Vol. 420. Springer, Berlin, 1974.

\bibitem[Str76]{street:catval-2lim}
Ross Street.
\newblock Limits indexed by category-valued {$2$}-functors.
\newblock {\em J. Pure Appl. Algebra}, 8(2):149--181, 1976.

\bibitem[Str80]{street:fibi}
Ross Street.
\newblock Fibrations in bicategories.
\newblock {\em Cahiers Topologie G\'eom. Diff\'erentielle}, 21(2):111--160,
  1980.

\bibitem[Str83]{street:enr-cohom}
Ross Street.
\newblock Enriched categories and cohomology.
\newblock In {\em Proceedings of the Symposium on Categorical Algebra and
  Topology (Cape Town, 1981)}, volume~6, pages 265--283, 1983.
\newblock Reprinted as Repr. Theory Appl. Categ. 14:1--18, 2005.

\bibitem[Tho07]{tholen:coursenotes}
Walter Tholen.
\newblock Lax-algebraic methods in general topology.
\newblock {\em Summer School on Contemporary Categorical Methods in Algebra and
  Topology (course notes)}, 2007.

\bibitem[Tho09]{tholen:ordered-topstruct}
Walter Tholen.
\newblock Ordered topological structures.
\newblock {\em Topology Appl.}, 156(12):2148--2157, 2009.

\bibitem[Ver92]{DomThesis}
Dominic Verity.
\newblock Enriched categories, internal categories, and change of base.
\newblock {\em Cambridge University (PhD Thesis)}, 1992.

\bibitem[Web05]{weber:operads}
Mark Weber.
\newblock Operads within monoidal pseudo algebras.
\newblock {\em Appl. Categ. Structures}, 13(5-6):389--420, 2005.

\bibitem[Web07]{weber:2toposes}
Mark Weber.
\newblock Yoneda structures from 2-toposes.
\newblock {\em Appl. Categ. Structures}, 15(3):259--323, 2007.

\bibitem[Woo82]{wood:proarrows-i}
R.~J. Wood.
\newblock Abstract proarrows. {I}.
\newblock {\em Cahiers Topologie G\'eom. Diff\'erentielle}, 23(3):279--290,
  1982.

\end{thebibliography}

\end{document}